\documentclass[11pt]{amsart}
\usepackage{amssymb, amsmath, amsthm,ulem}
\usepackage[latin1]{inputenc}
\usepackage{graphicx}
\usepackage[hidelinks]{hyperref}
\usepackage{color}
\usepackage{ulem}
\usepackage{epstopdf}
\usepackage{cancel}
\usepackage{tikz}

\usepackage[a4paper, centering]{geometry}
\usepackage{color}
\usepackage{graphicx}
\usepackage{scalerel,stackengine}

\usepackage{mathtools}

\newtheorem{theorem}{Theorem}
\newtheorem{proposition}[theorem]{Proposition}

\newtheorem{lemma}[theorem]{Lemma}
\newtheorem{corollary}[theorem]{Corollary}

\theoremstyle{definition}
\newtheorem{definition}[theorem]{Definition}
\theoremstyle{remark}
\newtheorem{remark}[theorem]{Remark}

\definecolor{verde}{RGB}{20,150,100}
\definecolor{purple}{RGB}{200,30,200}

\newcommand{\EEE}{\color{black}}

\stackMath
\newcommand\reallywidecheck[1]{%
\savestack{\tmpbox}{\stretchto{%
  \scaleto{%
    \scalerel*[\widthof{\ensuremath{#1}}]{\kern-.6pt\bigwedge\kern-.6pt}%
    {\rule[-\textheight/2]{1ex}{\textheight}}
  }{\textheight}%
}{0.5ex}}%
\stackon[1pt]{#1}{\scalebox{-1}{\tmpbox}}%
}

\def\R{\mathbb{R}}

\def\N{\mathbb{N}}

\def\A{\mathcal{A}}

\newcommand{\sm}{\setminus}
\newcommand{\vps}{\varepsilon}
\newcommand{\Om}{\Omega}
\newcommand{\om}{\omega}

\newcommand{\sq}{\subseteq}
\newcommand{\ra}{\rightarrow}

\def \d{\delta}
\def \e{\varepsilon}

\newcommand{\abs}[1]{{\left|#1\right|}}
\newcommand{\norma}[1]{{\left\Vert#1\right\Vert}}
\newcommand{\vphi}{\varphi}

\begin{document}
\title[]{ 
  The geometric size of the fundamental gap  
}

\bigskip\bigskip

\author[]{Vincenzo Amato, Dorin Bucur, Ilaria Fragal\`a}

\thanks{}

\address[Vincenzo Amato]{
	Holder of a research grant from Istituto Nazionale di Alta Matematica "Francesco Severi" at Dipartimento di Matematica e Applicazioni "R. Caccioppoli", Via Cintia, Complesso Universitario Monte S. Angelo, 80126 Napoli, Italy.}
\email{amato@altamatematica.it}

\address[Dorin Bucur]{
Universit\'e  Savoie Mont Blanc, Laboratoire de Math\'ematiques CNRS UMR 5127 \\
  Campus Scientifique \\
73376 Le-Bourget-Du-Lac (France)
}
\email{dorin.bucur@univ-savoie.fr}

\address[Ilaria Fragal\`a]{
Dipartimento di Matematica \\ Politecnico  di Milano \\
Piazza Leonardo da Vinci, 32 \\
20133 Milano (Italy)
}
\email{ilaria.fragala@polimi.it}

\date{\today}

\begin{abstract}

  The fundamental gap conjecture  proved by Andrews and Clutterbuck in 2011 provides
 the sharp lower bound for the fundamental gap of the Dirichlet Laplacian
 on any convex set in $\R^N$
  in terms of the diameter. We strengthen this seminal result by proving that  the excess 
of the fundamental  gap compared to the diameter,   can be quantified in terms of flatness. In particular, this answers affirmatively an open question  about rigidity raised by Yau in 1990.     The proof  is of variational nature and takes a different route from 
Andrews and Clutterbuck parabolic approach. We exploit a  one dimensional reduction issued from  collapsing cells of a convex partition,
the advantage being that we can decrypt separately the local contributions to the gap coming from  the eigenfunction first and second order gradient,  
and from the geometric degeneration. After 
quantifying the excess of the gap in each of these one dimensional collapsed cells, 
a  thorough geometric    analysis of their
collective behaviour in the partition  leads to the quantitative inequality   in two and  higher dimensions. As a by-product of our 
method we  prove,
 in a stronger version,  a conjecture  from 2007 by Hang-Wang on the  quantitative form of
 the classical  Payne-Weinberger  inequality.

\end{abstract} 

\maketitle 

\tableofcontents

\section{Introduction and statement of the results}
 Given an open bounded domain $\Om \subset \R ^N$,   the difference $\lambda _2 (\Om)- \lambda _ 1 (\Om)$ between its second and first Dirichlet Laplacian  eigenvalues  is usually referred to as  the
fundamental gap of $\Om$.  
It has several important implications in different areas of 
both mathematics and physics, e.g.\ heat diffusion, 
statistical mechanics, quantum field theory, numerical analysis.
Finding sharp lower bounds for the fundamental gap is a problem whose history covers several decades, so that we summarize it without any attempt of completeness. 
In the pioneering work \cite{vdB83}, van  den Berg   first  observed that, for many convex domains $\Om$, the gap is bounded from below 
by $\frac{ 3 \pi ^ 2}{ D_\Om ^ 2}$,   where  $D_\Om$  is the diameter of $\Om$. The validity of such inequality for any convex domain $\Om$
 was then    conjectured 
by  Yau \cite{Y85}  and by Ashbaugh-Benguria   \cite{AB89}, in the more general case of a Schr\"odinger operator of the form $- \Delta + V$, 
being $V$ a convex  potential on $\Om$. 
This more general formulation of the problem is meaningful also in the one-dimensional case, which was solved by Lavine \cite{Lav94}
see also \cite{AB89, H03}. 
A breakthrough in higher dimensions is due to Singer-Wong-Yau-Yau \cite{SWYY}, who obtained the lower bound $\frac{ \pi ^ 2}{ 4 D _\Om ^ 2}$, 
later improved into $\frac{\pi ^ 2}{ D _\Om ^ 2}$ by Yu-Zhong \cite{YZ86} and Smits \cite{smits}. These lately  non-optimal lower bounds rely
on the earlier fundamental result by Brascamp-Lieb \cite{BL76}, which states that the first Dirichlet eigenfunction is log-concave on any convex domain
(different proofs were given by Korevaar \cite{K83} and Singer-Wong-Yau-Yau \cite{SWYY}). Let us also mention the lower bound  $\frac{(\log 2 ) ^ 2}{D_\Om^ 2}$ 
 obtained by Bobkov \cite[inequality (2.8)]{B99}  when the Lebesgue measure in the Rayleigh quotient is replaced by any absolutely continuous measure
with log-concave density. In the early 2000s, the expected optimal lower bound $\frac{3 \pi ^ 2}{D_\Om ^ 2}$ 
has been obtained in some particular cases  when $\Om$ satisfies specific geometric assumptions 
\cite{BK01, BM00, davis}. 
An excellent overview up to that date is the paper by Ashbaugh \cite{ash}, where more related references can be found.  Further advances based on upper bounds  for $\nabla ^ 2 \log u _ 1$  were given in \cite{Y03,  ling}. 

The conjecture  was finally  proved  in 2011 by Andrews-Clutterbuck in  \cite{AC11}
(see  also the survey paper \cite{Areview} and the subsequent related works \cite{carron, DSW, lurowlett}).  
The groundbreaking new idea by Andrews-Clutterbuck  is    the   following refinement of
Brascamp-Lieb result into an improved log-concavity inequality for the first Dirichlet eigenfunction
\begin{equation}\label{f:improved}\big ( \nabla \log u _ 1 ( y) - \nabla \log u _ 1 ( x) \big ) \cdot \frac{ y - x } { \| y - x \| } \leq - 2 \frac{\pi}{D_\Om} \tan \Big (  \frac{\pi}{D_\Om} 
\frac{\| y - x \| }{2}  \Big ) \quad \forall x, y \in \Om \,.\end{equation}
 This  estimate is obtained by a parabolic approach and,
 combined with a method to control  the modulus of continuity of solutions to parabolic equations, allows them
to   prove  the conjectured lower bound.  
 Afterwards, Ni 
  recovered the sharp control of the gap   by  
 an elliptic argument, which still exploits the improved-log-concavity estimate \eqref{f:improved}. 
  
 However,   though the fundamental gap   conjecture was solved,
neither the parabolic nor the elliptic proofs could answer  the rigidity question, which   was 
formulated in  1990  by Yau himself as problem no.\ 44 in his ``Open problems in geometry" paper \cite{YAMS}: 
 
 \smallskip
 \centerline{ {\it Is the gap inequality always {\it strict} in dimension $N \ge 2$?} }
 \smallskip
 
 In case of 
an affirmative answer,  the   next  challenging question is of quantitative nature:  

\smallskip
\centerline{{\it  Is it possible to  evaluate  the excess of the gap in terms of the flatness of the convex set?}}
\smallskip

 The main difficulty arises from the lack of an optimal domain, which places the problem far away from 
quantitative spectral inequalities solvable via the use of second order shape derivatives.  

In this paper we answer   the above   questions. Our main result reads:

\begin{theorem}\label{t:quantitative} Let $N \geq 2$. 
There exists a  dimensional constant  $\overline c>0$ such that, for every open bounded convex domain $\Omega$ in $\R ^N$ with diameter $D_\Om$ and width $w_\Om$, we have
\begin{equation}\label{f:qgap}
\lambda _ 2 (\Omega)- \lambda _ 1 (\Om) \geq \frac{3 \pi ^ 2}{D_\Om^ 2}  + \overline c \frac{w_\Om^6}{D_\Om^ 8} \,.
\end{equation} 
\end{theorem}

  Our approach to obtain Theorem \ref{t:quantitative} is  variational. It stems from the  key  idea that  the fundamental gap can be   controlled by the collective behaviour of a  family of one dimensional problems.  
 In this process,  the role played by the improved log-concavity of the first Dirichlet eigenfunction
is totally different with respect to Andrews-Clutterbuck proof, and is intimately related  with the geometry of the domain, eventually 
leading us to encompass the previous information on  the size of the fundamental gap.

As a starting point, swe adopt the perspective originally due to Thompson and Kac \cite{TK69}
of viewing the fundamental gap as a first weighted Neumann eigenvalue  (see also \cite{KS87, S84, smits}).  Setting, for any positive log-concave  weight 
$p \in L ^ 1 (\Om)$, 
   \begin{equation}\label{def:mu}
\mu _ 1 (\Om, p ):= \inf
 \Big \{ { \frac{   \int _ \Om |\nabla u | ^ 2  p \, dx }{ \int_\om u^ 2 p \, dx }   } \ :\ u  \in H ^ 1_{\rm loc} (\Om) \cap L ^ 2 (\Om, p \,  dx)  \,, \ \int _\om u p \, dx = 0 \Big \}\,,
\end{equation} 
the Dirichlet fundamental gap can be recast by choosing the specific weight $p = u _ 1 ^ 2$: 
denoting by $u _1, u _2$ the first two Dirichlet eigenfunctions  in of $\Om$,  it holds
\begin{equation*}
 \lambda _ 2 (\Omega) - \lambda _ 1 (\Omega) = \mu _ 1 (\Om, u_ 1 ^ 2) , \qquad \text{ with eigenfunction } \overline u := \frac{u _ 2}{ u _ 1}\,.
\end{equation*}

The synergy between the above equality and the improved log concavity  of the weight $u _1 ^2$  has  been guessed in the literature as a possible way to attack the fundamental gap inequality: such a relationship, firstly addressed by Smits \cite{smits}, was disclosed by Andrews himself in 
\cite{An12}, see also \cite{LR12} for a statement going in this direction.
Specifically, as soon as $p$ satisfies  the improved log-concavity condition
\begin{equation}\label{f:improved2}\big ( \nabla \log p( y) -\nabla \log p ( x) \big ) \cdot \frac{ y - x } { \| y - x \| } \leq - 
 \frac{\pi}{D_\Om} \tan \Big (  \frac{\pi}{D_\Om} 
\frac{\| y - x \| }{2}  \Big ) \quad \forall x, y \in \Om \,,\end{equation}
the classical Payne-Weinberger inequality  (see \cite{PW,beb, ENT, FNT12})
   \begin{equation}\label{abf41.1}
\mu_1(\Om, p) \ge \frac{\pi^2}{D_\Om^2}
\end{equation}
  can be replaced by the stronger one 
 \begin{equation}\label{abf41}
\mu_1(\Om, p) \ge \frac{3\pi^2}{D_\Om^2}.
\end{equation}
Nevertheless,  though    the above  inequality can be  derived by the techniques of Andrews and Clutterbuck, 
  using their way the inequality is achieved  asymptotically,  with few room to make it quantitative, nor to get its rigidity.

 Thus we adopt a different approach.  The key new idea  is to exploit,  in place of parabolic pdes,  
the global behaviour of   an interacting family of one dimensional eigenvalue problems. These are  Schr\"odinger eigenvalue 
  problems with measure potential, set  on the 
 collapsed cells of a (suitably modified)
 convex partition \`a la Payne-Weinberger. 
The cells carry new weights which involve, besides the    squared   improved log-concave eigenfunction $u_1^2$, a new geometric term issued from the collapsing procedure, which is power-concave by Brunn-Minkowski Theorem. 
From these fine concavity properties, the measure potential inherits a particular structure. We take advantage of this structure in order  to decrypt separately the local contribution to the fundamental gap coming from different quantities, specifically
the first and second order gradient of $u _ 1$, and  the collapsed geometry.  This requires a major analytic work in one space dimension, 
which eventually equips us with a one-dimensional refined estimate. With this estimate in our hands,
  we are able to
localize and hence strengthen  the inequality \eqref{abf41} as follows:  we prove that there exists a universal constant $C$ such that, 
for every (possibly low dimensional) convex set $\om \subseteq \Om$ of diameter $D _\om$,  it holds 
\begin{equation}\label{f:diffd.1} \mu _ 1 (\om, p ) \geq \frac{3 \pi ^ 2}{D _\Om ^2}  + C \frac{(D _\Om -D_\om ) ^ 3}{D_\Om ^ 5} \,,
\end{equation} 
 see Proposition \ref{p:inside}. 
 The above inequality leads quite directly to answer affirmatively Yau rigidity question. 
    It also serves as a picklock to obtain the refined inequality  \eqref{f:qgap}. 
To that aim, we need to implement a more geometric point of view:   roughly, we have to 
         to get an insight into the geometry of the partition  in order to 
estimate the excess of the gap   over the diameter     generated by the  combined assessment  of the cells.  Here the local control of the second Dirichlet eigenfunction plays a crucial role. 
 
 \smallskip 
 
 We believe that the method of proof  used to obtain Theorem \ref{t:quantitative} can be fruitfully employed also for studying the geometric size of the fundamental gap on Riemannian manifolds. We stick to the Euclidean space to enlighten the ideas in the simplest setting. On the other hand, 
 we wish to present   a   related outcome of the approach used 
  to prove Theorem \ref{t:quantitative}:
taking   $p \equiv 1$ in    \eqref{f:diffd.1},    
   leads  to the rigidity and to   a sharp    quantitative form of  the   classical    Payne-Weinberger  inequality    for the first nontrivial
Neumann eigenvalue $\mu _ 1 (\Om)$.
  Since  $\mu _ 1 (\Om)$ can also be seen as the Neumann fundamental gap  (see  \cite{ash}),  
  there is a clear analogy  with 
the Dirichlet fundamental gap. 
  In the Neumann case, 
 the    saturation   question 
was asked by Sakai \cite{sakai} 
(and settled for  smooth compact Riemannian manifolds  with nonnegative Ricci curvature \cite{hangwang, Valtorta}), 
  while the quantitative question has been formulated by Hang-Wang  in 2007 \cite{hangwang}, along with the conjecture
  of a lower bound of the type
$\frac{\pi^2}{D_\Om ^ 2 } + \overline c  \frac{w _\Om ^ 2}{D_\Om^4}$, being 
$w_\Om$ the width of $\Om$.  
Our result  
encompasses Hang-Wang conjecture,  as  it shows that their expected lower bound holds  in any space dimension with  the second largest John semi-axis in place of the width. 
Recall that, up to a translation and rotation, for any convex domain $\Om \subset \R ^N$ there exists an ellipsoid $\mathcal E  = \{ \sum _{i=1} ^N \frac{x_ i ^ 2 }{a_i ^ 2 } < 1 \}$, called {\it John ellipsoid}, such that $\mathcal E \subseteq \Om \subseteq N \mathcal E$ (see e.g. \cite{guzman}).  

In the same vein of Theorem  \ref{t:quantitative},  we obtain Theorem \ref{t:quantitativeN} below
(whose proof is simpler  because we have to replace the weight $u_1^2$ by  a constant).

\begin{theorem}\label{t:quantitativeN}
Let  $N\ge 2$. There exists a dimensional constant $\overline c>0$ such that, 
 for every open bounded convex domain $\Omega$ in $\R ^N$ with diameter $D_\Om$
and John ellipsoid of semi-axes  $a_1 \geq \dots \geq a_N$, we have
	\begin{equation}
		\label{quantitative3d}
		\mu_1(\Omega)\ge \frac{\pi^2}{  D_\Om ^ 2} + \overline c \frac{a_2^2}{ D _\Om ^4}.
	\end{equation}
\end{theorem}

\begin{remark}   
In our proofs of Theorems \ref{t:quantitative} and \ref{t:quantitativeN} there is no evident loss of sharpness at any step.  
This leads to the power $6$ for the width in \eqref{f:qgap} and  to the power $2$ for the second dimension of the John ellipsoid   in \eqref{quantitative3d}. 
In the latter case, the power $2$ is optimal: taking $\Omega_\vps = (0,d)\times (0,\vps)^{N-1}$,  the second John  axis  of $\Omega_\vps$ equals $\e$, and we have
$$\mu_1(\Omega_\vps)= \frac{\pi^2}{d^2}= \frac{\pi^2}{D_{\Omega_\vps}^2}+ \frac{\pi^2}{D_{\Omega_\vps}^4} (N-1)\vps^2 + o(\vps^2)\,.$$
On the other hand, in the former case,  we can neither prove, nor disprove,  the optimality of the power $6$. 
In the Dirichlet case, the loss of sharpness {\it if true}, might be related to a possibly suboptimal knowledge of the geometry of the first Dirichlet eigenfunction near the boundary and of its localization.
More specifically, 
 the reason why the width is  the  only geometric size  of a convex set $\Om$ that we are able to relate  with the excess of its Dirichlet fundamental gap,   seems to   be the fact that  
the John ellipsoid  does not encode the dimensions of 
 the relevant parallelepiped contained in $\Om$ which is known to carry over 
most of the information about the spectrum   (as it is  well described in dimension $2$ in \cite{Je95,GJ98,GJ09,beck}).   
\end{remark}

\begin{remark}     An explicit estimate of the constant $\overline c$ appearing in \eqref{quantitative3d}, 
without any attempt  of optimality,
 might be rather easily given just by tracking it in all steps of the proof. 
A similar  target   for the constant $\overline c$ in \eqref{f:qgap} seems to be more delicate. 
 \end{remark} 

\begin{remark} 
We point out that Theorem \ref{t:quantitative} does not hold unaltered 
for the Schr\"odinger equation with a convex potential $V$ 
$$u \in H^1_0(\Om), \qquad-\Delta u +V u = \lambda_k(\Om,V) u \quad \mbox{ in } \Om \,.$$
Actually, while the gap inequality $\lambda_2(\Om,V)-\lambda_1(\Om, V) \ge  \frac{3 \pi ^ 2}{D_\Om^ 2}  $
is still true \cite{AC11}, its  quantitative form \eqref{f:qgap}   cannot hold keeping the same positive constant $\overline c$ independent of  $\Omega$ and $V$.  
This can be easily seen by looking at the following example in dimension $N = 2$. 
Let $\Om= \{(x,y)\in  \R^2 : |x|+|y|<1\}$ and $V_{\vps, \delta}(x,y) = \frac 1\delta \big (|y|-\vps\big )^+$.  
  We write the inequality \eqref{f:qgap}  and  we first 
pass to the limit as $\delta \ra 0^+$ 
at fixed $\vps $.   Since 
$\lambda_k(\Om, V_{\vps, \delta}) \ra \lambda_k(\Om _\e)$, where  $\Om _\e:= 
\{\Om \cap \{(x,y) : |y|<\vps\}$, we obtain $ \lambda_2(\Om _\e) -\lambda_1(\Om _\e)  \geq  \frac{3 \pi ^ 2}{4}  +  \frac{\overline c}{32}$.  Then, 
by  using the monotonicty of the eigenvalues with respect to inclusions and passing to the limit as $\vps \ra 0^+$, we obtain 
$$\frac{3 \pi ^ 2}{4}  +  \frac{\overline c}{32} 
 \le \lambda_2(\Om _\e) -\lambda_1(\Om _\e)  \le  \frac{4\pi^2}{ ( 2-2\vps)^2} + \frac{\pi^2}{4\vps^2}- \Big( \frac{\pi^2}{ 4} + \frac{\pi^2}{4\vps^2}   \Big) \to \frac{3\pi^2}{4} \,,$$
which leads to  $\overline c=0$.
\end{remark}

{\bf Outline of the proof.} We give hereafter a short overview of  the proof of Theorem \ref{t:quantitative}. 
 We refer to Section \ref{sec:Neumann}  for the specific modifications required for the proof of Theorem \ref{t:quantitativeN},  which include
  some nontrivial ones, in particular a  geometrically explicit  
$L ^ \infty$ estimate for Neumann eigenfunctions,  see Proposition \ref{ambu25}.  

 As a starting point,    having in mind that we aim at an inequality such as \eqref{abf41},    \EEE 
 we decompose   $\Om$ \`a la Payne-Weinberger, namely as the union of 
$n$ mutually disjoint convex cells  of equal measure, obtained by successively ``cutting'' $\Omega$ by hyperplanes parallel to a fixed direction,
on which the eigenfunction $\overline u$ has zero integral mean with respect to the measure $u _ 1 ^ 2 \, dx$. 
In the limit as $n \to + \infty$,   this operation allows 
to estimate 
 from below $\mu _ 1 (\Om, u _ 1 ^ 2)$  in terms a  one dimensional eigenvalue of the type
$\mu _ 1 (I, p)$, where  $I$ is a line segment contained into $\Om$, 
and  $p= h u _ 1 ^ 2$, being $h$ the $\mathcal H ^ {N-1}$ measure of the cell's section orthogonal to $I$.  

 Then the proof is developed 
along    following  lines.

\medskip

{\it I. A sharp $1D$ lower bound.}
 The  first step  is  to obtain a lower bound for  $\mu _1 ( I, p)$,
  with exhibits an explicit dependence on the structure of the weight $p$.  
The first basic feature, coming  from the log-concavity of both $u _ 1$ and $h$, is that
$p$  itself  is log-concave. 
This yields that $m _p: =   \frac{3}{4} \big (\frac{p'}{p} \big ) ^ 2 - \frac{1}{2} \frac{p''}{p}$ 
is a positive measure and that the inequality $\mu _ 1 (I, p ) \geq \lambda _ 1 (I, m_p)$ holds (see \cite{PW}),  being, for any positive measure $q$, 
$$\lambda _ 1 (I, q ):= \inf \Big \{ \frac{  \int _ {I}  |u' | ^ 2   dx  + \int_{I } |u| ^ 2  \, dq  }{\int _{ I} u  ^ 2  \, dx } \ :
\ u  \in H ^ 1_0 (I)\cap L ^ 2 (I, q)  \Big \}\,. $$
But the more subtle  key feature coming from the  factor $u _ 1 ^ 2$ in the weight $p$, 
is that  $p$ itself satisfies  the  improved log-concavity estimate \eqref{f:improved}. Hence the function $\psi_p := -  \frac{1}{2} (\log p  ) ' $
belongs  to the following class of functions, which are thought for convenience  as functions with extended real values defined on the fixed interval
$I _ \pi = ( - \frac \pi 2, \frac \pi 2 ) $:
$$\mathcal A ( I _\pi):=  \Big \{ \psi \text{ increasing},\ \psi (y) - \psi ( x) \geq 2 \tan \Big ( \frac {y-x}{2}\Big  ) \  \text{ if }  x < y  \text{ in } {\rm dom} (\psi)  \Big \}\,.
$$
Since  the measure $m _p$ can be written as 
$m _p = 
\psi _p ' + \psi  _p ^ 2 $, we arrive at the following novel $1D$ problem, which encodes  in the class of competitors the {\it log-concavity modulus}  of $u _ 1$: 
\begin{equation}\label{abf41.2}
\inf \Big \{ \lambda _ 1  (I _\pi, \psi' + \psi ^2 ) \ :\ \psi \in \mathcal A ( I _\pi) \Big \}\,.
\end{equation}  
Surprisingly, the above infimum can not only be estimated, but exactly computed:  in Theorem \ref{t:1d} we prove that it equals  $3$, 
and it is attained uniquely at the function  $\overline \psi (x) =  \tan (x)$  
(which corresponds to the weight $\overline p = \phi _ 1 ^ 2$, being $\phi _ 1 ( x) = \cos x$ the first Dirichlet eigenfunction of $I _\pi$).
A noticeable feature  is that the optimal function $\overline \psi =  \tan (x)$ does not saturate pointwise the equality sign in the definition of  $
\mathcal A ( I _\pi)$. 
 The key of the proof is a new sophisticated  ad-hoc   procedure of {\it stratified rearrangement}
(see Definifion \ref{d:strat}) which allows to handle the 
modulus of concavity constraint imposed on the 
 functions in the admissible class $
\mathcal A ( I _\pi)$.

The value $3$ given by Theorem \ref{t:1d}
 is clearly the good one in order to recover Andrews-Clutterbuck gap inequality for a convex  domain $\Om$, taken for convenience  
 of diameter $\pi$.
   In order to refine their inequality and proceed with the estimate of the gap excess,  
in Theorem \ref{t:1dextra} we establish two distinct refinements of Theorem \ref{t:1d},   holding 
  when the weight $p = h u _ 1 ^ 2$ enjoys some additional properties.
     
The first refinement has the target of handling cells of ``small'' diameter: for these cells the additional property 
of the weight is that the 
function $\psi =  -  \frac{1}{2} (\log p  ) '  $  in $\mathcal A (I _\pi)$  has finiteness domain of length $d < \pi$.
This leads to an improved  lower bound of the following type,  for an absolute constant  $C$:
\begin{equation}\label{f:ref1}  \lambda _ 1 (I_\pi, \psi ' + \psi ^ 2 ) \geq 3+ C ( \pi -d ) ^ 3 \,.
 \end{equation}
Here the power $3$ is obtained via a perturbation argument. Within the class $\mathcal A (I _\pi)$   the power $3$ is  sharp:
this is precisely the point leading to the power $6$ of the width in  Theorem \ref{t:quantitative}.  
 
The second refinement has the target of handling cells of ``large'' diameter.  A careful  analysis of the polygonal structure of the partition  will reveal that, in $2D$,  it is enough to analyse    only such cells for which the height  $h$ in orthogonal direction to a  diameter is an affine function
away from the endpoints.  
Then, denoting by $h _{min}$ and $h _{max}$  the extrema of the affine function $h$,  we obtain    an improved lower bound of the following type, for an absolute constant $K$:
 \begin{equation}\label{f:ref2}
 \lambda _1 ( I _ d,  \psi ' + \psi ^ 2 ) \geq 3 + K  \Big ( 1 - \frac{h _{min}} {h _{max}} \Big ) ^2\,.
 \end{equation}

\bigskip
{\it II.  Localized  version  and rigidity of Andrews-Clutterbuck inequality.}  
Exploiting the one-dimensional estimate \eqref{f:ref1}, we   prove   the above mentioned  
 inequality \eqref{abf41}    and the 
new  localized version   of Andrews-Clutterbuck  inequality  \eqref{f:diffd.1}.

  Enforcing  \eqref{f:diffd.1},  
 in Theorem \ref{t:rigidity} we derive the rigidity of Andrews-Clutterbuck inequality. 
The idea is the following: 
if by contradiction
 $\lambda _ 2 (\Om) - \lambda _ 1 (\Om) =  \frac{3 \pi^ 2}{D ^ 2}$, taking 
a partition of $\Om$ into $n$ mutually disjoint convex sets $\Om _i$  having the mean value property $ \mu _ 1 (\Om, u _ 1 ^ 2) \geq \frac{1}{n} \sum _{i = 1} ^n \mu _ 1 (\Om_i, u _ 1 ^ 2 )$,  
we get 
 \begin{equation*}\frac{3 \pi ^ 2}{D_\Om^ 2}  = \mu _ 1 (\Om, u _ 1 ^ 2) \geq \frac{1}{n} \sum _{i = 1} ^n \mu _ 1 (\Om_i, u _ 1 ^ 2 ) \geq \frac{3 \pi ^ 2}{ D_\Om^2}  + C  \sum _{i = 1} ^n
\frac{(D_\Om-D_{\Om_i}) ^ 3}{ D_\Om ^ 5}  \,.\end{equation*}
This implies that all the diameters  $ D_{\Om_ i}$  are equal to $ D_\Om$, which yields a contradiction: 
for $N= 2$, the contradiction comes from a geometric argument,  because the equality of all the diameters 
forces $\Om$ to  be a circular sector, for $N \geq 3$
 the same argument applies  to a suitable two-dimensional section of $\Om$.
  
We stress that, in order to gain the  above mentioned mean value property of the cells, we need 
to work with a new kind of partitions,    which are   
distinct  from the classical ones by  Payne-Weinberger not only for the presence of the weight $u _1 ^ 2$, but also
for the equipartition request: 
the  measure equipartition condition  $|\Om _ i| = \frac{1}{n} |\Om|$  is replaced by the  
$L^2$-equipartition condition $ \int_{ \Om _i} \overline u ^ 2 u _ 1 ^ 2     = 
  \int_{ \Om _i} u _ 2 ^ 2  = \frac{1}{n}$.

\medskip
 {\it III. }  {\it The collective behavior of the convex partition yielding the quantitative inequality}.  Evaluating the excess of the gap demands to get an insight about the geometric display of the cells
of a partition of $\Om$, meant as a  weighted  $L^2$ 
equipartition of $\Om$  which enjoy the  mean value property quoted above. Using this type of partitions
offers the advantage that  the  required estimate of the excess is fulfilled   as soon as, for a fixed proportion of cells, the eigenvalue $\mu _1 (\Om _i, u _ 1 ^ 2)$ is sufficiently large with respect to $\frac{3\pi^2}{D_\Om^2}$, with a controlled increment.
 So we are led to assess the geometry of the cells.

For $N = 2$, we  distinguish a list of binary crossroads in cascade, depending on different geometric
features holding for  a fixed proportion of  cells. 
The main distinction is made by looking at the size of the cell's  diameter: if  most of the cells
have ``small'' diameter
(in the sense that the difference between the diameter of $\Om$ and the diameter of the cell
is controlled from below by the width), 
applying  to such cells the localized inequality  \eqref{f:diffd.1}  we get the  estimate of the excess. 
Otherwise, if 
most of the cells have ``large'' diameter,  
assuming that the quantitative inequality does not occur, a contradiction is obtained through a 
geometric argument  which 
can be intuitively  sketched as follows.  
Since we are dealing with the situation in which
most cells are thin and long, 
they can be vertically piled over a diameter of $\Om$, and they have a profile function
which is affine
away from the endpoints (if this was not the case, the quantitative inequality would hold as well, 
by analyzing the position of vertices in the partition and the consequent presence of other cells with small diameter). Then,  the
one-dimensional refined inequality \eqref{f:ref2} applies to  such cells, 
i.e., to the one dimensional problem set on  their diameter and, if the extra term is small, 
we get the geometric information that 
the cells have to be ``almost" rectangular (the precise meaning is given in Section \ref{s5.1}).
  At this point,   
 we obtain a uniform control on the height of each cell, related to the non localization of the second eigenfunction (see Remark \ref{r:propsection}). 
 This leads to  the conclusion that the pile of the cells is, in a sense, ``too high'', 
 because the actual diameter would be strictly larger than $D_\Om$,  finally yielding a contradiction.

For $N \geq 3$, the result is obtained by a partial slicing procedure, reducing ourselves 
to a two-dimensional analysis
involving a modified weight,   which can be carried over by similar arguments as the ones used to treat the case $N = 2$.

\medskip 
The paper is organized as follows. Section \ref{sec:1D} is devoted to the analysis of the 1D-eigenvalue problem associated with the measure potentials issued from   restrictions of first eigenfunctions to line segments contained into a convex set. In Section \ref{sec:partitions} 
we introduce the modified Payne-Weinberger partitions.    In Section \ref {sec:localized} 
we prove the  localized strengthened version of  Andrews-Clutterbuck inequality and  we establish the rigidity  of the gap inequality.    Section \ref{sec:construction} contains the proof of Theorem \ref{t:quantitative} while in Section \ref{sec:Neumann} we prove Theorem \ref{t:quantitativeN}.  
  In the Appendix we collect some useful results about eigenfunctions associated with  weighted Neumann eigenvalues.

\section{The sharp one-dimensional lower bound}\label{sec:1D}

Let $I$ be  a one-dimensional  open bounded interval.
Given  a positive weight $p$  in $L ^ 1 ( I )$ and a nonnegative   Borel measure $q$  possibly taking the value $+\infty$, consider the following {weighted} Neumann eigenvalue and Dirichlet eigenvalue  with {potential}: 
\begin{eqnarray*}& \displaystyle \mu _ 1 (I, p ):= \inf \Big \{ \frac{  \int _ {I} | v' | ^ 2  p \, dx  }{\int _ {I} v  ^ 2 p \, dx } \ :\ v  \in H ^ 1 _{\rm loc} (I)\cap L ^ 2 (I, p \, dx)  \, , \ \int _{I} v p = 0 \Big \}   & \\  \noalign{\medskip} 
&  \displaystyle \lambda _ 1 (I, q ):= \inf \Big \{ \frac{  \int _ {I}  |v' | ^ 2   dx  + \int_{I } |v| ^ 2  \, dq  }{\int _{ I} v  ^ 2  \, dx } \ :
\ v \in H ^ 1_0 (I)\cap L ^ 2 (I, q)  \Big \} \,.   
\end{eqnarray*}   
Let us mention that eigenvalues associated with potentials which are measures, such as $\lambda _ 1 (I , q)$,  have been extensively studied in the context of shape optimization in dimension $N\ge 2$, see  for instance
\cite[Section 4.3]{bubu}.  

 If $p$ is log-concave, 
we can   introduce   the positive measure $m _p$ defined  by 
 \begin{equation}\label{f:measuremp}
 m _p
: =  \Big [ \frac{3}{4} \Big (\frac{p'}{p} \Big ) ^ 2 - \frac{1}{2} \frac{p''}{p} \Big ] = \psi _p ' + \psi  _p ^ 2 \, , \qquad \text{ with } \psi_p := -  (\log p ^ {\frac{1}{2}} ) '\, .
\end{equation}

Here $\psi' _p=  
  [\frac{1}{2} \Big (\frac{p'}{p} \Big ) ^ 2 - \frac{1}{2} \frac{p''}{p}]$ is the distributional derivative of the non-decreasing function $\psi _p$, which is a nonnegative measure thanks to the log-concavity of $p$, while $\psi _ p ^ 2$ denotes with a slight abuse of notation the nonnegative measure $\psi_p ^ 2 \, dx$. 

For simplicity,  also in the sequel  we denote measures which are absolutely continuous 
simply by writing their density
with respect to the Lebesgue measure.

 \begin{lemma}\label{l:ml} \ 
For any  positive log-concave   weight $p \in L ^ 1 ( I )$, 
if $m _p$ is given by \eqref{f:measuremp}
 it holds 
 $$\mu _ 1 (I, p) \geq \lambda _ 1 ( I, m _ p)\,.$$
    \end{lemma}

\proof     Under the additional assumptions $p \in W^{1, \infty}(I)$ and $\inf _{x \in I}p(x) >0$, an eigenfunction  $v$ for $ \mu _ 1 (I_\pi, p )$  exists in $H^2(I)$ and satisfies $$
\begin{cases}
- ( pv' ) ' =\mu _ 1 (I, p )   pv  &  \text{ in }  I 
\\ \noalign{\medskip} 
pv' (- \frac{\pi}{2}  ) = pv' ( \frac{\pi}{2}  ) = 0 \,. & 
\end{cases}
$$ 
Then the function $w := p ^ { 1/2}v'$ belongs to $H^1_0(I)$ and solves
$$
\begin{cases}
- w'' + \Big [ \frac{3}{4} \Big (\frac{p'}{p} \Big ) ^ 2 - \frac{1}{2} \frac{p''}{p} \Big ]  w = \mu _ 1 (I_\pi, p )  w & \text{ in } I 
\\ \noalign{\medskip} 
w(- \frac{\pi}{2}   ) = w ( \frac{\pi}{2}  ) = 0 \,, & 
\end{cases}
$$ 
yielding the inequality $\mu _ 1 (I, p) \geq \lambda _ 1 ( I, m _ p)$. 

 Assume now $p \in L ^ 1 ( I )$ is positive and log-concave. 
Letting $I^{\vps}$  be  intervals compactly included in $I$ and increasingly   converging to  $I$,   we have   $\inf _{x \in I^\e}p(x) >0$ and $p|_{I^\e}\in W^{1, \infty} (I)$. Consequently,
the Neumann eigenvalue problem $\mu _ 1 ( I ^\e, p)$ is well-posed  and  the inequality
$\mu _ 1 (I^\e, p) \geq \lambda _ 1 ( I^\e, m _ p)$ is satisfied. 
Then we obtain that the same inequality holds true for the interval $I$   by observing that 
$$\lambda_1(I, m_p)= \lim _{\vps \ra 0} \lambda_1(I^{\vps}, m_{p})\quad \mbox{ and } \quad  \mu _ 1 (I, p)  \ge \limsup _{\vps \ra 0} \mu _1 (I^{\vps}, p).$$
Indeed, the first  assertion follows from the inclusion $I^\e \sq I$ and the fact that $I^\e$ is increasingly converging to $I$. 
The second  assertion follows from the monotone convergence theorem, 
since, for  any admissible test function for  $\mu_1(I,p)$, its restriction to
$I^\vps$, corrected by a small constant so to make it orthogonal to $p$ in $L ^ 2 (I ^ \e)$, becomes an admissible test function for 
$\mu_1(I^ \e,p)$.
 \qed

\bigskip

By Lemma \ref{l:ml},  we are led
to deal with Dirichlet eigenvalues of the 
type $\lambda _ 1 (I, m_p)$, with $p = h u _ 1 ^ 2$ and $m _p$ given by \eqref{f:measuremp}.
The corresponding function $\psi _p$ has the form
$$\psi _p =- (\log u _ 1)'  - \frac{1}{2}  (\log h)' \,.$$

The heart of the matter is that, by the improved log-concavity estimate \eqref{f:improved}, 
 the function $\psi _p$  turns out to belong to  the  class of 
 functions  defined  hereafter in \eqref{f:defclass}. 
 To formulate the problem, without any loss of generality we work on the interval 
$$I _ \pi := \big ( - \frac {\pi}{ 2} , \frac{ \pi}{ 2 } \big )\,,$$ and 
we  introduce the following class of functions  defined on $I _\pi$ 
with values in $\overline \R = \R \cup \{ \pm \infty\}$ and finiteness domain  $ {\rm dom} (\psi)$
\begin{equation}\label{f:defclass} \mathcal A ( I _\pi):=  \Big \{ \psi \text{ increasing},\ \psi (y) - \psi ( x) \geq 2 \tan \Big ( \frac {y-x}{2}\Big  ) \  \text{ if }  x < y  \text{ in } {\rm dom} (\psi)  \Big \}\,.
\end{equation}

For functions  $\psi \in \mathcal A (I _\pi)$, we tacitly extend the measures  $\psi ^ 2 $ 
and $\psi'$ to $+ \infty$ in $I _\pi \setminus {\rm dom} (\psi)$.  
Then our target can  be precisely expressed as the study of the minimization problem 
\begin{equation}\label{f:minpb}
\min \Big \{ \lambda _ 1  (I _\pi, q) \ :\ q = \psi' + \psi ^2  \text{ for some } \psi \in \mathcal A ( I _\pi) \Big \}\,,
\end{equation} 

\smallskip
The remaining of this section is entirely devoted to that goal. 
It is divided in two parts: 

\smallskip
-- in the first part  we give the sharp lower bound  for $\lambda _1  (I _\pi, \psi' + \psi ^ 2)$ for $\psi \in \mathcal A (I _\pi)$, namely we fully solve the minimisation problem \eqref{f:minpb}, see Theorem \ref{t:1d}; 

\smallskip
-- in the second part,  for $\psi \in \mathcal A (I _\pi)$  with ${\rm dom}(\psi) = ( - \frac{d}{2}, \frac{d}{2} )= :I_d$,   being $d<\pi$,   we  give some lower bounds for $\lambda _1  (I _d, \psi' + \psi ^ 2)$ with extra terms involving the    difference  $(\pi-d)$,   see Theorem \ref{t:1dextra}.
Here and in the sequel, for $d < \pi$, if $m$ is a nonnegative Borel measure on $I_d$, we identify the eigenvalue $\lambda _ 1 (I_d, m)$ with 
$\lambda _ 1 (I _\pi, \tilde m)$,  where $\tilde m$  is the measure obtained extending $m$ to $+ \infty$ on $I_\pi \setminus I_d $.

\smallskip

Let us start with an elementary observation:

\begin{remark}\label{r:tan} The function $\overline \psi  ( x) = \tan x$ belongs to $\mathcal A ( I _\pi)$. The corresponding weight $\overline q(x) = 1 + 2 \tan ^ 2 \! x$ is equal to $m_{\overline p}$ for 
$\overline p = \phi _ 1 ^ 2$, being $\phi _ 1 ( x) = \cos x$ the first Dirichlet eigenvalue on $I _\pi$, and  we have the following equalities for eigenvalues, 
all of them with eigenfunction $ \cos ^ 2 \!x$:
 $$\begin{array}{ll}
& \lambda _ 1 (I_\pi , 2 \overline \psi   ^ 2 ) = \lambda _ 1 (I_\pi , 2\tan ^ 2 \! x) = 2 \\  \noalign{\medskip} 
& \lambda _ 1 (I_\pi , 2  \overline \psi ' )  = \lambda _ 1 (I _\pi, 2 (1+ \tan ^ 2 \! x)) = 4 \\  \noalign{\medskip} 
 & \lambda _ 1 (I _\pi, \overline \psi'+ \overline \psi ^ 2 )  = \lambda _ 1 (I_\pi, 1 + 2 \tan ^ 2 \!x) = 3 \,.
 \end{array}$$ 
\end{remark}

It is somehow natural to wonder whether the weight $\overline q$ associated with the one-dimensional eigenfunction  $\phi _1$ 
is optimal for the minimization problem \eqref{f:minpb}. Our result below states that this exactly is the case. 
In addition, $\overline q$ is the {\it unique} solution.

\begin{theorem}\label{t:1d} 
Let  $q = \psi' + \psi ^ 2$, with  $\psi \in \mathcal A ( I _\pi)$. Then 
\begin{equation}  \label{f:stima3}
 \lambda _ 1  (I _\pi, q) \geq 3 \,,
\end{equation}
with equality if and only if $q(x) = 1 + 2 \tan ^ 2 \!x$ (and in this case an  eigenfunction is $\cos ^ 2 \!x$). 
\end{theorem}

\medskip 
The proof of Theorem \ref{t:1d} is built upon the following two independent propositions. 

\begin{proposition}\label{p:prima} 
For every $\psi \in \mathcal A ( I _\pi)$, it holds  
\begin{equation}\label{f:stima2}
\lambda _ 1 (I_\pi,  2 \psi ^ 2) \geq 2 \,, 
\end{equation} 
with equality if  and only if $\psi ( x) = \tan x$ (and in this case an  eigenfunction is $\cos ^ 2 \!x$).  
\end{proposition}

\begin{proposition}\label{p:seconda}
For every $\psi \in \mathcal A ( I _\pi)$, it holds  
\begin{equation}\label{f:stima4} \lambda _ 1 (I_\pi,  2 \psi ') \geq 4\,.
\end{equation}  
In particular, equality occurs if $\psi ( x) = \tan x$  (and in this case an  eigenfunction is $\cos ^ 2 \!x$).  
\end{proposition}

Let us  first show how the above propositions imply Theorem \ref{t:1d}, and then turn back to their proof. 

\bigskip {\it Proof of Theorem \ref{t:1d}}. It is easy to check 
that 
the map $q \mapsto \lambda _ 1(I_\pi,  q)$ is concave. Indeed, for every pair of weights $q_1, q_2$, every $t \in [0, 1]$, and every $v \in H ^ 1 _0 ( I_\pi)$ we have
$$ \begin{array}{ll} \displaystyle \frac{\int _{I _\pi} (v') ^ 2 +  v ^ 2 \, ((1-t)d q_1 + t d q _ 2 )  }{\int _{I_\pi} u ^ 2 }   
& \displaystyle =  ( 1- t) \frac{\int _{I _\pi} (v') ^ 2 +  v ^ 2 \, d q_1  }{\int _{I_\pi} v^ 2 } 
 + t  \frac{\int _{I _\pi} (v') ^ 2 +  v ^ 2 \, d q_2   }{\int _{I_\pi} v ^ 2 } \\ \noalign{\medskip} & \geq (1-t) \lambda _ 1 (I _\pi, q _1) + t \lambda_1 ( I_\pi, q_2) \,.
 \end{array}$$

Therefore, for every $\psi \in \mathcal A ( I _\pi)$,  we have
 \begin{equation}\label{f:separation}
 \lambda _ 1  (I_\pi, \psi' + \psi ^ 2 ) \geq \frac{\lambda _ 1 (I _\pi, 2\psi') + \lambda  _ 1 (I_\pi,  2 \psi ^ 2)}{2} \,.
 \end{equation}
 Then  the inequality \eqref{f:stima3}   follows immediately from Proposition \ref{p:prima} and Proposition \ref{p:seconda}.

Concerning the equality case,  when $q ( x) = 1 + 2 \tan ^ 2 \!x$, we  have 
$\lambda _ 1 (I _\pi, q) = 3$, with eigenfunction $\cos^2 \!x$ (cf.\ Remark \ref{r:tan}).
Viceversa, if $\lambda _ 1 (I _\pi, q) = 3$, for some $q = \psi ' + \psi ^ 2$ with $\psi  \in \mathcal A (I _\pi)$, it follows from the inequalities \eqref{f:stima2}, \eqref{f:stima4}, and \eqref{f:separation},  that $\lambda _ 1 (I_\pi, 2 \psi') = 4$ and $\lambda _ 1 (I_\pi, 2 \psi ^ 2) = 2$.  
 By the last assertion in Proposition \ref{p:prima}, we conclude that $\psi ( x) = \tan x$, and hence $q ( x) = 1 + 2 \tan ^ 2 x$. 
\qed

\bigskip
For later use, we state below a rigidity result for  the Neumann eigenvalue $\mu _ 1 (I _\pi, p)$, which 
 is  a straightforward by-product of Theorem \ref{t:1d}. 
  
 \begin{corollary}\label{c:sharpmu}
If $\mu _ 1 (I_\pi, p) = 3$ for a positive log-concave weight $p$ such that $\psi _p \in \mathcal A ( I _\pi)$, we have $\psi _ p (x)= \tan x$, and $p (x) = k \cos ^ 2 \!x$ for some   positive constant $k$. 
\end{corollary}

 \proof   If $\mu _ 1 (I_\pi, p) = 3$, by Lemma \ref{l:ml} and Proposition \ref{p:seconda}, we have
$$3 = \mu _ 1 (I_\pi, p) \geq \lambda _ 1 (I_\pi, \psi ' _p + \psi_p ^ 2 ) \geq \frac12 \big ( \lambda _ 1 (I_\pi, 2 \psi ' _p) +\lambda _ 1 ( I_\pi, 2 \psi_p ^ 2 ) \big )    \geq 
 2 + \frac 12 \lambda _ 1 ( I_\pi, 2 \psi_p ^ 2 ) \big )  
\geq 3 \,.$$
We infer that $\lambda _ 1 ( I_\pi, 2 \psi_p ^ 2 )  = 2$. Therefore, by Proposition \ref{p:prima} we conclude that 
 $\psi _p (x) = \tan x$, and hence $p ( x) = k \cos ^ 2 \! x$ for some positive constant $k$. \qed  

\bigskip
We now provide the proofs of Propositions \ref{p:prima} and \ref{p:seconda}.

\bigskip
{\it Proof of Proposition \ref{p:prima}}. 
We first prove the following claim: for every function $\psi \in \mathcal A (I _\pi)$, there exists another function $\widetilde \psi \in \mathcal A (I _\pi)$  which 
changes sign once in $I _\pi$, and satisfies 
$\lambda _ 1 (I _\pi, \psi ^ 2) \geq \lambda _  1 ( I _\pi , \widetilde \psi ^ 2) $. We search for  $\widetilde \psi$ in the subclass of 
$\mathcal A (I _\pi)$ given by functions of the type  $ \psi + c$, for $c \in \R$.  If $v$ denotes an eigenfunction for $\lambda _ 1 (I_\pi, \psi ^ 2)$, normalized in $L ^ 2$, we have 
$$\lambda _ 1 (I _\pi, \psi ^ 2 ) = 
  \int _ {I_\pi}  |v' | ^ 2   dx  + \int_{I _\pi}  \psi ^ 2 |v| ^ 2  \, dx   
\geq \inf _ {c \in \R } 
\int _ {I_\pi}  |v' | ^ 2   dx  + \int_{I _\pi}  (\psi + c) ^ 2 |v| ^ 2  \, dx  \,.
$$
By differentiating with respect to $c$, we see that the above infimum is attained at 
$\tilde c = - {\int _{I _\pi} \psi | v| ^ 2 }$. 
The function $\widetilde \psi:= \psi + \tilde c$ satisfies $\int  _{I _\pi} \widetilde \psi |v | ^ 2 = 0$, so that
it changes sign at least one time,  and exactly one time, because $\widetilde \psi \in \mathcal A (I _\pi)$, so that it is strictly increasing. Moreover,  
$$\lambda _ 1 ( I _\pi, \psi ^ 2) \geq 
\int _ {I_\pi}  |v' | ^ 2   dx  + \int_{I _\pi}  (\psi + \tilde c) ^ 2 |v| ^ 2  \, dx  
 \geq \lambda _ 1 (I _\pi, \widetilde \psi)\,.$$ 

In view of the claim just proved,   it is not restrictive to prove the inequality \eqref{f:stima2} under the assumption 
that the function $\psi$ changes sign exactly one time in $I _\pi$.

Consider the function $|\psi|$.  Thanks to the assumption that $\psi$ has exactly one zero $x _0 \in I _\pi$, we know that the function $|\psi|$  vanishes at $x_0$, is strictly decreasing for $x \leq x _0$, and strictly increasing for $x \geq x _0$. So, for almost every $t >0$, the level set $\{ |\psi| < t \}$ is an interval (containing $x_0$). 
We  rearrange the function $|\psi|$ into the even function $|\psi| _*$ defined by
 $$\{ |\psi|_*< t \} :=    \{|\psi| < t \}  ^* \qquad \forall t \in (0, \|\psi \| _\infty]\,,$$
where $  \{ |\psi| < t \} ^*$   denotes the translation of the interval  $\{ |\psi| < t \}$ which sends its midpoint to the origin.
(Notice that this  is a kind of ``symmetric increasing rearrangement'', which  is the analogue of the classical  symmetric decreasing rearrangement, just replacing super-levels with sub-levels). 
 By construction, for every $x \in \big  (0, \frac{{|\rm dom}  (\psi)  |}{2} \big )  $,   we have 
 $$|\psi| _* (x) =| \psi|( b_x) =  |\psi| ( a_x)  \qquad \text{  } x = \frac{b_x-a_x}{2},  \text{ with } a_x<x_0 <b_x \text{ and }  \psi (a_x) = - \psi ( b_x)  \,.$$

Thus the assumption $\psi \in \mathcal A ( I _\pi)$  be expressed as a pointwise inequality for $|\psi| _ *$:
\begin{equation}\label{f:eqri} 
|\psi| _* (x)  = \frac{1}{2} \big(  \psi ( b_x) - \psi ( a_x) \big )  \geq \tan \Big (\frac{b_x - a _x}{2}  \Big )  = \tan x \qquad \forall  x \in  \big (0, \frac{{|\rm dom}  (\psi)  |}{2} \big ) \,.
\end{equation} 

Let now $v \in H ^ 1 _0 ( I _\pi)$. Denoting by $v ^*$ its classical symmetric decreasing rearrangement,  
and by  ${\rm dom} (\psi) ^*$ the interval of length $| {\rm dom} (\psi) | $ centred at the origin, 
it holds
\begin{equation}\label{f:PS}\int _{{\rm dom} (\psi)} |v'| ^ 2 \geq \int _{{\rm dom} (\psi)^*} |(v^*)' |  ^ 2 \qquad \text{ and } \qquad
 \int _{{\rm dom} (\psi)}  \psi ^ 2 |v | ^ 2   \geq \int _{{\rm dom} (\psi) ^*}   |\psi |_*  ^ 2 |v^* | ^ 2\,.
 \end{equation}
Indeed,  
the first inequality is the classical P\'olya-Szeg\"o inequality for the decreasing rearrangement (see e.g.\ \cite[Section II.4]{K}), while the second one 
follows from the classical Hardy-Littlewood inequality (see e.g.\ \cite[Chapter 10]{HLP}), the sign of the inequality being reversed because for one of the two involved functions  ($|\psi|$)  we take the  increasing rearrangement $|\psi |_*$ defined above in place of $|\psi|^*$.  Then, we have 
\begin{equation}\label{f:globalall}
 \frac{\int _{{\rm dom} (\psi)} \!\! |v'| ^ 2 \!+ \! 2 \psi ^ 2 v ^ 2  }{\int _{{\rm dom} (\psi)} \!  |v| ^ 2 }\!  \geq  \! \frac{\int _{{\rm dom}  (\psi) ^*} \!\!  |(v^*)' |  ^ 2 \! +\!  2 |\psi |_*  ^ 2 |v^* | ^ 2  }{\int _{{\rm dom} (\psi)^*}  \!\! | v^* | ^ 2 }  \!  \geq  \!
 \frac{\int _{I _\pi}  \! |(v^*)' |  ^ 2 \!+ \! 2\tan^2  \!x   |v^* | ^ 2  }{\int _{I _\pi} \! | v^* | ^ 2 }  \!  \geq 2  \,, 
\end{equation}
where in the second and third inequality we have used respectively the estimate \eqref{f:eqri} and Remark \ref{r:tan}.

Concerning the equality case, if $\psi (x)  =  \tan  x$, we have  $\lambda _ 1 (I_\pi, 2 \psi^2) = 2$, with eigenfunction $ \cos ^ 2 \! x$  (cf.\ Remark \ref{r:tan}). 
Viceversa, assume that the equality $\lambda _ 1 (I_\pi,  2 \psi ^ 2) = 2$ holds for some function $\psi$. 
Then,  if $v$ is an eigenfunction for
$\lambda _ 1 (I_\pi,  2 \psi ^ 2)$, all the inequalities in \eqref{f:PS} and \eqref{f:globalall} must hold with equality sign. 
From the fact  that the last inequality in \eqref{f:globalall} holds with equality sign, 
we infer that $v ^* = k \cos ^ 2 \, x$.

In particular, $v ^*$ does not have critical level sets of positive Lebesgue measure. 
Then the fact that the first inequality in \eqref{f:PS} holds with equality sign implies that
$v = v ^*$ \cite[Theorem 1.1]{BZ} (see also \cite[Theorem 4.1]{Hajaiej}).     

In turn, since the second inequality in \eqref{f:globalall} holds with equality sign,  we infer that $ |\psi|_* ^ 2
= \tan ^ 2 \!x $ a.e. By the monotonicity of $\psi$, we have $\psi (x )= \tan  x$ for every $x \in \R$.   
\qed

\bigskip

We now turn to the proof of Proposition \ref{p:seconda}, which requires some preliminaries. 
Given a function $\psi \in \mathcal A ( I _\pi)$, 
in order to estimate from below   $\lambda _ 1 ( I _\pi, 2 \psi')$, 
we need to find lower bounds for integrals of the following type, for $v \in H ^ 1 _ 0 ( I _\pi)$:  
 $$\int _{I _\pi} v ^ 2 \, d {\nu _{\psi'}}  \,. $$
 Here and in the sequel,  we use the notation 
 $d {\nu _{\psi'}}$ when writing an integral with respect to the measure $\psi'$.   
  We have: 
\begin{equation}\label{r:area} 
\begin{array}{ll} \displaystyle \int _{I _\pi} v ^ 2 \, d \nu _{\psi'}  & \displaystyle = \int _{I _\pi} \int _0 ^ {+ \infty} \chi _{\{v^ 2 > s \} } \, ds \,  d \nu_{\psi'} = 
\int _0 ^ {+ \infty}\!\!\! \int _{I _\pi}  \chi _{\{v ^ 2 > s \} } \,  d \nu_{\psi'}\, ds   \\ \noalign{\medskip} & \displaystyle= \int_0 ^ {+ \infty} \!\!\! \nu _{\psi'}  \big ( \{v  >  \sqrt s \}   \big )  \, ds 
= \int _0 ^ {+ \infty} \!\!\!
\nu _{\psi'} \big ( \{v  > t \}   \big ) \, 2t \, dt \,. \end{array} 
\end{equation} 
    
On the other hand,  the constraint  $\psi \in \mathcal A ( I _\pi)$ can be equivalently expressed as   the following inequality holding for all intervals $[a, b]\subset I _ \pi  \cap {\rm dom}  (\psi)$: 
\begin{equation}
\label{r:constraint} 
\nu _{\psi'} ( [a, b])    \geq   2 \tan \big ( \frac{b-a}{2} \big )= \int _ a ^ b \Big ( 1 + \tan ^ 2  \big ( x-\frac{a+b}{2} \Big  ) \big ) \, dx   \,.\end{equation}

Condition \eqref{r:constraint} cannot be applied directly to estimate the integral \eqref{r:area}, because 
not all level sets of $u$ are intervals. Thus we are going to exploit condition \eqref{r:constraint}  in a more subtle way, 
passing through the introduction of new notions of {\it stratified rearrangement} and {\it stratified potential}; they  are obtained by finitely many applications of elementary constructions
that we call respectively   {\it blocked rearrangement}
and {\it blocked potential}.  

Below we introduce these definitions on a generic bounded open interval $I$, 
for  $v$ in   the space $\mathcal X (I)$ of functions which are continuous on the closure of $I$, attain their global  minimum   
at both the endpoints of $I$,  and have finitely many local minima in $I$, each one attained at a single point.  
(Notice that $H ^ 1 _0 ( I ) \cap \mathcal X (I)$  is a dense subspace of $H ^ 1 _0 ( I)$).

We  denote by $[v ^*, I] $ the symmetric decreasing rearrangement of $v$ with respect to  the mid-point $x_I$ of $I$.

 \begin{definition}[blocked rearrangement] \label{d:block} Let  $v \in \mathcal X(I)$. 
  The blocked rearrangement of $v$ in $I$ is the function $[v^ {\flat}, I]$ defined on $I$ as follows: 
 \begin{itemize}
 \item[(i)]  If $v$ does not have any local minimum in $I$,  $$ [v ^{\flat}, I] := [v ^*, I] \,.$$

\smallskip
  \item[(ii)] If $v$ has some local minimum in $I$, letting $\ell$  be
  the smallest level of local minimum (so that $\{ v > \ell \}$ is the union of two consecutive open intervals),
 and denoting respectively by  $x_I$ and $x _\ell$ the midpoints of $I$ and  of the closed interval $\overline{\{ v > \ell \}}$ , 
  $$[v ^ { \flat}, I] (x):=  \begin{cases}
   [v ^* , I] (x) & \text{ if  }    [v ^* , I] (x) \leq \ell 
   \\
   v ( x- x_\ell + x_I ) & \text{ otherwise } 
   \end{cases} 
   $$ 
   \end{itemize}
  \end{definition}



\begin{definition}[stratified rearrangement] \label{d:strat} 
Let $v \in \mathcal X( I )$.    
Its stratified rearrangement  is the function $\tilde v$ defined on $I$  via
a finite number of blocked rearrangements as follows:
\begin{itemize}
\item[--] 
Set $I ^1 := I $, and   
$v _1:=[v ^{ \flat},  I ^1]$. 
Two cases may occur:

\smallskip
(i) If $v$ does not have any local minimum in $ I ^1$, define $\tilde v:= v _1$ in the whole interval $ I ^1$, and the definition stops here. 

\smallskip
(ii) If $v$ has some local minimum in $ I ^1$, letting $\ell_1$ be the smallest value of local minimum for $v$  in $I^1$, define $\tilde v := v _1$ only on the set $[v ^ {*},  I ^1]  \leq \ell _ 1$, and for the definition on its complement in $I ^ 1$ go to the next step.

\smallskip
\item[--]  Set $I ^{1, j}$ the two consecutive open intervals such that
$\big \{  v _1 > \ell _1 \big \}  = I  ^{ 1,1}  \cup  I  ^{1,2}$,  and 
$v _{1, j}:= [ v _ { 1} ^{\flat} ,  I ^{1,j}]$, for $j = 1, 2$. 
For a fixed $j \in \{1, 2 \}$,
 two cases may occur:

\smallskip
 
 (i) If $v _ { 1}$  does not have any local minimum in $I^{1,j}$,  define 
  $\tilde v:= v _{1, j}$ in the whole interval $ I ^{1, j}$, and the definition on the interval $I ^{1, j}$ stops here. 

\smallskip

(ii) If $v_{1}$ has some local minimum in $ I^{1, j}$, letting $\ell_2$ be the smallest value of local minimum for $v_{1}$  in $I^{1, j}$, define $\tilde v := v _{1, j}$ only on the set $[v _1^ {*},  I ^{1, j}]  \leq \ell _ 2$, and for  the definition on its complement in $I ^ {1, j}$ go to the next step.  
 \smallskip 
 \item[--]  Set $I ^{1, j, 1}$ and  $I ^{1, j, 2}$ the two consecutive open intervals such that
 $\big \{  v _{1, j} > \ell _2 \big \}  = I  ^{ 1,j,1}  \cup  I  ^{1,j,2}$, and  proceed  as in the previous steps. After finitely many steps, 
the procedure stops because by assumption the number of local minima is finite. 
\end{itemize}
 \end{definition}

\begin{remark}\label{r:multiindex}   
The open intervals constructed in Definition \ref{d:strat} can be labeled by a family $\mathcal F$ of  multi-indices 
 $\alpha= ( \alpha _1, \dots, \alpha _k)$, with  $\alpha _1 = 1$ and $\alpha _ i \in \{1, 2 \}$  for $i \geq 2$, and  $k$ less than or equal to the number of  levels of  local minimum of $v$ in $\overline I$.   We denote by $ \Gamma $ the subfamily of $\mathcal F$ of multi-indices $\alpha$ such that $I ^ \alpha$ contains a local minimum of $\tilde v$. 
 Equivalently, we set 
 \begin{equation}\label{f:gamma}
 \Gamma := \Big \{  \alpha  \in \mathcal F \text{ such that }  (\alpha, 1) \text{ and } (\alpha, 2) \text{ belong to  } \mathcal F \Big \}\,.
 \end{equation} 
Then by construction the stratified rearrangement $\tilde v$ enjoys the following symmetry property with respect to the mid-point $x _\alpha$ of each interval $I ^ \alpha$:
 if $\alpha\not \in \Gamma$,    
 $\widetilde v$  is symmetric with respect to   $x _\alpha$ on the whole interval $I ^ \alpha$; 
 if $\alpha \in \Gamma $,  $\widetilde v$  is symmetric with respect to   $x _\alpha$ 
 just on 
 $I ^\alpha\setminus (   I ^{\alpha, 1} \cup   I ^{\alpha, 2} )$. 
  \end{remark}

\begin{definition}[blocked potential]
Let  $v \in \mathcal X (I)$,  and assume $|I| \leq \pi$.  The blocked potential of $u$ on $I$ is defined as the function
 $1 + \tan ^ 2  ( x - x_I )$,
 with  definition   domain equal to the subset of $I$ where $ [v ^{\flat}, I] = [v ^*, I]$.
 \end{definition}

\begin{definition} [stratified potential] Let  $v \in \mathcal X (I)$,  and assume $|I| \leq \pi$.
Let $I ^\alpha$ be the family of open intervals constructed in Definition \ref{d:strat},    labeled as in 
Remark \ref{r:multiindex}.  The stratified  potential of $v$ is the  function $V$ defined on $I$ by glueing all the blocked potentials of $v$ on the intervals $I ^ \alpha$.  Equivalently, 
for any $x \in I$, pick the longest  multiindex $\alpha$ such that 
$x \in I^ \alpha$, and  set 
$$V ( x) = 1 + \tan ^ 2 \Big ( \frac{  x - x _\alpha}{2} \Big )\,, \qquad x _\alpha:= \text{midpoint of  } I _\alpha. 
 $$
\end{definition}

\begin{remark}  With the notation introduced in Remark \ref{r:multiindex},  the stratified potential $V$ can be identified as follows:   
 if $\alpha \not \in \Gamma $,    then $V ( x) = 1 + \tan^2 ( x- x_\alpha)$ 
on the whole interval $I ^ \alpha$; 
 if $\alpha \in \Gamma $,   then  $V ( x) = 1 + \tan^2 ( x- x_\alpha)$ 
 just on 
 $I ^\alpha\setminus (   I ^{\alpha, 1} \cup   I ^{\alpha, 2} )$. 
  \end{remark}

The next two lemmas, relying on the definitions of stratified rearrangement and potentials, provide the intermediate results needed for
the proof of Proposition \ref{p:seconda}.

\begin{lemma}\label{l:lemmatilde}  Let $\psi \in \mathcal A ( I _\pi)$ and  let $v\in  H ^ 1 _0 (I _\pi) \cap \mathcal X( I _\pi)$. Denoting by  
$\tilde v$ and $V$ respectively the stratified rearrangement and potential associated with $v$, it holds 
\begin{equation}\label{f:bbelow0} 
\int_{I_\pi}  v ^ 2 (x)  \, d\nu _ {\psi'}   \geq \int _{I_\pi} V ( x) \tilde v ^ 2 (x)    \, dx  \,.
\end{equation}

\end{lemma} 
\proof   If ${\rm dom}( \psi)$ is strictly contained into $I _\pi$ and $v  \not \in H ^ 1 _ 0  ({\rm dom } (\psi))$, then 
the l.h.s.\ of  \eqref{f:bbelow0} is $+ \infty$, and the inequality is trivially true. So we only have to consider the situation when  
$v \in H ^ 1 _ 0  ({\rm dom } (\psi))$. 
Denoting by $\nu _{\psi'}$ and $\nu _V$ the absolutely continuous measures with densities $\psi'$ and $V$, 
and recalling \eqref{r:area}, we have  
 $$ \int_{I_\pi}   v ^ 2 (x)  \, d\nu _ {\psi'}     - \int _{I_\pi} V ( x) \tilde v ^ 2 (x)    \, dx  =  \int_0 ^{\| v \| _\infty   }  2t \Big  [ \nu _{\psi'} 
 \big (  \{ v  > t \} \big )  - \nu _ V (  \{ \tilde v  > t \} \big ) \Big  ]  \, dt  $$ 
 Let us distinguish two cases, according to whether the family of multi-indices $\Gamma$ introduced in \eqref{f:gamma} is empty or not. 
 
When  $\Gamma  = \emptyset$ (or equivalently,  $v$ does not have any local minimum in $I _\pi$), the set $\{ v>t \}$ is an interval for every $t\in (0, \| u \| _\infty)$. 
 Then \eqref{f:bbelow0} is satisfied because, by \eqref{r:constraint}, it holds
$$\nu _{\psi'} \big ( \{ v>t \}   \big )  \geq 
 2 \tan \Big ( \frac {   | \{ v>t \} | }{2} \Big ) =   2 \tan \Big ( \frac {   | \{ \tilde v>t \}  | }{2} \Big )   = \nu _V\big ( \{ \tilde v>t \}   \big ) \,  . $$ 

When  $\Gamma \neq \emptyset$, 
not all level sets $\{ v > t \}$ are intervals. Therefore,  in order to exploit the estimate \eqref{r:constraint}, 
we decompose the set $\{ v > t \}$  (and accordingly $\{ \tilde v  > t \} $) as a finite union of disjoint intervals $J _n (t)$
(and of their translations $\widetilde J _n (t)$)
$$ \{ v  > t \}  = \bigcup _{n= 1  } ^ {N(t)}   J _n (t) \,, \qquad    \{ \tilde v  > t \}  = \bigcup _{n= 1  } ^ {N(t)}    \widetilde J _n (t) 
\,.$$   
  
We focus attention on a fixed interval $\widetilde J _n (t)$.
Since we  are   working under the assumption that $\Gamma \neq \emptyset$, 
  we have that 
$\widetilde J _n ( t)$
is  contained into some interval $ I ^\alpha$ with $\alpha  \in \Gamma$. 
Among these intervals $I ^\alpha$ containing $\widetilde J _n ( t)$, we choose the one  with multiindex $\overline \alpha = \overline \alpha (n, t)$ 
having the maximum number of components.  

By applying the estimate \eqref{r:constraint} to the interval $J _ n ( t)$, we obtain
\begin{equation}\label{f:spsi}   \nu _{\psi'} \big ( J _n (t)   \big )  \geq 
 2 \tan \Big ( \frac {   | J _n ( t) | }{2} \Big ) =   2 \tan \Big ( \frac {   |  \widetilde J _n ( t) | }{2} \Big )   , \end{equation} 

We observe that
$$ \begin{array}{ll} &  \displaystyle 2 \tan \Big ( \frac {   |  \widetilde J _n ( t) | }{2} \Big )   =  2 \tan \Big ( \frac {   |  \widetilde J _n ( t) | }{2} \Big )    - 
 2 \tan \Big ( \frac {|  I    ^{\overline \alpha, 1} \cup  I ^{\overline \alpha, 2} | }{2} \Big ) +  
2 \tan \Big ( \frac { |  I  ^{\overline \alpha, 1} \cup  I ^{\overline \alpha, 2} | }{2} \Big )  
\\  \noalign{\medskip} 
& \displaystyle = 2 \tan \Big ( \frac {   |  \widetilde J _n ( t) | }{2} \Big )    - 
 2 \tan \Big ( \frac {|  I    ^{\overline \alpha, 1} \cup I ^{\overline \alpha, 2} | }{2} \Big ) +  
2 \tan \Big ( \frac { |  I  ^{\overline \alpha, 1} | }{2} \Big )    + 2 \tan \Big ( \frac { |  I  ^{\overline \alpha, 2} | }{2} \Big )   + \mathcal R _{\overline \alpha} \,,
\end{array} 
$$ 
where 
$$\mathcal R _{\overline \alpha} := 2 \tan \Big ( \frac { | I  ^{\overline \alpha, 1} \cup I ^{\overline \alpha, 2} | }{2} \Big )  - 2 \tan \Big ( \frac { |  I  ^{\overline \alpha, 1} | }{2} \Big )    - 2 \tan \Big ( \frac { |  I  ^{\overline \alpha, 2} | }{2} \Big ) \geq 0 \,,$$ 
the last inequality being due to the super-additivity of the tangent function on $(0, \frac{\pi}{2})$.
 
In case the multiindices $(\overline \alpha, 1)$ and $(\overline \alpha, 2 )$ do not belong to $\Gamma$, we have 

\begin{equation}\label{f:sU} 
\begin{array}{ll}  
& \displaystyle \nu _{V} \big ( \widetilde J _n (t)   \big )    = 
\nu _{V} \big ( \widetilde J _n (t)  \setminus  ( I ^{\overline \alpha, 1}  \cup I ^{\overline \alpha, 2} ) \big )   + 
\nu _ V ( I ^{\overline \alpha, 1}   ) + \nu _ V ( I ^{\overline \alpha, 2}   ) \\ 
\noalign{\medskip} 
& \displaystyle = 
  2  \Big [ \tan \Big ( \frac {   |  \widetilde J _n ( t) | }{2} \Big )    - 
  \tan \Big ( \frac {|  I    ^{\overline \alpha, 1} \cup   I ^{\overline \alpha, 2} | }{2} \Big )  \Big ] +  
2 \tan \Big ( \frac { |  I  ^{\overline \alpha, 1} | }{2} \Big )    + 2 \tan \Big ( \frac { |  I  ^{\overline \alpha, 2} | }{2} \Big ) \,.
\end{array} 
\end{equation}

From \eqref{f:spsi} and \eqref{f:sU} we obtain the estimate
 $$ \nu _{\psi'} \big ( J _n (t)   \big ) -  \nu _{V} \big ( \widetilde J _n (t)   \big ) \geq \mathcal R _{\overline \alpha}\,.$$ 

Otherwise, if one or both the multiindices $(\overline \alpha, 1)$ and $(\overline \alpha, 2 )$ belong to $\Gamma$,  the left hand side of \eqref{f:sU} does not correspond to  the measure $ \nu _{V} \big ( \widetilde J _n (t)   \big )$. Assume for definiteness that $(\overline \alpha, 1)  \in \Gamma$ and $(\overline \alpha, 2 ) \not \in \Gamma$. Then we split  $ I ^{\overline \alpha , 1}$ as  $I ^{\overline \alpha , 1,1} \cup   I ^{\overline \alpha , 1,2}$,    and we rewrite the left hand side of \eqref{f:sU}  as 
\begin{equation}\label{f:sUbis}  \begin{array}{ll} 
& 2 \tan \Big ( \frac {   |  \widetilde J _n ( t) | }{2} \Big )    - 
 2 \tan \Big ( \frac {| I    ^{\overline \alpha, 1} \cup   I ^{\overline \alpha, 2} | }{2} \Big ) +  
2 \tan \Big ( \frac { |  I  ^{\overline \alpha, 1} | }{2} \Big )    + 2 \tan \Big ( \frac { |  I  ^{\overline \alpha, 2} | }{2} \Big ) =  
\\ \noalign {\medskip} &  
2 \tan \Big ( \frac {   |  \widetilde J _n ( t) | }{2} \Big )    - 
 2 \tan \Big ( \frac {|  I    ^{\overline \alpha, 1} \cup   I ^{\overline \alpha, 2} | }{2} \Big ) + 2 \tan \Big ( \frac { |  I  ^{\overline \alpha, 2} | }{2} \Big ) 
\\ \noalign {\medskip} &   
+ 2 \tan \Big ( \frac { |  I  ^{\overline \alpha, 1,1} | }{2} \Big )    + 2 \tan \Big ( \frac { |  I  ^{\overline \alpha, 1, 2} | }{2} \Big ) 
+ \mathcal R _{\overline \alpha, 1}\,,   
\end{array} 
 \end{equation}  
with 
$$\mathcal R _{\overline \alpha,1} := 2 \tan \Big ( \frac { |  I  ^{\overline \alpha, 1,1} \cup  I ^{\overline \alpha, 1,2} | }{2} \Big )  - 2 \tan \Big ( \frac { |  I  ^{\overline \alpha, 1,1} | }{2} \Big )    - 2 \tan \Big ( \frac { |  I  ^{\overline \alpha, 1, 2} | }{2} \Big ) \geq 0 \,.$$ 
 In case the multiindices  $(\overline \alpha, 1,1)$ and $(\overline \alpha, 1,2 )$ do not belong to $\Gamma$, we have  
$$ \begin{array}{ll}  &  \displaystyle 2 \tan \Big ( \frac {   |  \widetilde J _n ( t) | }{2} \Big )    - 
 2 \tan \Big ( \frac {|  I    ^{\overline \alpha, 1} \cup   I ^{\overline \alpha, 2} | }{2} \Big )  \\ \noalign{\medskip} & \displaystyle + 2 \tan \Big ( \frac { |  I  ^{\overline \alpha, 2} | }{2} \Big ) 
 + 2 \tan \Big ( \frac { |  I  ^{\overline \alpha, 1,1} | }{2} \Big )    + 2 \tan \Big ( \frac { |  I  ^{\overline \alpha, 1, 2} | }{2} \Big ) 
=\nu _{V} \big ( \widetilde J _n (t)   \big ) \end{array} $$ 

 Hence from \eqref{f:spsi} and \eqref{f:sUbis} we obtain the estimate
 $$ \nu _{\psi'} \big ( J _n (t)   \big ) -  \nu _{V} \big ( \widetilde J _n (t)   \big ) \geq \mathcal R _{\overline \alpha }+ \mathcal R _{\overline \alpha, 1}\,.$$ 
Otherwise, we continue the procedure by splitting one or both the intervals $I ^{\overline \alpha , 1,1}$ and $I ^{\overline \alpha , 1,2}$.  
 In a finite number of steps, we arrive at the conclusion that 

 $$ \nu _{\psi'} \big ( J _n (t)   \big ) -  \nu _{V} \big ( \widetilde J _n (t)   \big ) \geq  \delta _n ( t):= \sum _\alpha \mathcal R _\alpha \geq 0 \,,$$ 
 where the sum is extended to all indices $\alpha \in \Gamma$ of the form $(\overline \alpha (n, t), \dots)$.  
Then the inequality \eqref{f:bbelow0} is proved because 
each term  $\mathcal R _\alpha$ is non-negative, and hence 
$$
\int_0 ^{\| v \| _\infty   }  2t   \sum _{n=1} ^ { N (t)}  \delta _n ( t) \, dt  \geq 0 \,.
$$
\qed

\medskip

\begin{lemma} \label{l:lemmaeta}   Let $V$ be the stratified  potential associated with some function 
$v \in H ^ 1 _0 (I _\pi) \cap \mathcal X( I _\pi)$,  and let   $\Gamma $ denote the set of indices in \eqref{f:gamma}. For any $\alpha \in \Gamma $, 
set
$ I ^ {\alpha , 1} \cup I ^ {\alpha, 2} = (a^\alpha _1, a ^ \alpha_2) \cup (a ^\alpha _2, a^\alpha_3)$, and consider the following  eigenvalue problem  with stratified potential 
$$   \eta _ 1 (I _\pi, 2V):= \inf \left\{   \frac{\int _{I _\pi} 
(\varphi ') ^ 2 + 2 V  \varphi ^ 2    }{\int _{I_\pi} \varphi  ^ 2 } \! :\!\varphi  \in H ^ 1 _0 ( I _\pi) \text{ s.t., }  
\forall \alpha \in \Gamma , \ 
\varphi  ( a ^ \alpha _ 1) =    \varphi  ( a ^ \alpha _ 2) =  \varphi  ( a ^ \alpha _ 3) \! \right\}\!.$$ 
Then it holds
 \begin{equation}\label{f:eta4} 
 \eta_ 1 (I_\pi, 2 V) \geq 4\,,
\end{equation} 
with equality if and only if $\Gamma = \emptyset$, so that $V (x) = 1 +  \tan ^ 2 x$  and an eigenfunction is $\cos ^ 2 \!x$.
\end{lemma}  

\proof  
  Throughout this proof,   
since we are going to work with different stratified potentials, in order to avoid any confusion  we denote by $\Gamma_V$ the family of multi-indices associated with 
a potential $V$ according to \eqref{f:gamma}.   Moreover, we write for shortness that $\varphi$ is an eigenfunction associated with $V$, if it is an eigenfunction for the eigenvalue problem $\eta _ 1 (I _\pi, 2 V)$. 
We argue by contradiction.  Assume  that the family $\mathcal S$ of stratified potentials such that  $\eta _ 1 (I _\pi, 2V) < 4$ is nonempty. We proceed in two steps.

 \medskip
 {\it Step 1}.  Since by assumption $\mathcal S \neq \emptyset$, we can select a stratified potential $V$ for which the cardinality of  $\Gamma _V$ is minimal among potentials in $\mathcal S$. Clearly, since $V\in \mathcal S$,  ${\rm card} (\Gamma_{V }) \geq 1$
(recall that, for $V (x) = 1 +  \tan ^ 2 x$, we have $\eta _ 1 (I _\pi, 2V) = 4$, with eigenfunction $\cos ^ 2 \!x$).
By suitably perturbing $V$, we are going to find  a potential $V ^\e$ such that
 $$ \eta _ 1 (I_ \pi, 2V^\e ) = 4  \qquad \text{ and  }\qquad     {\rm card} (\Gamma_{V ^ \e}) =  {\rm card} (\Gamma_{V }) \,. $$
 
 To that aim, let  $- \frac{\pi}{2} = a_0 < a_1 < \dots < a_K = \frac{\pi}{2} $ denote  an increasing  relabelling of the family of points  $\{ a ^ \alpha_i, \ i = 1,2,3 \}$ associated with $V$.  In particular, for some $2 \leq j \leq K -2$ we have $I ^ {1, 1} = (a_1, a_j)$, 
 $I ^ {1, 2} = (a_j, a_{K-1})$, and $a_{K-1} = - a _ 1$.  
 We consider a one-parameter family $\{ V ^ \e \}$ of continuous perturbations of $V$,  given by stratified potentials such that
 ${\rm card} (\Gamma_{V }) = {\rm card} (\Gamma_{V^ \e })  $,   in  which in particular the points 
 $a_1$ and $a _{K-1}$  are replaced by  $a _1 -\e$ and $a _{K-1} + \e$ (inside the interval 
  $[a _1 -\e,  a _{K-1} + \e]$, the potential $V ^ \e$ can be built by rescaling the family of points $a _i$ into the new
  family of points $a _i ^ \e = a _i \frac{ a_{K-1}  + \e } { a_{K-1} }$).    
  
This operation can be carried over for $\e \in (0,\e ^* )$, with  $\e ^*=  \frac{\pi}{2} - a _{K-1}$.   
Clearly, there exists $\e$ in such interval  such that 
$ \eta _ 1 (I _\pi, 2V ^ \e) = 4$. 
Otherwise, it would be 
$ \eta _ 1 (I _\pi, 2V ^ {\e ^*}) \leq 4$, which implies that  an eigenfunction $\varphi$ 
for $V ^ {\e ^*}$  must vanish at $a_1 ^ { \e ^*} = - \frac{\pi}{2}$, $a_{K-1} ^ { \e ^*} =  \frac{\pi}{2}$, and at $a _j ^ { \e ^*}$. 
Then 
the interval $( - \frac{\pi}{2}, \frac{\pi}{2})$ would be disconnected into the union of the  two disjoint intervals
$ I ^ { 1, 1} _* := (  - \frac{\pi}{2} , a_j ^{\e*})$ and $I ^ { 1, 2} _* := (   a_j ^{\e*},  \frac{\pi}{2} )$. The restriction of the potential to
one of the two intervals, say $I ^ { 1, 1} _*$, would give an eigenvalue less than or equal to $4$. 
Since the length of the interval $I ^ { 1, 1} _*$  is strictly less than $\pi$,  
this proves the existence of a stratified potential $\widehat V$ on $I _ \pi$ with 
${\rm card} (\Gamma _{\widehat V}) < {\rm card} (\Gamma _V)$ and  $\eta _ 1 ( I _\pi, \widehat V) < 4$.   
(The stratified potential $\widehat V$ can be obtained by centering $I ^ { 1, 1} _*$ at the origin, and
extending $\widehat V$ by setting it equal to $1 + \tan ^ 2 x$  on its complement in  $I _ \pi$.)

We conclude that it is possible to freeze $\e$ so that  $ \eta _ 1 (I _\pi, 2V ^ \e)  = 4$. The potential $V ^ \e$ satisfies
  \begin{equation} \label{f:requests} 
  \eta _ 1 (I_ \pi, 2V^\e )  = 4 \qquad \text{ and  }\qquad    
  {\rm card} (\Gamma_{V })  \geq    {\rm card} (\Gamma_{V ^ \e})   \ \forall V \in \mathcal S \,. \end{equation}  
 
\medskip
{\it Step 2}. Let $V ^ \e$ be a stratified potential satisfying \eqref{f:requests}. For simplicity of notation, in the remaining of the proof we drop the index $\e$ and we denote it simply by $V$.  
Let $\varphi$ be
an  eigenfunction for $ \eta _ 1 (I _\pi, 2V)$.  By optimality,  if
$- \frac{\pi}{2} = a_0 < a_1 < \dots < a_K  = \frac{\pi}{2} $ denotes  an increasing  relabelling of the family of points  $\{ a ^ \alpha_i, \ i = 1,2,3 \}$ associated with $V$,  
 $\varphi$ satisfies the system of 
equations
\begin{equation*}
- \varphi'' + 2 V \varphi =    4 
\varphi  \qquad \text{ on } (a_k, a _{k+1}) \qquad \forall k = 0, \dots K-1\,, 
\end{equation*} 
and the following global equality of the  boundary terms, where  $\varphi' _+ (a_{k}) $,  $\varphi' _- (a_{k+1})$ are the right and left derivatives  of $\varphi$ respectively at $a_k$ and $a_{k+1}$:

\begin{equation}\label{f:energy} \sum _{k = 0} ^{K-1}\big [  \varphi(a_{k+1}) \varphi' _- (a_{k+1}) - \varphi(a_{k}) \varphi' _+ (a_{k})     \big ]  = 0 \,,
\end{equation}

We are now ready to reach a contradiction, by distinguishing  the two cases ${\rm card} (\Gamma_{V}) = 1$ and ${\rm card} (\Gamma_{V }) >1$. 
 
{\it Case a. }  If ${\rm card} (\Gamma_{V }) = 1$, a first eigenfunction for $ \eta _ 1 (I _\pi, 2V)$  is explicitly determined as

$$\varphi (x)  = 
\begin{cases} 
\displaystyle C \frac{ \cos ^2 x } {\cos ^ 2 (a_1) } & \text{ on } (a_0, a_1)
\\  \noalign{\medskip} 
\displaystyle C  \frac{\cos ^2   (x- (\frac{a_1 + a_{2}}{2}  )  )}{\cos ^ 2 (\frac{ a _{2}- a _1 }{2}) }  & \text{ on } (a_1, a_2) 
\\   \noalign{\medskip}  
\displaystyle C  \frac{\cos ^2  (x- (\frac{a_2 + a_{3}}{2}  )  )}{\cos ^ 2 (\frac{ a _{3}- a _2 }{2}) } & \text{ on } (a_2, a_3)  
\\ \noalign{\medskip} 
\displaystyle C \frac{ \cos ^2 x } {\cos ^ 2 (a_3) } & \text{ on } (a_3, a_4)\,.
 \end{cases}
$$

Then the sum in \eqref{f:energy} can be written as
$$2  C^2   \Big [ - \tan ( a_1) - 2 \tan \Big (\frac {a_2 - a_ 1 }{2} \Big ) -   2 \tan \Big ( \frac{a_3 - a_ 2 }{2} \Big ) -  \tan (  a _{3} ) \Big ]\,.$$  
Using the elementary inequality 
$$\tan (a_{k+1}) - \tan (a_k) > 2 \tan \Big (\frac {a_{k+1} - a_ k }{2} \Big )  \qquad \text{ for } a_{k+1} > a _k\,$$ 
we see that  the equality \eqref{f:energy} cannot hold, contradiction. 

\medskip
{\it Case b. }  If ${\rm card} (\Gamma_{V }) > 1$, we   consider the highest level among the local minima of $\varphi$, and denote by $a_p$ the point where it is attained. Then  there are two consecutive intervals, say $(a_{p-1}, a_{p})$  and $(a_p, a _{p+1})$ such that a  first eigenfunction 
 for $ \eta _ 1 (I _\pi, 2V)$  satisfies, for some positive constant $C_p$ 
$$\varphi ( x) = \begin{cases}
\displaystyle C _p  \frac{\cos ^2  \Big (x- \big (\frac{a_{p-1} + a_{p}}{2}  )\Big ) }{\cos ^ 2 (\frac{ a _{p}- a _{p-1} }{2}) } \qquad \forall x \in (a_{p-1}, a _{p}) 
\\
\displaystyle C_p   \frac{\cos ^2  \Big (x- \big (\frac{a_{p} + a_{p+1}}{2}  )\Big ) }{\cos ^ 2 (\frac{ a _{p+1}- a _p }{2}) } \qquad \forall x \in (a_p, a _{p+1}) \,.
\end{cases} 
$$
 Since we are assuming that $ \eta _ 1 (I _\pi, 2V) = 4$, we have
 $$ \frac{\int _{I _\pi} (\varphi') ^ 2 + 2 V  \varphi ^ 2    }{\int _{I_\pi} \varphi ^ 2 } = 4  \quad \text{ and } \quad   \sum _{k = 0} ^{K-1}\big [  \varphi  (a_{k+1}) \varphi'  _- (a_{k+1}) -
  \varphi (a_{k}) \varphi' _+ (a_{k})     \big ]     = 0  \,.$$

We then modify the 
potential $V$  into a new potential $\overline V$ which differs from it uniquely on the interval $(a_p, a_{p+1})$ by setting 
$$\overline V ( x) = 1 + \tan ^ 2 \big ( x - \frac{a_p+ a_{p+1}}{2} \big ) \qquad \forall x \in ( a_p, a _{p+1})\,.$$  
Accordingly, we modify the function $\varphi$ uniquely on the interval $(a_p, a_{p+1})$ by setting  
$$\overline \varphi ( x) = 
C _p  \frac{\cos ^2  \Big (x- \big (\frac{a_{p-1} + a_{p+1}}{2}  )\Big ) }{\cos ^ 2 (\frac{ a _{p+1}- a _{p-1} }{2}) } \qquad \forall x \in (a_{p-1}, a _{p+1}) \,.
$$
Then,  on each of the intervals associated with the potential $\overline V$, the function $\overline \varphi$ still satisfies  the same PDE as $v$, namely we have
\begin{equation}\label{f:eqinterno} - \overline \varphi'' +  2 \overline V \overline \varphi = 4 \overline  \varphi \qquad \text{ on } (a_1,  a_ 2) \cup  \dots \cup (a_{p-1}, a_{p+1}) \cup \dots
\cup ( a_{K-1}, a _K) \,.\end{equation}
On the other hand, thanks to \eqref{f:energy} and the super-additivity  of the tangent function on $(0, \frac{\pi}{2})$, it holds 
 \begin{equation}
\label{f:vbordo}\begin{array}{ll}
  & \displaystyle   \sum _{k = 0} ^{K-1}\big [  \overline \varphi  (a_{k+1})  \overline \varphi'  _- (a_{k+1}) - \overline \varphi  (a_{k}) \overline \varphi'_+ (a_{k})     \big ]    \\ 
  \noalign{\smallskip} =  & \displaystyle  
  \sum _{k = 0} ^{K-1}\big [ \overline  \varphi  (a_{k+1})    \overline \varphi'  _- (a_{k+1}) - \overline  \varphi  (a_{k})  \overline \varphi' _-  (a_{k})     \big ]   \ 
\\ \noalign{\medskip} 
& + 2 \tan \Big (\frac {a_p - a_ {p-1} }{2} \Big ) +   2 \tan \Big ( \frac{a_{p+1} - a_ p }{2} \Big )  - 2 \tan \Big (\frac {a_{p+1} - a_ {p-1} }{2} \Big )  <0 \,.
\end{array} 
\end{equation}
By combining \eqref{f:eqinterno} and \eqref{f:vbordo}, we infer that the Rayleigh quotient with potential $2 \overline V$ of the function $\overline \varphi$
is strictly smaller than $4$, i.e. 
$$   \frac{\int _{I _\pi} (\overline \varphi') ^ 2 +  2 \overline V \overline \varphi ^ 2    }{\int _{I_\pi} \overline \varphi ^ 2 } < 4 \,.$$ 
We conclude that $ \eta _ 1 (I _\pi, 2 \overline V) < 4$.  Since ${\rm card} (\Gamma _{\overline V}) <  {\rm card} (\Gamma _{V})$,  this contradicts condition \eqref{f:requests} satisfied by $V$, and 
the proof of inequality \eqref{f:eta4} is achieved. 
 
\smallskip
 Concerning the equality case, 
assume that  $\eta _ 1 (I _\pi, 2V) = 4$ holds, and assume by contradiction that $\Gamma _V \geq 1$. If $\Gamma _V = 1$ by arguing as in Case a.\ above we obtain a contradiction; if $\Gamma _V >1$, by arguing as in Case b.\ above we arrive to  contradict,  for another potential   $\overline V$,  the inequality $\eta _ 1 (I, \overline V) \geq 4$  (that we have already proved). We conclude that $\Gamma_V = 0$, namely that  $V (x) = 1 +  \tan ^ 2 x$.

\qed

\bigskip
{\it Proof of Proposition \ref{p:seconda}}.   Let $\psi \in \A (I _\pi)$. Assume by contradiction that 
$$\lambda _ 1 (I_\pi, 2 \psi') = \inf  _{ v \in H ^ 1 _0 ( I _\pi) } \frac{\int _{I _\pi} (v') ^ 2 \, dx + 2  \int _{I _\pi}  v ^ 2 \, d \nu_{\psi'}    }{\int _{I_\pi} v ^ 2 }  < 4  \,.$$ 

By a density argument, we can find a function  $v \in  H ^1 _0( I _\pi)\cap \mathcal X (I _\pi)$  such that 
$$\frac{\int _{I _\pi} (v') ^ 2 \, dx + 2  \int _{I _\pi}  v ^ 2 \, d \nu_{\psi'}   }{\int _{I_\pi} v ^ 2 } < 4\,.$$

Let $\tilde v$ and $V$ denote the stratified rearrangement and potential associated with $v$. 

 By exploiting, at each step of the construction of $\tilde v$, the well-known behaviour under decreasing rearrangement  of the $L ^2$-norm of a function and of its first derivative,  we infer that 
$$
 \int _{ I _\pi} v ^ 2 = \int _{I_\pi} \tilde v ^ 2 \qquad \hbox{ and } \qquad \int _{ I _\pi} (v' ) ^ 2 \geq\int _{I_\pi} (\tilde v') ^ 2 \,.
$$
 
 Then, by
Lemma  \ref{l:lemmatilde}, we have 
$$\frac{\int _{I _\pi} (\tilde v') ^ 2 + 2 V \tilde v ^ 2    }{\int _{I_\pi} \tilde v ^ 2 } < 4\,.$$
 Therefore, for the stratified potential $V$, the eigenvalue  $ \eta _ 1 (I _\pi, 2V)$  introduced in Lemma \ref{l:lemmaeta} would be strictly smaller than $4$, contradicting such lemma.   
 \qed

\bigskip

We now turn attention the problem of estimating  from below $\lambda _ 1 (I _ d, q)$, where  the weight is still of the form $q= \psi' + \psi ^ 2$ for $\psi \in \mathcal A ( I _\pi)$, and the finiteness domain of $\psi$ is an interval
$I _ d= ( - \frac{d}{2}, \frac{d}{2})$ with  $d \leq \pi$.  

\begin{theorem}\label{t:1dextra}  

Inequality \eqref{f:stima3} can be refined as follows: 

\smallskip
\begin{itemize}
\item[(i)]  There exists an  absolute constant $C>0$  such that, 
for any  $q = \psi' + \psi ^ 2$, being $\psi \in \mathcal A ( I _\pi)$  with ${\rm dom } (\psi) = I _ d$ ($d \leq \pi$), it holds 
\begin{equation} \label{f:stima3length}  
\lambda _ 1  (I _{d}, q) \geq 3  + C(\pi - d) ^ 3 \,. 
\end{equation} 

\medskip 
\item[(ii)] There exists an  absolute constant $K$  such that, 
for any  $q = \psi' + \psi ^ 2$, with  $\psi \in \mathcal A ( I _\pi)$ of the form 
 $ \psi =  ( f + \frac{g}{2}   \big )$, being 
$f \in \mathcal A ( I _\pi)$ with  ${\rm dom} ( f) = I _ d$ ($d \leq \pi$), and 
 $g =  - (\log h)'  $ on $I _d$, for some  $({ \frac{1}{m}})$-concave function $h$,  
affine on an interval $[a,b]$  with 
$I_{ \frac{\pi}{4} } \subseteq [a, b] \subseteq I _d$,
 the following implication holds: 
\begin{equation}\label{f:stima3delta} 
\lambda _1 ( I _ d, q ) \leq 7 \ \Rightarrow \ \lambda _1 ( I _ d, q ) \geq 3 +\frac{ 8 K}{ m  \pi ^ 2} 
  \Big [ 1 - \frac{ \min\{ h ( a), h (b)  \}   }{ \max \{ h ( a), h (b)  \} } \Big ] ^2    
\,.
 \end{equation}
 \end{itemize}
 \end{theorem} 

\bigskip
Theorem \ref{t:1dextra} is obtained by an analogue strategy as Theorem \ref{t:1d}, replacing the use of Propositions \ref{p:prima} and \ref{p:seconda} 
by their refined versions stated respectively  in Propositions \ref{p:quarta} and \ref{p:terza} below. 
More precisely, one needs first to apply
the inequality $\lambda _1 ( I _ d, q) \geq  \frac{1}{2}\big [  \lambda _1 ( I _ d, 2 \psi' ) + \lambda _1 ( I _ d, 2 \psi^2 ) \big ]$. 
Then: inequality \eqref{f:stima3length} follows by using  Proposition \ref{p:seconda} to estimate $\lambda _ 1 (I _ d, 2 \psi') $ and
 Proposition \ref{p:quarta} below to estimate $\lambda _1 ( I _ d, 2 \psi ^ 2)$; inequality \eqref{f:stima3delta} follows by using Proposition \ref{p:prima} to estimate $\lambda _ 1 (I _d, 2 \psi ^ 2)$, 
and Propostion  \ref{p:terza} below to estimate $\lambda _1 ( I _ d, 2 \psi')$.

\begin{proposition}\label{p:quarta} 
There exists an  absolute constant $C >0$  such that, 
for any   $\psi \in \mathcal A ( I _\pi)$  with ${\rm dom } (\psi) = I _ d$ ($d \leq \pi$), it holds
\begin{equation}\label{f:stima2length}
\lambda _ 1  (I _{d}, 2 \psi ^ 2 ) \geq 2  + C (\pi - d ) ^ 3 \,.
\end{equation} 
\end{proposition} 
 
\proof By following the proof of Proposition \ref{p:prima}, we arrive at the inequality
$$ \lambda _ 1  (I _{d} , 2 \psi ^ 2 )  \geq \lambda _ {d}:= \min _{ v \in H ^ 1 _0 ( I _{d}) } 
 \frac{\int _{I _{d}} (v') ^ 2 +  2(\tan ^2\!  x)     v ^ 2  }{\int _{I_{d}} v ^ 2 } \,. 
$$ 
 Thus we are reduced to prove that there exists an  absolute constant $C>0$  such that $\lambda _ d \geq 2 + C (\pi -d ) ^ 3$ for every $d \in (0, \pi]$.  We observe that it is enough to prove that there exists $\overline \e>0$ and an absolute constant $C>0$ such that 
 $\lambda _ d \geq 2 + C (\pi -d ) ^ 3$ for every $d \in [\pi -  \overline \e, \pi]$.    Indeed in this case,
since
the map $d \mapsto \lambda _d$ is nonincreasing,  for every $d \in  (0, \pi - \overline \e)$ we have  $$\lambda _ d \geq \lambda _{ \pi - \overline \e}  \geq 2 + C\,  \overline \e ^ 3 = 2 + C' ( \pi-d) ^3\qquad \text{ with } C' = C \frac{ \overline \e ^ 3 }{ ( \pi - d) ^ 3 }  \geq C \frac{ \overline \e ^ 3 }{ \pi  ^ 3 } \,,  $$ 
so that  the required inequality is satisfied (for another 
 absolute constant) also for $d \in  (0, \pi - \overline \e)$.    
 Hence, in the remaining of the proof, we focus attention on the estimate of $\lambda _ { \pi - \e}$, for $\e$ sufficiently small. Denoting by $v _ \e \in H ^ 1 _0 ( I _{\pi-\e}) $  an eigenfunction for $\lambda _{\pi-\e}$,  we have 
\begin{eqnarray}
& - ( \cos^2 \! x ) ''  + 2 (\tan^2 \! x)  (\cos^2 \! x)   = 2  \cos^2  \! x \qquad \text { in }  I _\pi & \label{f:deb1}
\\ \noalign{\medskip}
& - v _\e '' + 2 (\tan^2  x)  v _ \e = \lambda _{\pi-\e}  v _ \e \qquad \text { in }  I _{\pi-\e}\,. & \label{f:deb2}
\end{eqnarray} 
We multiply \eqref{f:deb1} by $v _\e$ (extended to $0$ on $I _\pi \setminus I _ {\pi-\e}$), and \eqref{f:deb2} by $\cos^2 x $, and we integrate, respectively, on $I _\pi$ and on $I _ {\pi-\e}$. We get: 
  $$ \begin{array}{ll}
&  \displaystyle \int _{ I _\pi} ( \cos^2  x ) ' v _ \e '   + 2 \int _ {I _\pi}   (\sin^2 x)   v _ \e  = 2   \int _{I_\pi} (\cos^2 x)   v _ \e   
\\ \noalign{\medskip}
 \displaystyle &\displaystyle    \int _{ I _{\pi-\e}}  \!\!\!  ( \cos^2 x ) ' v _\e '   - 2 v ' _\e    \big (\frac{\pi - \e }{2}  \big ) \cos ^ 2\big (\frac{\pi - \e}{2} \big )    +   \int _{ I _{\pi- \e}} \!\!\! \!\!\! 2 (\sin^2 x)  v _ \e = \lambda _{\pi-\e}  \int_{I _ {\pi-\e}}  \!\!\!  (\cos^2 x)   v _ \e \,.  
\end{array} $$
By subtraction, we obtain  
$$ (\lambda _{\pi-\e}-2)   \int_{I _{\pi-\e}}  (\cos^2 \! x)   v _ \e   =  - 2  \big (\sin ^ 2 \! \frac{\e}{2} \big )  v ' _\e  \big (\frac{\pi - \e}{2}  \big )\,. $$ 
Since it is easily checked that $v _ \e$ converges weakly to $\cos^2\!  x$ in $H ^ 1 _0 ( I_\pi)$, 
we have that
$$\lim _{\e \to 0 } \int_{I _ {\pi- \e}}  (\cos^2 \! x)   v _ \e    =\int_{I _ \pi}  (\cos^4 \! x)    \in (0, + \infty)\,.$$ 
Therefore, to prove \eqref{f:stima2length}, it is enough to show that 
\begin{equation}\label{f:ordre1}
\lim _{\e \to 0 } \frac{1}{\e}  v ' _\e  \big (\frac{\pi - \e}{2}  \big )  \in (0, + \infty)\,.
\end{equation} 
To that aim, we multiply equation \eqref{f:deb2} by $\frac{1}{\cos x}$, and we integrate on $I _ {\pi-\e}$. We obtain 
$$\int _{ I _{\pi-\e}}    \frac{\sin x}{ \cos^2 \! x }  v _\e '   - 2 \frac{ v ' _\e  \big (\frac{\pi-\e}{2}  \big )  }{\cos \big (\frac{\pi -\e}{2} \big )  }  + 2  
\int _{ I _{\pi-\e}} \frac{\sin ^2 \! x}{\cos ^ 3 \! x}  v _ \e = \lambda _{\pi-\e} \int_{I _{\pi-\e}} \frac{   v _ \e}{ \cos  x} $$ 
Since  
$$\int _{ I _{\pi-\e}}    \frac{\sin x}{ \cos^2 \! x }  v _\e '   = -  \int _{ I _{\pi-\e}}   \frac{v _\e }{\cos x}     - \int _{ I _{\pi-\e}} 2 \frac{\sin ^2 \! x}{\cos ^ 3 \! x} v _ \e \,, $$ 
we end up with
 \begin{equation}\label{f:endup}- 2 \frac{ v ' _\e  \big (\frac{\pi -\e}{2} \big )  }{\cos \big (\frac{\pi -\e}{2} \big )  } =      (\lambda _{\pi-\e} + 1 )  \int _{ I_{\pi-\e} }   \frac{v _\e }{\cos x}  \,.
\end{equation}
Finally, we observe that 
\begin{equation}\label{f:ordre1bis}
\lim _{\e \to 0 }  \int _{ I _{\pi-\e}}    \frac{v _\e }{\cos x}    \in (0, + \infty)\,.
\end{equation}  
Indeed, since $v _ \e$ converges to $\cos ^ 2 x$ a.e. on $I _\pi$, we have 
$$\liminf _{\e \to 0 }  \int _{ I _{\pi-\e}}    \frac{v _\e }{\cos x}   \geq  \int _{ I _\pi}    {\cos x}   >0\,.$$ 
On the other hand,  by  H\"older inequality we have 
$$  \limsup _{\e \to 0}  \int _{ I _{\pi-\e}}    \frac{v _\e }{\cos x}   \leq  \pi ^ {\frac{1}{2}}   \limsup _{\e \to 0} 
 \Big [  \int _{ I _{\pi-\e}}    \frac{v _\e ^2 }{\cos^2 \!  x} \Big ] ^ {\frac{1}{2}}   < + \infty\,,$$ 
where the last inequality is obtained by observing that,  normalizing $v _ \e$ in $L ^\infty$,
and recalling from 
 \eqref{f:deb2} that 
$\sup _ \e \int  _{I _{\pi-\e}}  2 (\tan^2  x)  v _ \e ^ 2  \leq   \pi  \sup _ \e  (   \lambda _\e  v _ \e   ^ 2 ) < + \infty$, 
it holds
 $$ \begin{array}{ll}
\displaystyle \int _{ I _{\pi-\e}}    \frac{v _\e ^2 }{\cos^2 \!  x}   & \displaystyle \leq     \int _{ I _{\pi-\e} \cap  \{ |x| \leq \frac{\pi}{4}  \}}  \frac{1 }{\cos^2 \!  x}  + 
\int _{ I _{\pi-\e} \cap  \{ |x| > \frac{\pi}{4}  \}}  ( \tan^2\!  x ) \frac{ v _ \e ^ 2}{  \sin ^ 2 \!  x}  
\\
 & \displaystyle \leq     \int _{ I _{\pi-\e} \cap  \{ |x| \leq \frac{\pi}{4}  \}}  \frac{1 }{\cos^2 \!  x}  + 
\int _{ I _{\pi-\e} \cap  \{ |x| > \frac{\pi}{4}  \}}  2 (\tan^2\!  x )  v _ \e ^ 2   
\leq C\,.
\end{array}$$ 
From \eqref{f:endup} and \eqref{f:ordre1bis}, we see that  \eqref{f:ordre1} is satisfied, so that our proof is achieved. 
\qed

\bigskip

\begin{proposition}\label{p:terza} 
There exists an absolute constant $K$ such that, 
for any   $\psi \in \mathcal A ( I _\pi)$ of the form 
 $ \psi =  ( f + \frac{g}{2}   \big )$, being 
$f \in \mathcal A ( I _\pi)$ with  ${\rm dom} ( f) = I _ d \subseteq I _\pi$, and 
 $g =  - (\log h)'  $ on $I _d$, for some  $({ \frac{1}{m}})$-concave function $h$,  
affine on an interval $[a,b]$  with 
$I_{ \frac{\pi}{4} } \subseteq [a, b] \subseteq I _d$,
 the following implication holds: 
\begin{equation}\label{f:stima3delta2} 
\lambda _1 ( I _ d, 2 \psi' ) \leq 15 \ \Rightarrow \ \lambda _1 ( I _ d, 2 \psi' ) \geq 4 +\frac{ 16 K}{ m  \pi ^ 2} 
 \Big [ 1 - \frac{ \min\{ h ( a), h (b)  \}   }{ \max \{ h ( a), h (b)  \} } \Big ] ^2  
\,.
\end{equation} 
\end{proposition}

\proof  We claim that,  for $\psi$ as in the assumptions, denoting by $\overline v \in  H ^ 1 _0(I _d)$ an  eigenfunction for  $\lambda _ 1 (I _d, 2\psi')$, normalized in $L ^ 2(I _d)$, 
it holds  
\begin{equation}\label{f:stima4delta}
\lambda _ 1 (I_d,  2 \psi ')  \geq 4+ \delta _h \, , \qquad \text{ with } \delta _h :=    
 {\frac{1}{m}   \int _{I _d} \overline v ^ 2 \big ( \frac{h'}{h}  \big ) ^ 2}  \,. 
\end{equation}

Indeed,  we have
\begin{equation}\label{f:spezz} 
\lambda _ 1 (I _d, 2\psi') = {\int _{I _d} (\overline v') ^ 2 + (2f' + g' ) \overline v ^ 2  }\geq 4 + {\int _{I _d }g' \overline v ^ 2  } \,,  
\end{equation} 

where  the inequality follows by applying Proposition \ref{p:seconda} to the function  $f \in \mathcal A ( I _\pi)$ 
(actually, $\overline v$ is an admissible test function for $\lambda _ 1 ( I  _\pi, 2 f')$ when extended to $0$ on $I _\pi\setminus I _ d$). 
Then the  inequality \eqref{f:stima4delta} follows from  \eqref{f:spezz}  provided 
\begin{equation}\label{f:g'1} g'  \geq 
  \frac{1}{m}    \Big ( \frac{h'}{h}  \Big ) ^ 2 \qquad \text{ on } I _ d  \,. \end{equation}
 From the assumption $ g= - (\log h)'$ on $I _d$, we have (in the sense of measures)
 \begin{equation}\label{f:g'2} g' =  \Big ( \frac{h'}{h}  \Big ) ^ 2 - \frac{h''}{h}\qquad \text{ on } I _ d    \,.
\end{equation} 
The power-concavity assumption  on  $h$ implies that
$(h ^ { \frac{1}{m}})  ^{''} \leq 0$. 
The latter inequality, by an elementary computation, implies that the right hand side of \eqref{f:g'2} is larger than or equal to the right hand side of \eqref{f:g'1}. 
In view of \eqref{f:stima4delta},  we have 

$$\begin{array}{ll}  \lambda _1 ( I _ d, 2 \psi' )   & \displaystyle \geq   4+ 
 \frac{1}{ m  } \Big [ \frac{h ( a) - h ( b) }{ a-b  }  \Big ] ^ 2 \frac{1}{ [ \max \{ h ( a), h (b) \}   ] ^ 2}  \int _{I_{ \frac{\pi}{4} } }    \overline v ^ 2  
 \\ \noalign{\bigskip}   &\displaystyle  \geq 
 4+  \frac{16}{ m  \pi ^ 2 } \Big [ \frac{h ( a) - h ( b) }{ \max \{ h ( a), h (b) \}  }  \Big ] ^ 2   \int _{I_{ \frac{\pi}{4} } }    \overline v ^ 2 \,. 
   \end{array} 
$$

Hence, to conclude the proof of \eqref{f:stima3delta2}, it is enough to show that there exists  an absolute constant $K>0$  such that 
$$\lambda _1 ( I _ d, 2 \psi' ) \leq 15  \ \Rightarrow  { \int  _{I_{ \frac{\pi}{4}}}  \overline v ^ 2 }\geq K\,.$$ 
Assume by contradiction this is false. Then it would be possible to find a sequence  of  functions $\psi _n$ and a sequence of segments $I _ { d_n}$ 
as in the assumptions of the Lemma such that
the  eigenfunctions $\overline v_n \in  H ^ 1 _0(I _{d_n})$ for  $\lambda _ 1 (I _{d_n}, 2\psi' _n)$, 
normalized in $L ^ 2 ( I _ {d_n})$,  satisfy 
$$ { \int _{I_{ \frac{\pi}{4}}}   \overline v _n ^ 2 }\leq \frac{1}{n}\,.$$  

By the assumption  $\lambda _ 1 ( I _ { d_n},  2 \psi' _n ) \leq 15$, 
 we have $\int _{ I _\pi} |\overline v ' _n| ^ 2 \leq 15$, and hence 
up to subsequences  $\overline v _n$ converges, weakly in $H ^ 1 (I_\pi)$ and strongly in $L ^ 2 (I _\pi)$, to a function  $\overline v _\infty$ which has unit norm in $L ^ 2 (I _\pi)$ and vanishes on   ${I_{ \frac{\pi}{4}}}$.  This leads to a contradiction, as 
$$15  \geq \liminf _n \int _{ I _\pi} |\overline v ' _n|^2 \geq  \int _{ I _\pi} |\overline v ' _\infty| ^ 2 \geq \lambda _ 1 (  {I_{ \frac{\pi}{4}}} ) = 16\,.$$
\qed

\bigskip

\section{Weighted equipartitions \`a la Payne-Weinberger}\label{sec:partitions}  

The key idea in the proof of Payne-Weinberger inequality in \cite{PW} is a partition procedure of the set $\Om$ into convex cells of equal measure, such that, on each cell, the first Neumann eigenfunction has to have zero integral mean. 
While this procedure is still useful if adapted to a weighted Neumann problem, 
in order to obtain a quantitative estimate it is  necessary  to keep track of the $L^2$-norm of the eigenfunctions rather than of the measures of the cells.  Consequently, we are going  to  work with two different types of $p$-weighted equipartitions, each one playing a specific role in the estimate of
weighted Neumann eigenvalues in higher dimensions. Such estimate will be based the one-dimensional lower bounds given in Theorem \ref{t:1dextra}. Thus
we are going to handle line segments contained into $\R ^N$: {\it for simplicity, in this section and in the remaining of the paper, the notation $I _ \ell$  is 
adopted for
any  line segment of length $\ell$ in $\R ^N$, say with generic direction and  not necessarily centred at the origin
(as it was the case in Section \ref{sec:1D}).  In the few cases when we have to consider a centred interval in a fixed frame,  
this will be explicitly indicated, by writing  e.g. $( -\frac{\ell}{2}, \frac{\ell}{2}) \times \{ 0 \}$. }

\begin{definition}
Given an open bounded convex set $\om \subset \R ^N$, a positive weight $p \in L ^ 1 (\om)$, and a function $u\in L ^ 2  (\om, p \, dx)$ satisfying $\int _\om u p = 0$, 
we call:

\medskip
\begin{itemize}

\item[--]  a {\it $p$-weighted measure equipartition  of $u$  in $\om$} a family $\mathcal P _n = \{ \om _ 1, \dots, \om _n \}$, where $\om _i$ are mutually disjoint convex sets  such that $\om = \om _1 \cup \dots \cup  \om _n $ and 
$$\int _{\om _i } u  p = 0 \qquad \text{ and } \qquad |\om _i|    = \frac{1}{n} |\om|  \qquad \forall i = 1, \dots, n\,. $$

\medskip

\item[--]  a {\it $p$-weighted $L^2$ equipartition of $u$  in $\om$} a family $\mathcal P _n = \{ \om _ 1, \dots, \om _n \}$, where $\om _i$ are mutually disjoint convex sets  such that $\om = \om _1 \cup \dots \cup  \om _n $ and 
$$\int _{\om _i } u  p = 0 \qquad \text{ and } \qquad \int _{\om _i }   u   ^ 2  p   = \frac{1}{n} \int _{\om  }   u   ^ 2  p  \qquad \forall i = 1, \dots, n\,. $$

\end{itemize}

\end{definition}

\begin{remark}\label{r:PW} 
 (i) The existence of $p$-weighted  $L^2$ (or measure) equipartitions of $u$  in $\om$ given by  two cells  is obtained by the analogue argument as in \cite{PW} (see also \cite{beb}).  Namely, for every $\alpha \in [0, 2 \pi]$, there exists a unique hyperplane with normal $(\cos \alpha_1, \sin \alpha_1 , 0, \dots, 0)$ which divides $\omega$ into two subsets $\omega _\alpha ' $ and $\omega _\alpha ''$ such that  $\int _{\om _\alpha' }   u   ^ 2  p  = \int _{\om _\alpha'' }  u   ^ 2  p $. Since the function $\mathcal I (\alpha) = \int _{\om _\alpha' }   u     p$ is continuous and satisfies $\mathcal I ( \alpha ) = - \mathcal I (\alpha + \pi)$,  there exists an angle $\overline \alpha$ such that $ \mathcal I (\overline \alpha) = 0$. 
Applying repeatedly the above argument
yields the existence  of a $p$-weighted $L^2$ equipartition  of $u$  in $\om$ given by  $n$  cells, each of them being contained into a narrow strip, 
 determined by two hyperplanes with normal of the form $(\cos \alpha_1, \sin \alpha_1, 0, \dots, 0)$ at infinitesimal distance from each other as $n \to + \infty$. 
The latter property follows from the fact that the volume of all the elements of the partition is infinitesimal as $n \to  + \infty$,   thanks to the assumption $p>0$ a.e. in $\om$. 

(ii) If the above procedure is repeated overall $(N-2)$  times, 
using as a last package of cutting hyperplanes those with normals of the form 
$ (0, \dots, 0, \cos \alpha _{N-2}, \sin \alpha _{N-2}, 0 )$, 
we obtain
a $p$-weighted $L^2$ equipartition  of $u$  in $\om$  into  mutually disjoint convex cells which are  narrow in $N-2$ directions, each one  orthogonal to $e _N$.  
If the procedure is repeated once more, we arrive at  a $p$-weighted $L^2$ equipartition   
of $\om$ into mutually disjoint convex cells of one dimensional type, being narrow in $(N-1)$ orthogonal directions. 
\end{remark}

Motivated by the above remark, we state the following single-cell estimate, holding for a set  $\om _\e \subset \om$  which is narrow in  $(N-1)$ orthogonal directions.

\begin{lemma}\label{l:stimaPW} Let $\om\subset \R ^N$ be an open bounded convex set,  and
let $p$ be a positive uniformly continuous weight defined in $\om$. 
Given $\e>0$,  let $\om_\e \subset \om$ be an open bounded convex set of diameter $d_\e$, 
which satisfies, 
in a suitable orthogonal coordinates system and for some $\e>0$, the inclusion
\begin{equation}\label{f:orizzontale} \om_\e \subset \Big \{ (x_1, y) \in \R \times \R ^ { N-1} \ : \ |x _1 | \leq   \frac{d_\e}{2}   ,\  |y _j| \leq \e \ \forall j = 1, \dots, N -1 \}\,.\end{equation}  
For every function $u\in W ^ {2, \infty}( \om)$ whose restriction to $\om _\e$  is an admissible test function for $\mu _ 1 (\om_\e ,p)$,  
setting  $h ( x) := \mathcal H ^ {N-1} ( \om_\e \cap \{ x _ 1 = x \} )$, 
it holds
 \begin{equation}\label{f:boundfilo} \frac{ \int _{\om_\e} |\nabla u| ^ 2    p}{\int _{\om _\e}| u| ^ 2    p  }     \geq   \Big \{ \mu _ 1 ( I _  {d_\e}, h p) - \frac{ \alpha (\e) |\om _\e| }{\int _{\om_\e}  |u|^2   p  } \Big  [ 1 + \mu _ 1 ( I _ {d_\e } , h p)  \Big ( 1 + \beta (\e) |\om_\e|   \Big )  \Big ] \Big \} \,,
\end{equation}
where $\alpha(\e)$ and $\beta (\e)$ are infinitesimal as $\e\to 0 $, depending only from  
$\| u \| _{W ^ {2, \infty}( \om) }$, $\| p \| _{ L ^ \infty (\om) }$, and from 
the modulus of continuity of  $p$ at $\e$ in $\om $.

 \end{lemma} 

\proof  Let $M$ be a positive constant such that 
$
\|  u \| _{W ^ {2, \infty}( \om) } \leq M$  and $\| p \| _{ L ^ \infty (\om) } \leq M$, 
 and let $\delta _\e>0$ be such that  $|p (x) - p ( y) | < \e$ for  every $x, y \in \om$ with $|x-y |  < \delta _ \e$.

We have  
$$
\begin{array}{ll}
&  \displaystyle \Big | \int _{\om_\e} \Big  (\frac{\partial u }{\partial x _1}  \Big ) ^ 2  p - \int _{I_{d_\e}}   [u (x, 0) '] ^ 2  h   p   \, dx \Big | \leq (   2 M ^ 3  \e  + M ^ 2 \delta _\e)  |\om_\e| =:\delta_1  
\\ \noalign{\medskip}
& \displaystyle \Big | \int _{\om_\e}  u ^ 2  p - \int _{I_{d_\e}}  u (x, 0)  ^ 2  h   p   \, dx \Big | \leq  
(2 M ^ 3  \e  + M ^ 2 \delta _\e) |\om_\e| =: \delta_2  ( = \delta _1) 
\\  \noalign{\medskip}
& \displaystyle   \Big  | \int _{I_{d_\e}}   u (x, 0)    h   p  \, dx  \Big | =   \Big | \int _{\om_\e}  u   p - \int _{I_{d_\e}}   u (x, 0)    h   p  \, dx \Big | \leq     (M^ 2 \e  + M \delta _\e )|\om _\e| =: \delta_3 \,.
\end{array} 
$$
Then, setting $\overline u := \int _{I_{d_\e}}  u (x,  0) h p$, we have
$$\begin{array}{ll} \displaystyle \int _{\om _\e} |\nabla u| ^ 2    p   &  \geq \displaystyle \int _{\om _\e} \Big  (\frac{\partial u }{\partial x _1}  \Big ) ^ 2  p \geq \int _{I_{d_\e}}    [u (x, 0) '] ^ 2  h   p   \, dx -  \delta_1 
\\ \noalign{\smallskip} 
& \geq \displaystyle \mu _ 1 ( I _ {d_\e}, h p) \Big [  \int_{I_{d_\e}}    u (x, 0) ^ 2 h p  - \overline u ^ 2 \Big ] - \delta_1
\\ \noalign{\smallskip} 
& \geq \displaystyle \mu _ 1 ( I _ {d_\e}, h p) \Big [ \int _{\om_\e}  u ^2      p   \, dx  - \delta _2  - \delta_3 ^ 2 \Big ]  - \delta _ 1
\\ \noalign{\smallskip} 
& = \displaystyle \mu _ 1 ( I _ {d_\e}, h p)  \Big[ \int _{\om_\e}  u ^2      p   \, dx \Big] - 
 \delta _ 1   \Big [   1+  \mu _ 1 ( I _ {d_\e}, h p) \Big ( 1  + \frac{ \delta_3 ^ 2}{\delta _1} \Big )  \Big ] \,.
 \end{array}$$ 

The result follows by inserting the expressions of $\delta_1$ and $\delta _3$ in the above estimate. 
\qed

\bigskip

The next two lemmas contain   lower bounds of different nature for $\mu_ 1 (\om, p)$ (or, more generally, for Rayleigh quotients).
The first lower bound,  stated in
 Lemma \ref{l:preinside},   
is given in terms of 
 one dimensional eigenvalues: it   is in fact obtained working with measure equipartitions 
 and applying the single cell estimate of Lemma \ref{l:stimaPW}.   The second lower bound, stated in  Lemma \ref{l:goodparti}, is given in terms of the average of the eigenvalues of the cells of the partition, and  is  obtained working with $L^2$ equipartitions.

\begin{lemma}\label{l:preinside}
Let $\om$ be an open bounded convex set, and let $p$ be a positive uniformly continuous weight in $\om$.  
Let $u\in W ^ { 2, \infty}  (\om)$ satisfy $\int _\om u p = 0$. 
Assume that,  for  every $n$ sufficiently large,  there exists a 
$p$-weighted measure equipartition  $\mathcal P _n = \{ \om _ 1, \dots, \om _n \}$   of $u$ in $\om$ such that  $\mu _ 1 ( I _ { d_i}, h _ i p) \geq c$ for every $i = 1, \dots, n$, where  $I _ { d_i}$  is a diameter for $\om _i$, and   $h _ i (x)$  is the $\mathcal H ^ {N-1}$-measure of the sections of $\om _i$ 
as in Lemma \ref{l:stimaPW}.  
Then
$$  \frac{ \int _{\om} |\nabla u| ^ 2    p}{\int _{\om}| u| ^ 2    p  }     \geq  c \,.$$ 
In particular, in case $\om$ is smooth and $p$ is smooth and strictly positive in $\om$, taking $u$ equal to a first eigenfunction for $\mu_1 (\om, p )$, we obtain that $\mu_1 (\om, p ) \geq c$.  
\end{lemma} 
\proof  For a given $\e>0$,  for $n$ large enough 
each of the sets $\omega _i$ satisfies in a suitable orthogonal coordinates system the inclusion \eqref{f:orizzontale}  (cf.\ Remark \ref{r:PW}). Moreover,  the restriction of $u$ to $\om _i$ is an admissible test function for $\mu _ 1 (\om_i ,p)$.
  Then by   
Lemma \ref{l:stimaPW} the inequality \eqref{f:boundfilo}  
is fulfilled for every $i = 1, \dots , n$, for some infinitesimal  $\alpha (\e) $ and $\beta (\e)$ 
which are independent of $i$.  (Indeed, 
as stated in Lemma \ref{l:stimaPW}, $\alpha (\e) $ and $\beta (\e)$  depend  only  on $\| u \| _{W ^ {2, \infty}( \om) }$, $\| p \| _{ L ^ \infty (\om) }$, and 
on the the modulus of continuity of  $p$ at $\e$ in $\om$). 

Writing the inequality \eqref{f:boundfilo} for each of the sets $\om _i$, and  recalling that  $|\om _i|=  \frac{|\om| }{n}$, we obtain 
$$   {\int _{\om _i } | u| ^ 2    p  }   \leq \frac{1}{ \mu _ 1 ( I _ {d_i}, h_i p) }   {   \int _{\om_i} |\nabla u| ^ 2    p  }      +\frac{1}{ \mu _ 1 ( I _ {d_i}, h_i p) }  \Big  ( {   \alpha ( \e)  \frac{|\om| }{n}   }  \Big  )   +  { \alpha ( \e)   \frac{|\om| }{n}   }   \Big ( 1 + \beta (\e)  \frac{|\om| }{n}   \Big )    \,. 
$$ 
  By using  the assumption $ \mu _ 1 ( I _ { d_i}, h _ i p) \geq c $ for every $i = 1, \dots, n$, 
and summing over $i = 1, \dots, n$, we get 
$$   {\int _{\om  } | u| ^ 2    p  }   \leq \frac{1}{ c }   {   \int _{\om} |\nabla u| ^ 2    p  }      +
\frac{1}{c }  \Big  ( { \alpha (\e)  |\om|   }  \Big  )   +  {\alpha (\e)  }   \Big ( 1 +   \beta (\e)  {|\om| }  \Big )    \,. 
$$  
The statement follows by letting $\e$ tend to $0$. 
   \qed

\bigskip

\begin{lemma}\label{l:goodparti} 
Let $\om$ be an open bounded convex set of diameter $d$ and let $p$ be a positive weight in $L ^ 1 (\om)$.  
Let $u  \in H ^ 1_{\rm loc} (\Om) \cap L ^ 2 (\Om, p  dx)$ satisfy $\int _\Om u p \, dx = 0$.  
If   $\mathcal P _n = \{ \om _1 , \dots,   \om _n\}$ is a $p$-weighted  $L^2\!$ equipartition  of $u$ in $\om$, it holds 
\begin{equation*}   \frac{  \int_\om |\nabla   u  | ^ 2 p  }{
 \int_\om |u  | ^ 2 p} 
 =  \frac{1}{n}    \sum _{i = 1} ^n\frac{  \int_{\om _i }  |\nabla   u  | ^ 2 p  }{  \int_{\om _i} |u  | ^ 2 p }  \geq 
  \frac{1}{n}  \sum _{i = 1} ^n  \mu _ 1 (\om _ i, p) \,.  
\end{equation*}
 In particular, if $u$ is  an eigenfunction for $\mu _1 (\om,p)$,   we have $\mu _ 1 (\om, p) \geq  \frac{1}{n}  \sum _{i = 1} ^n  \mu _ 1 (\om _ i, p) $ and, in case 
 $\mu _ 1 (\om, p) = \mu _ 1 (\om_i, p)$ for every $i = 1, \dots, n$, 
$u$ is necessarily an eigenfunction also for 
each $\mu _ 1 (\om_i, p)$. 
\end{lemma} 

\proof  
The statement is an immediate consequence of the two facts   that  $\int _{\om _i }   u  ^ 2  p   = \frac{1}{n} \int _{\om  }   u   ^ 2  p$ and 
 the restriction of $u$  to $\om _i$ is admissible as a test function for $\mu _ 1 (\om _1, p)$. 
In case $u$ is an eigenfunction for $\mu _ 1 (\om, p)$, and $\mu _ 1 (\om, p) = \mu _ 1 (\om_i, p )$ for every $i = 1, \dots, n$,  we see that none of the inequalities
$ \int_{\om _i }  |\nabla   u  | ^ 2 p   \geq 
   \mu _ 1 (\om _ i, p) \int_{\om _i }  |   u  | ^ 2  p$ can be strict.  
\qed

\section{Local strengthened version and rigidity of Andrews-Clutterbuck inequality}\label{sec:localized}

We are now in a position to give the  key  localized inequality \eqref{f:diffd.1} announced in the Introduction (and its weaker form \eqref{abf41}). 

 We stress that  the role of $\Om$ in the inequality \eqref{f:diffd.1}  
is just to precise the modulus of log-concavity of the weight. We keep this formulation
  since it fits 
  our purposes,  but we warn the reader that   the assertion  could   
be rephrased in terms of $(\om, p)$ only,   simply asking that 
 \eqref{f:improved2}  holds, for points $x, y \in \om$, with 
$D_\Om$ replaced by a length not smaller than $D _\om$.

  Let us also mention   that an inequality somehow similar to \eqref{abf41} appears in 
\cite[Theorem 4]{LR12}, without proof (merely invoking the techniques by Andrews-Clutterbuck), and seemingly with some inconsistencies in the assumption on the weight.

\begin{proposition}\label{p:inside}  
There exists an absolute constant $C>0$ such that for any   bounded open convex set $\Om\subset \R ^ N$ ($N \geq 1$) of diameter $D _\Om$,   and 
 any positive weight  $ p \in L ^ \infty (\Om)$ satisfying the improved log-concavity condition \eqref{f:improved2}, 
 if  $\Pi$ is an affine subspace of $\R^N$ of dimension $k \leq N$, and 
$\om  $ is a  relatively  open  convex subset of $\Om\cap \Pi$ of diameter $D_\om>0$ then inequality \eqref{f:diffd.1} holds, namely 
 \begin{equation*}\label{f:diffd} \mu _ 1 (\om, p ) \geq \frac{3 \pi ^ 2}{D_\Om^2}  + C \frac{(D_\Om-D_\om) ^ 3}{D _\Om^ 5} \,.
\end{equation*}
\end{proposition}

{\it Proof of Proposition \ref{p:inside}}.  The one dimensional case was discussed in Section \ref{sec:1D}. Let $N \ge 2$, and let us assume  
without loss of generality that $D_\Om=\pi$, in order to  easily handle the application of 
the one dimensional results from Section \ref{sec:1D} on the interval $I_\pi$.  
 We first prove the inequality \eqref{f:diffd.1}  in case 
$\omega$ is a smooth convex set  such that $\overline \om \subset \Om$.  In this case the restriction of $p$ to $\omega$ can be approximated by a sequence of functions $p _\eta$ which are  
$C ^ 2 (\overline \om)$ bounded away from $0$.   Then, 
there exists  a first eigenfunction for $\mu _1 (\om, p_\eta)$, of class $W ^ { 2, \infty} (\om)$.
So we are in a position to apply Lemma \ref{l:preinside}.  Specifically, we consider a
$p_\eta$-weighted measure equipartition  $\mathcal P _n = \{\om _1, \dots, \om _n \}$   
 of
a first eigenfunction for $\mu _1 (\om, p_\eta)$. 
We denote  by $I _ { d_i}$ a diameter for $\om _i$, and by    $h _ i (x)$ the $\mathcal H ^ { k-1}$-measure of the sections of $\om _i$  orthogonal to $I _ { d_i}$
as in Lemma \ref{l:stimaPW}.  
  We observe that, by the assumption made on $p_\eta$,  
\  each weight $h _ i p_\eta$    is of the form $\psi _i ' + \psi _i ^ 2$ for the function $\psi _i \in \mathcal A ( I _{\pi})$ given by 
$\psi_i = - \frac{1}{2}  ( (\log p_\eta)'  +(\log h_i ) ')$.  
By Lemma \ref{l:ml}  and Theorem  \ref{t:1dextra} (i), we have
$$
\mu _ 1 ( I _ {d_i}, h_i p_\eta) \geq  \frac{3 \pi ^ 2}{D_\Om^2}  + C \frac{(D_\Om-d_i) ^ 3}{D _\Om^ 5} \geq  \frac{3 \pi ^ 2}{D_\Om^2}  + C \frac{(D_\Om-D_\om) ^ 3}{D _\Om^ 5}  \,.
$$

By applying  Lemma \ref{l:preinside} to an eigenfunction for $\mu _ 1 (\om, p_\eta)$, the above inequality yields 
\eqref{f:diffd.1} for the weight $p _\eta$. On the same set $\om$, the convergence of $p _\eta$ to $p$ leads to 
\eqref{f:diffd.1}  for the weight $p$. We now prove \eqref{f:diffd.1}  for $\om$ and $p$ as in the statement. Let   
 $\omega _n$  be an increasing sequence of smooth open bounded convex sets, compactly contained into $\omega$, which converges to $\omega$ in Hausdorff distance as $n \to +\infty$. 
 Since we have already shown that   \eqref{f:diffd.1}   holds true for $\mu _ 1 (\om_n,  p)$, 
  to conclude our proof it is enough to observe that 
\begin{equation}\label{f:scs} \mu _ 1 (\om, p) \geq \limsup _n \mu _1 (\om_n, p) \,.
\end{equation} 
The argument to obtain the above inequality is similar  as the one used in the proof of Lemma \ref{l:ml}.
Namely, for every $\e >0$, we consider a function  $u _\e \in H ^ 1 _{\rm loc} (\om) \cap L ^ 2 (\om, p \, dx)$ with $\int _\om u _ \e p = 0$,  such that
$$\mu _ 1 (\om, p) \geq  \frac{ \int _{\om} |\nabla u_\e| ^ 2    p}{\int _{\om}| u_\e| ^ 2    p  }     - \e\,;$$ 
setting $c _{\e, n}:=  \frac{1}{|\om _n|} u _ \e p $, we have $\lim _n c _{\e , n } = 0$. Hence, taking
$u _{\e, n}:= u _\e - c _{\e, n}$ as a test function for $\mu _ 1 (\om _n, p )$, we obtain 
$$\mu _ 1 (\om, p) \geq  \frac{ \int _{\om} |\nabla u_\e| ^ 2    p}{\int _{\om}| u_\e| ^ 2    p  }     - \e
= \lim _n   \frac{ \int _{\om_n} |\nabla u_{\e, n}| ^ 2    p}{\int _{\om_n}| u_{\e, n}| ^ 2    p  }     - \e 
\geq \limsup _n   \mu _ 1 (\om _n, p  )    - \e 
\, .$$ 
The inequality \eqref{f:scs} follows by the arbitrariness of $\e >0$. 
\qed

\bigskip

 Relying  on  the previous result, we are ready to prove:

\begin{theorem} [Rigidity]  \label{t:rigidity} Let $N \geq 2$. For every open bounded convex domain $\Omega$ in $\R ^N$ of diameter 
 $D_\Om$,   we have
$$\lambda _ 2 (\Omega)- \lambda _ 1 (\Om) > \frac{3 \pi ^ 2}{  D_\Om^ 2  } \,.$$ 
\end{theorem}  
\bigskip

{\it Proof}.   For brevity, we are going to write $D$ in place of $D _\Om$.   Assume by contradiction that $\Omega\subset \R ^N$, $N \geq 2$,  is an open bounded convex domain    such that 
$ \mu _ 1 (\Om, u _ 1 ^ 2) = \frac{3 \pi ^ 2}{D ^ 2}$ 
We distinguish the case of dimension $N=2$ from the case of arbitrary dimension.

\medskip
{\it Case $N = 2$}.   
For every $n  \geq 1$, let 
$\mathcal P _n = \{ \Om _1 , \dots,   \Om _n\}$ be a $u _ 1 ^2$-weighted $L ^2$ equipartition of $\overline u$ in $\Om$,
being $\overline u:= \frac{u _ 2}{ u _ 1}$ an eigenfunction for $\mu _1 (\Om, u _ 1 ^ 2)$ (cf.\ 
Proposition \ref{p:yau}).  
  Denoting by $D_i$ the diameter of $\Om _i$, by Lemma   \ref{l:goodparti} and  Proposition \ref{p:inside} (applied with $p = u _ 1 ^ 2$), we have 
$$\frac{3 \pi ^ 2}{D^ 2}  = \mu _ 1 (\Om, u _ 1 ^ 2) \geq \frac{1}{n} \sum _{i = 1} ^n \mu _ 1 (\Om_i, u _ 1 ^ 2 ) \geq \frac{3 \pi ^ 2}{D^2}  + C  \sum _{i = 1} ^n
\frac{(D-D_i) ^ 3}{D ^ 5}  \,, $$
and $\overline u$ is an eigenfunction for  each
$\mu _ 1 (\Om_i, u_ 1 ^ 2)$. 
 
We infer that 
\begin{eqnarray}
& D = D _i \qquad \forall i = 1, \dots, n\, , & \label{f:eqdiam}
\\ 
& \mu _ 1 (\Om _i, u _ 1 ^ 2) =  \frac{3  \pi ^ 2}{D ^ 2}\qquad \forall i = 1, \dots, n\, . & \label{f:eqeig}
 \end{eqnarray}
Moreover, by the last assertion in Lemma \ref{l:goodparti},  $\overline u$ is an eigenfunction for  each
$\mu _ 1 (\Om_i, u_ 1 ^ 2)$. 

 We claim that $\Om$ is a circular sector (and hence $\Om _i$ are circular subsectors of the same diameter,  the supremum of whose opening angles is infinitesimal as $n \to + \infty$). Indeed,  consider  a straight line $r$ which determines  the $u _ 1 ^2$-weighted $L ^ 2$ equipartition $\mathcal P _2  = \{ \Om _1 ,  \Om _2\}$.
 Necessarily, $r \cap \Om$  must be a diametral segment of $\Omega$.  Indeed, let $I _ D$ be a fixed diametral segment of $\Om$: if $I_D \neq (r \cap \Om)$, either $I_D\cap  (r\cap \Om) = \emptyset$, or $I_D$ is transverse to $r  \cap \Om$. 
In the first case, assume with no loss of generality that $I_D$ is contained  into $\Om _1$. Then, since also  $\Om _2$ has diameter $D$,  the set $\Omega$ would contain two distinct diametral segments, which is not possible since a diagonal of the quadrilateral having vertices at the endpoints of the two diametral segments would have length strictly larger than $D$. In the second case, since $I_D$ is not entirely contained neither in $\Om _1$, nor in $\Om _2$, each among $\Om _1$ and $\Om _2$ would contain a segment of length $D$, and again this would  imply the existence of a segment of length strictly larger than $D$ in $\Om$.  By repeating the same argument for the successive straight lines which determine  the $u _ 1 ^2$-weighted $L ^ 2$ equipartition $\mathcal P _n = \{ \Om _1 , \dots,   \Om _n\}$
we infer that 
each cutting line is a diametral segment.  Since, as already noticed above, there cannot be  two disjoint diametral segments in $\Om$, we conclude that $\Om$ contains, for every $n \geq 1$, $n$ diametral segments having a common endpoint. By convexity, we infer that $\Om$ contains their convex envelope, and,  by passing to the limit as $n \to+ \infty$, we conclude that $\Om$ contains a circular sector. 
Actually, $\Om$ must coincide with such circular sector, because any point outside the sector cannot lie into
any cell of 
$\mathcal P _n$.

Once we know that $\Om$ is a circular sector and  that any cell of $\mathcal P _n$ a circular subsector, we  exploit the 
fact that $\overline u := \frac{u _ 2}{ u _ 1}$ is an eigenfunction for  each
$\mu _ 1 (\Om_i, u_ 1 ^ 2)$.  Since this holds true for every $n$, 
we infer that the function $\overline u$ satisfies, 
 on every radius of the circular sector $\Om$, 
the Neumann condition $\frac{\partial \overline u}{\partial \nu} = 0$, being $\nu$  the normal direction to the radius.  
This implies that the function $\overline u$ is radial on $\Omega$. But the expression of $u _ 2$ and $u _ 1$ (which  is explicitly known for
a circular sector, see for instance \cite[Section 3.2]{greb}) shows that  this is not the case. So we have reached a contradiction. 

\medskip
{\it Step 2} ($N$ arbitrary). We consider again 
 a $u _ 1 ^2$-weighted $L ^2$ equipartition $\mathcal P _n = \{ \Om _1 , \dots,   \Om _n\}$  of $\overline u$ in $\Om$.   By the same arguments 
 used in case $N = 2$,  all the cells of the partition satisfy  \eqref{f:eqdiam}-\eqref{f:eqeig}. Moreover, 
by arguing as in Remark \ref{r:PW}, we can assume that, for each each cell $\Om _i \in \mathcal P _n$, it holds (in a suitable orthogonal coordinates system associated with the cell) 
\begin{equation}\label{f:coord}\Om _i \subset \Big \{ ( y, x) \in \ \R ^ { N-2}\times \R ^ 2 \ : \ \  |y _j| \leq \e_n \ \forall j = 1, \dots, N-2;\  |x _j | \leq \frac{ D}{2}  \ \forall j = 1, 2    \}\,,\end{equation}
 with $\e _n \to 0$ as $n \to + \infty$. 
 This implies  that, in the limit as $n \to + \infty$, any sequence of cells  $\Om _{k (n) }  \in \mathcal P _n$ 
 converges, up to a subsequence, to a degenerate convex set, which may be 
either one-dimensional or two-dimensional. Any such limit set has diameter $D$, because
all the cells have diameter $D$. 
This allows to exclude that all the limit sets are one-dimensional. 
Otherwise,  by fixing two distinct segments $S'$ and $S''$ contained in $\Omega$, both parallel to $e _N$, and considering the sequences of domains
$\Om _{k' (n)} \in \mathcal P _n$ and $\Om _{k'' (n)} \in \mathcal P _n$  which contain them, 
in the limit as $n \to + \infty$ we would find two parallel
diametral segments  contained in $\Omega$. 
Let $\Om _{k (n )} \in \mathcal P _n$ be a sequence which converges to a two-dimensional convex set $\Om _0$. 
We claim that, for a suitable $(\frac{1}{N-2})$-concave function $h$, it holds 
\begin{equation}\label{f:equpianapre} 
\liminf _{n\to + \infty} \frac{ \int_{\Om _{k (n)}} |\nabla \overline u| ^ 2 u _ 1 ^ 2  }{ \int_{\Om _{k (n)}} |\overline u| ^ 2 u _ 1 ^ 2 }  \geq 
\frac{ \int_{\Om _0} |\nabla \overline u| ^ 2  h u _ 1 ^ 2   }{ \int_{\Om _{0} }|\overline u| ^ 2 h u _ 1 ^ 2   } \,.\end{equation} 
Indeed, in a coordinate system such that $\Om _{k ( n)}$ satisfies \eqref{f:coord}, with the change of variables 
  $T_n :  \R ^ {N-2} \times \R ^2  \to  \R ^ {N-2 } \times \R ^ 2$ defined by $T_n ( y, x) = (\e _n x', x )$, we have 
$$ 
\frac{ \int_{\Om _{k (n)}} |\nabla \overline u| ^ 2( y, x)  u _ 1 ^ 2 (y, x)  }{ \int_{\Om _{k (n)}} |\overline u (y, x)| ^ 2 u _ 1 ^ 2 ( y, x) } = 
\frac{ \int_{T_n ^ { -1} ( \Om _{k (n)})} |\nabla \overline u| ^ 2 (  \e_n   x', x)    u _ 1 ^ 2 ( \e_n   x', x)   }{ \int_{T_n ^ { -1} ( \Om _{k (n)})}|u| ^ 2 ( \e_n   x', x)  u _ 1 ^ 2  
(  \e_n   x', x)   } \,.$$ 
We now pass to the limit as $n \to + \infty$: by using Fatou's lemma and denoting by $h ( x)$ the $\mathcal H ^ {N-2}$-measure of the slices of the limit set of 
$T_n ^ { -1} ( \Om _{k (n)})$ at fixed $x = (x_1, x_2)$, we obtain 
$$
 \liminf _{n \to + \infty}  \int_{\Om _{k (n)}} |\nabla \overline u| ^ 2( y,x)  u _ 1 ^ 2 ( y, x)    \geq  \int_{\Om_0} |\nabla \overline u| ^ 2 ( 0, x)  h ( x)   u _ 1 ^ 2 ( 0, x)   \,;$$
on the other hand, recalling that $\overline u = \frac{u_2}{u _ 1}$,  by dominated convergence we get  
$$ 
\lim _{n \to + \infty}  \int_{T_n ^ { -1} ( \Om _{k (n)})}  \overline u ^ 2 ( \e_n   x', x)    u _ 1 ^ 2 ( \e_n   x', x)  =  \int_{\Om_0}|\overline u| ^ 2 (  0, x) h ( x)  u _ 1 ^ 2  (0, x)\,.
$$ 
Therefore, we have 
\begin{equation}\label{f:equpiana} 
\frac{3 \pi ^ 2}{ D ^ 2} = \liminf_n  \frac{ \int_{\Om _{k (n)}} |\nabla \overline u| ^ 2 u _ 1 ^ 2  }{ \int_{\Om _{k (n)}} |\overline u| ^ 2 u _ 1 ^ 2 }  \geq 
\frac{ \int_{\Om _0} |\nabla \overline u| ^ 2  h u _ 1 ^ 2   }{ \int_{\Om _{0} }|\overline u| ^ 2 h u _ 1 ^ 2   }  \geq \mu _ 1 (\Om_0,h  u _ 1 ^ 2  ) \geq  \frac{3 \pi ^ 2 }{D ^ 2 }\,, 
\end{equation}  
 where the first equality holds by \eqref{f:eqeig}, the second inequality holds  by \eqref{f:equpianapre}, the third inequality holds because $\overline u (x, 0)$ is an admissible test function for $\mu _1 (\Om _0, u_ 1 ^ 2 h)$, and the fourth inequality holds by Proposition \ref{p:inside} 
(applied with $\om = \Omega_0$ and  $p = h u _ 1 ^ 2$).

We now focus on the two-dimensional convex set $\Om _0$. 
Since   we know from \eqref{f:equpiana} that $\mu _ 1 (\Om_0, h u _ 1 ^ 2  ) =  \frac{3 \pi ^ 2 }{D ^ 2 }$, with eigenfunction $\overline u$, we can repeat the same arguments as in Step 1 of the proof,  
with the weight $u _ 1 ^ 2$ replaced by $h u _ 1 ^ 2 $. (More precisely, the arguments of Step 1 are repeated working with $(h u _ 1 ^ 2)$-weighted $L^2$ equipartitions of $\overline u$ in $\Om_0$,  and applying Proposition \ref{p:inside} with $p = h u _ 1 ^ 2$).

By this way, we obtain that $\Om_0$ is a circular sector
and, as well, all the cells of any  
 $(h u _ 1 ^ 2)$-weighted $L^2$ equipartition of $\overline u$ in $\Om_0$  are circular sectors $\Om_0 ^i$,  with 
$\mu _ 1 (\Om _0 ^ i , h u _ 1 ^ 2 ) = \frac{3 \pi ^ 2}{D^2}$.
 
Let us fix a radius in $\Om _0$  which does not belong to $\partial \Om _0$, and look at a sequence of  sectors $\Om_0 ^i$,  of infinitesimal opening angle,  having such radius
in their closure. 
On this radius, hereafter denote by $I _D$, in the limit of a large number of cells,  by arguing as in the proof of \eqref{f:equpiana},  we obtain 
$$ \mu _1  \big ( I _ D, p \big  ) =  \frac{3 \pi ^ 2}{D^2}\,, \qquad \text{ with } p(x):=  \big (x + \frac{D}{2} \big )   h u _ 1 ^ 2 \,.$$  

By Corollary \ref{c:sharpmu} we have, for some positive constant $k$,  
$p ( x) = k \cos ^2  \big  (\frac{\pi}{D} x \big )$.  
This implies in particular that
$p (x) \sim \big ( x + \frac{D}{2} \big ) ^ 2$  as $x \to - \frac{D}{2}$, 
which is not possible since 
$$ p( x) =  \big (x + \frac{D}{2} \big )   h u _ 1 ^ 2  = o   \big ( x + \frac{D}{2} \big ) ^ 2 \qquad \text{ as } x \to - \frac{D}{2}\,.$$ 
We have thus reached a contradiction and our proof is achieved.  \qed

\section{Estimate of the excess: the geometric play of cells in $N$ dimensions}\label{sec:construction}

In this section we prove Theorem \ref{t:quantitative}, first in dimension $N = 2$ and then in higher dimensions. 
 We keep the shortened notation $D$ in place of $D _\Om$.

\subsection{The case $N=2$}\label{s5.1}

We start by giving some preliminaries. 
For convenience, let us introduce a geometric quantity, related to the width,  which will be 
used throughout the proof. 
Once fixed a diameter  $I _D$   of $\Om$, we define the depth $\eta$ of  $\Omega$ with respect to  $I_D$  as
the maximum of the length of all sections orthogonal to it. 
We point out that, denoting by  $w^{\perp}$ the width of $\Om$ in direction orthogonal to  $I _D$, 
it holds $$
	\eta \leq w^{\perp } \leq 2\eta, \qquad 
	\frac{\eta D}{2}\leq|{\Omega}|\leq \eta  D, \qquad 
 |\Omega| \leq wD. 
$$ 
Moreover, the inequality $	\frac{w  D}{2}\leq  |{\Omega}|$ holds   provided the width is small with respect to the diameter, precisely  
$w\le \frac{\sqrt{3}}{2}  D$ (see \cite{kubota,scott}). In particular,  when the latter inequality is satisfied,  
the depth and the width are equivalent, in the sense that
\begin{equation}\label{f:eqetaw} 
\frac w 2 \leq \eta \leq 2 w \,.
\end{equation}  
 Next proposition ensures the existence of a dimensional constant,  related to the localization of $u_2$,   which will intervene in the proof of Theorem \ref{t:quantitative}.

\begin{proposition} \label{p:propsection} 
There exists a dimensional constant $\Lambda$ such that: if $\Omega \subset \R ^2$ is an  open bounded convex set with diameter $D$, 
width 
$w\le \frac{\sqrt{3}}{2}  D$, and Dirichlet eigenfunctions $u_1, u _2$, and $\omega$ is cell of a 
 $u _ 1 ^2$-weighted $L^2$ equipartition of $\overline u = \frac{u_2 }{u _1}$ in $\Om$ composed by $n$ cells, in an orthogonal coordinate system such that 
 a diameter of $\Om$ is the horizontal segment $\big (- \frac D 2, \frac D 2 \big ) \times \{ 0 \}$, the vertical sections 
 $\omega _{x_1}:= \om \cap (  \{ x _ 1\} \times \R  ) $ 
satisfy   
\begin{equation}\label{e:propsection}
\| \mathcal H ^ 1 (\om _ {x_1}) \| _{L^\infty ( - \frac D 2, \frac D 2 ) } \geq  \Lambda \frac{w}{n} \,.
\end{equation}

\end{proposition} 

\proof  Let us first prove the following claim: if $\Omega$ is an open bounded convex set of diameter $( - \frac D 2, \frac D 2 ) \times \R$ and depth $\eta$ with respect to such diameter, for any  open bounded convex set 
$\om\subseteq \Om$, setting  $\omega _{x_1}:= \om \cap ( \{ x _ 1 \} \times \R ) $, it holds
$$\int _\om | u | ^ 2  \leq \eta \, \| \mathcal H ^ 1 (\om _{x_1}) \| _{L^\infty( - \frac D 2, \frac D 2 ) } \int_\Om |\nabla u | ^ 2 \qquad \forall u \in H ^ 1 _0 (\Om) \,.$$ 
  It is not restrictive to prove the above inequality for $u \in \mathcal C ^ \infty _0 (\Om)$. Let $x = (x_1, x_2 ) \in \Om$, and set 
$\Omega _{x_1}:= \Om \cap ( \{ x _ 1\}  \times \R) $. We have $u ( x) = \int  _{- \infty} ^ { x _ 2} \frac{\partial u}{\partial x _ 2} (x_1, s ) \, ds$, and 
 $\| \mathcal H ^ 1 (\Om _{x_1}) \| _{L^\infty(0, D)} \leq \eta$, and hence, by H\"older inequality, 
$$|u (x)| ^ 2 \leq \eta \int _{\Om _{x_1}} \big  |\frac{\partial u}{\partial x _ 2} (x_1, s ) \big | ^ 2 \, ds \qquad \forall x = (x_1, x _2) \in \om \,.$$
Therefore, 
$$\begin{array}{ll}
\displaystyle \int _\om | u | ^ 2  & \displaystyle \leq \eta \int _\om   \Big [ \int _{\Om _{x_1}} \big  |\frac{\partial u}{\partial x _ 2} (x_1, s ) \big | ^ 2 \, ds   \Big ] \, dx_1 \, dx _ 2 
\\ \noalign{\medskip} 
& \displaystyle = \eta \int _\Om  \chi _{\om} ( x_1, x_2)   \Big [ \int _{\Om _{x_1}} \big  |\frac{\partial u}{\partial x _ 2} (x_1, s ) \big | ^ 2 \, ds   \Big ] \, dx_1 \, dx _ 2 
\\ \noalign{\medskip} 
& \displaystyle = \eta \int _{- \frac D 2 } ^ {\frac D 2}  \mathcal H ^ 1 (\om_{ x_1})   \Big [ \int _{\Om _{x_1}} \big  |\frac{\partial u}{\partial x _ 2} (x_1, s ) \big | ^ 2 \, ds   \Big ] \, dx_1  
\\ \noalign{\medskip} 
& \displaystyle \leq \eta \, \| \mathcal H ^ 1 (\om _{x_1}) \| _{L^\infty(- \frac D 2 , \frac D 2)}   \int _{\Om } \big  |\nabla u  | ^ 2\,. 
\end{array}$$ 

Now, we apply the claim just proved to the function $u _2$, assuming that $\omega$ is a cell of a 
 $u _ 1 ^2$-weighted $L^2$ equipartition of $\overline u$ in $\Om$. 
 Since 
 $$\int _\om  u _ 2 ^ 2 =  \int _\om  |\overline u | ^ 2  u _ 1 ^ 2  = \frac{1}{n} \int_ \Om  |\overline u | ^ 2  u _ 1 ^ 2    = 
  \frac{1}{n} \int_ \Om   u _2 ^ 2 = 
 \frac{1}{n} \,,  $$ 
 we obtain
 \begin{equation}\label{f:s1}
 \frac{1}{n} \leq \eta \, \| \mathcal H ^ 1 (\om _{x_1}) \| _{L^\infty(0, D)} \int_\Om |\nabla u _2| ^ 2 \leq  \eta \, \| \mathcal H ^ 1 (\om _{x_1}) \| _{L^\infty(0, D)}  \lambda _ 2 (\Om)\,.
 \end{equation}
Denoting by $\rho$ the inradius of $\Om$,  for a dimensional constant $\Lambda$ we have
 \begin{equation}\label{f:s2} 
  \lambda _ 2 (\Om) \leq \frac{ \Lambda}{\rho^ 2}\leq \frac{9 \Lambda}{w^ 2}  \,,
  \end{equation} 
where the first inequality follows from the fact that $\Omega$ contains a disk of radius $\rho$ (thanks to the monotonicity of $\lambda _ 2 (\cdot)$ under inclusions), and the second one from the  elementary
inequality $w \leq 3  \rho $.  
The conclusion follows by combining \eqref{f:s1} and \eqref{f:s2}, and recalling that $\eta$ and $w$ satisfy \eqref{f:eqetaw}
(thanks to the assumption $w\le \frac{\sqrt{3}}{2}  D$). 
 \qed 

\bigskip
\begin{remark}\label{r:propsection}
A consequence of Proposition \ref{p:propsection} is the following information on the mass distribution of the second eigenfunction. Assuming that the width of the set $\Omega$ is small and the diameter of the cell $\omega$ is large, we get from \eqref{e:propsection} that 
$$|\omega|\ge \Lambda' \frac{|\Omega|}{n} = \Lambda'  |\Omega| \int _\omega u_2^2,$$
so that   $u_2 ^2$   does not localize on $\omega$.  This means that, though in general   $u_2^2$  may localize (for instance on thinning triangles),  the concentration of mass cannot occur on cells with large diameter.
\end{remark}

\bigskip

{\it Proof of Theorem \ref{t:quantitative} in dimension $N = 2$}.  It is not restrictive to prove the statement under the hypotheses that $\Om$ is smooth, $D = \pi$, 
and a diameter of $\Om$ is the segment $\big ( - \frac \pi 2 , \frac \pi 2 \big ) \times \{0 \}$.  
Moreover, thanks to Theorem \ref{t:rigidity}, we can assume that 
$w\le  w_0 := \frac{\sqrt{3}}{2}  \pi$, so that \eqref{f:eqetaw}   holds.
Along the proof, the maximal admissible width will be diminished when necessary, and will be still denoted by $w _0$ (equivalently, we are going to 
indicate by $\eta _0$ a maximal admissible value for the depth $\eta$ with respect to the diameter we have fixed).

By Proposition \ref{p:yau}, we have $\lambda _ 2 (\Om)- \lambda _ 1 (\Om)  =  \mu _ 1 (\Om, u _ 1 ^ 2)$, with first eigenfunction  $\overline u = \frac{u _ 2}{u _ 1}$. 
For every $n$, we let $\{ \Om _1, \dots, \Om _n \}$ be a $u _ 1 ^ 2$-weighted  $L^2$ equipartition of $\overline u$ in $\Om$. We set $D _i:= {\rm diam} (\Om _i )$ and $h _i$ the profile function of $\Om _i$ in direction orthogonal to a fixed diameter of $\Om _i$. 

We let $c$ be some number in $(0, 1)$,
whose value will be diminished when necessary during the proof,  its final choice being postponed to the 
end of the proof.

For the benefit of the reader, before giving the detailed proof, we provide below the list of cases, along with a 
very short heuristic description for each of them.

\medskip
-- {\it Case 1}: 
For a sequence of integers $n$,  there exists $\mathcal I _n \subset  \{1, \dots, n\}$  with  ${\rm card} ( \mathcal I _n ) \geq \frac{n}{2}$, such that, $\forall i \in \mathcal I _n$,  the diameter of $\Om _i$ is ``small''. In this case, by applying 
  Proposition \ref{p:inside} to cells  $\Om _i$ for
$i \in \mathcal I _n$, we obtain the quantitative gap inequality.

\smallskip

-- {\it Case 2}: For a sequence of integers $n$, 
there exists $\mathcal I' _n \subset  \{1, \dots, n\}$  with  ${\rm card} ( \mathcal I' _n ) \geq \frac{n}{2}$, such that,  $\forall i \in \mathcal I' _n$,  
the diameter of $\Om _i$ ``large''. In this case we prove that the cells $\Om _i$ for $i \in \mathcal I ' _n$  must intersect the vertical walls of a strip based on the diameter of $\Om$, so that they
can be ordered vertically and they have a polygonal boundary inside the strip. 

\smallskip

-- {\it Case 2.1}: For a sequence of integers $n$,  there exists
$ \mathcal J _n \subset\mathcal  I' _n$, with  ${\rm card} ( \mathcal J_n ) \geq \frac{n}{4}$, such that,  $\forall i \in \mathcal J _n$,
$\Om _i$ has a vertex in  the strip. 
In this case we reduce ourselves back to  a similar situation as in Case 1, so we prove the quantitative inequality. 

\smallskip

-- {\it Case 2.2}: For a sequence of integers $n$, there exists
$ \mathcal  J' _n \subset \mathcal  I'_n$, with ${\rm card} ( \mathcal  J' _n ) \geq \frac{n}{4}$, such that, $\forall i \in \mathcal J' _n$,  
$\Om _i $  has no vertex in the strip.
In this case we prove that,  $\forall i \in \mathcal  J'_n$,   the 
 function $h _i$ representing the profile of $\Om _i$ in direction orthogonal to its diameter is affine. 
 Then we distinguish two further sub-cases.

\smallskip

-- {\it Case 2.2.1}: For 
 a sequence of integers $n$,  there exists 
$  \mathcal  S _n\subset J'_n$,  with  ${\rm card} ( \mathcal  S_n ) \geq \frac{n}{8}$, such that, $\forall i \in \mathcal S _n$, the  Rayleigh  quotient 
of $\overline u$ on $\Om _i$ with respect to the measure $u _ 1 ^ 2 \, dx$ is ``large''.  In this case we obtain the quantitative inequality.

\smallskip

-- {\it Case 2.2.2}: 
For 
 a sequence of integers $n$, there exists 
$  \mathcal  S _n' \subset J'_n$  with ${\rm card} (\mathcal   S'_n ) \geq \frac{n}{8 }$, such that,  $\forall i \in \mathcal S'_n$, 
the above mentioned  Rayleigh  quotient is ``small''. 
  In this case we prove that, for $i$ in a subfamily $S _n'' \subset \mathcal  S'_n$,  with ${\rm card} ( \mathcal  S'' _n ) \geq \frac{n}{16} $, 
an extra-term  $\delta _i>0$, proportional to $(1 -  \frac{h _ i ^ {max} }{ h _i ^ {min}})^2$, can be added in our lower bound for the Neumann eigenvalue
  of the diameter of $\Om _i$
  with weight $h _ i u _ 1 ^ 2$.

\smallskip

--  {\it Case 2.2.2.1}: For a sequence of integers $n$, 
there exists  $\mathcal  Z _n \subset \mathcal  S'' _n$,  with  ${\rm card} ( \mathcal  Z_n ) \geq \frac{n}{32}$, such that,  $\forall i \in \mathcal Z _n$, 
$\delta _i$ is ``large''.   In this case we obtain the quantitative  inequality. 

\smallskip

--  {\it Case 2.2.2.2}: For a sequence of integers $n$, 
there exists $\mathcal Z' _n \subset\mathcal   S''_n$,  with  ${\rm card} ( \mathcal   Z'_n ) \geq \frac{n}{32}$, such that, $\forall i \in \mathcal Z'_n$,  
$\delta _i$ is ``small''.  
Loosely speaking, this enables us to minorate the length of the two segments 
of intersection between any of these  cells and the vertical walls of the strip, implying that our pile is sufficiently high. And since it is composed by cells with ``large'' diameter, 
 by the Pythagorean Theorem  we eventually find in $\Om$ a segment larger than its diameter, reaching a contradiction.

\medskip
We now detail each of the cases above.

\bigskip
{\it Case 1.} 
For a sequence of integers $n$,  there exists $\mathcal I _n \subset  \{1, \dots, n\}$  with  ${\rm card} ( \mathcal I _n ) \geq \frac{n}{2}$, such that, 
$$D_i \leq \sqrt {\pi ^ 2 - c \eta ^ 2 } \qquad \forall i \in  \mathcal I   _n \subset \{1, \dots, n\}\,. $$
When applying Proposition \ref{p:inside} to cells $\Om _i$ for $i \in  \mathcal I _n$,  the extra-term at the right hand side of inequality
\eqref{f:diffd}  
 admits the following lower bound \begin{equation}\label{f:stimetta}
 \frac{(D- D_ i ) ^ 3 }{D ^ 5}   \geq  \frac {c^ 3}{8 \pi ^ 8} \eta ^ 6 \,.
 \end{equation}
Therefore  we can prove the validity of the inequality \eqref{f:qgap},  
with $\overline c =  C \frac{ c ^ 3}{2 ^ {10} \pi ^ 8 }$, being $C$ the  absolute constant appearing in Proposition \ref{p:inside}. 
Indeed, we have 
$$\begin{array}{ll}
\mu _ 1 (\Om, u _ 1 ^ 2 )  & \displaystyle \geq  \frac{1}{n} \sum _{i = 1} ^ n \mu _ 1 (\Om _i, u _ 1 ^ 2 )  
   =  \displaystyle \frac{1}{n}  \Big [ \sum _{i \not \in \mathcal  I _n}  \mu _ 1 (\Om _i, u _ 1 ^ 2 )  +  \sum _{i \in  \mathcal I _n}  \mu _ 1 (\Om _i, u _ 1 ^ 2 )  \Big ] 
\\  \noalign{\medskip}
&   \geq  \displaystyle \frac{1}{n}  \Big [ 3  \, {\rm card} (  \mathcal I _n^c )   +  \Big ( 3  + C  \frac {c^ 3}{8 \pi ^ 8 } \eta ^ 6 \Big )  {\rm card} ( \mathcal  I_n )   \Big ] 
\\  \noalign{\medskip}
& \geq 3   +  C \frac{ c ^ 3}{16 \pi ^ 8} \eta ^ 6 \geq 3   +  C \frac{ c ^ 3}{2 ^ {10} \pi ^ 8 } w ^ 6  \,,
\end{array}
$$
where the first inequality follows from Lemma \ref{l:goodparti}, the second one from Proposition  \ref{p:inside} (taking into account \eqref{f:stimetta}), and  the last one from the assumption ${\rm card} ( \mathcal I _n ) \geq \frac{n}{2}$ and \eqref{f:eqetaw}.

\bigskip
{\it Case 2}: For a sequence of integers $n$, 
there exists $\mathcal I' _n \subset  \{1, \dots, n\}$  with  ${\rm card} ( \mathcal I' _n ) \geq \frac{n}{2}$, such that, 
$$D_i \geq \sqrt {\pi ^ 2 - c \eta ^ 2 } \qquad \forall i \in     {\mathcal I }' _n\,.$$
Let $a$ be a parameter such that  $6\eta_0^2\leq a<\frac  \pi 4$.  We observe that for every $i \in \mathcal  I'_n$,  $\Om _i$ must intersect the vertical lines 
$ \{x_1= - \frac{\pi-a} {2} \}$ and $\{x_1=\frac{\pi-a} {2} \ \}$.  Indeed, assume by contradiction this is not true.
Since by the choice of $a$ it holds $\pi a - \frac{a ^ 2}{4} > 5 \eta _0 ^ 2$,  
the length of the projection of $\Om_i$ onto 
the vertical axis would be bounded from below by 
$$\sqrt { D _ i ^ 2 - ( \pi -\frac a2 ) ^ 2}  \geq \sqrt{ \pi ^ 2 - c \eta ^ 2  -\pi ^ 2 + \pi a - \frac {a^2}{ 4} } > \sqrt{ - c \eta _ 0 ^ 2 + 5 \eta _0 ^ 2 }  
\geq 2 \eta _0   \geq 
w ^ \perp \,,$$ 
yielding a contradiction.

As a consequence of this fact,  for  $i \in \mathcal  I'_n$,  the 
cells $\Om _i$ can be ordered in a vertical way,  and  inside the strip $\big ( - \frac {\pi-a} 2  , \frac {\pi-a} 2 \big ) \times \R$ 
none of their boundaries  (except for the bottom and the top cell) 
can contain portions of   $\partial \Om$.
Taking into account that $\Om _i$ are convex sets  obtained cutting by lines 
we infer that,  for every $i \in \mathcal  I'_n$ exception made for two indices,  the intersection  of $\overline \Om _i$ with the strip  $\big [ - \frac {\pi-a} 2  , \frac {\pi-a} 2 \big ]  \times \R$ is a polygon. 
	 
We proceed to analyse Cases 2.1 and 2.2. 	 
 
 \bigskip
  {\it Case 2.1}: For a sequence of integers $n$,  there exists
$ \mathcal J _n \subset\mathcal  I' _n$, with  ${\rm card} ( \mathcal J_n ) \geq \frac{n}{4}$, such that
$$\Om _i \text{ has a vertex in  }\big ( - \frac {\pi-2a} 2  , \frac {\pi-2a} 2 \big ) \times \R  
\quad \forall i \in \mathcal J _n.$$ 
  For $i \in  \mathcal J _n$,  denoting by $V_i$ one of its vertices  inside the strip $\big ( - \frac {\pi-2a} 2  , \frac {\pi-2a} 2 \big )  \times \R$
 (meant as   a vertex of the polygon $\overline 
 \Om _i \cap  ( \big [ - \frac {\pi-2a} 2  , \frac {\pi-2a} 2 \big ]  \times \R )$), there exists  
 a neighbouring cell $\Om _j$ such that $V _i \in \partial \Om _i \cap \partial \Om _j$. 
The diameter $D_j$ of  such 
$\Om _j$ satisfies 
  $$D_j \leq \sqrt{ (\pi-a) ^ 2 + 4 \eta ^ 2 } \leq \sqrt { \pi ^ 2 - 26 \eta_0 ^ 2 } \,,$$ 
 where the last inequality holds because we are assuming that $6\eta_0^2\leq a<\frac  \pi 4$.  
 Since the neighbouring cell $\Om _j$ can touch at most another cell  $\widetilde \Om _i$ with $i \in  \mathcal J _n$, we conclude that, in Case 2.1, for 
a sequence of integers $n$ 
  there exists   $\mathcal W _n \subset \{1, \dots, n\}$,  with ${\rm card} (  \mathcal W _n ) \geq \frac{n}{8}$, 
such that
 $$D_j  \leq  \sqrt { \pi ^ 2 - 26 \eta_0 ^ 2 }  \qquad \forall   j   \in  \mathcal W _n\,.$$
Then, by arguing 
as in Case 1,  we can prove \eqref{f:qgap}   (for a suitable  dimensional constant $\overline c$).

 \bigskip
 {\it Case 2.2}: For a sequence of integers $n$, there exists
$ \mathcal  J' _n \subset \mathcal  I'_n$, with ${\rm card} ( \mathcal  J' _n ) \geq \frac{n}{4}$, such that 
$$\Om _i \text{ has no vertex in  }\big ( - \frac {\pi-2a} 2  , \frac {\pi-2a} 2 \big ) \times \R \quad \forall i \in{\mathcal  J' _n }.$$ 
 For  $i \in  \mathcal J'_n$, we fix a diameter of $\Om _i$, and 
we observe that the angle $\alpha _i$ it forms with the  horizontal axis (i.e., with the diameter of $\Om$) 
does not exceed $\arcsin (\frac {2\eta}{D_i})\le \arcsin (\frac {2\eta}{\sqrt{\pi ^ 2- c\eta^2}})$. 
Then,  if we work in a new local coordinate system in which the fixed diameter of $\Om _i$ is the segment $\big (- \frac{ D_i}{2}, \frac{D_i}{2} \big ) \times \{0\}$, 
and we denote by $h_i$ the profile function  of $\Om _i$ in vertical direction, 
the function $h_i$  is necessarily  affine on $\big (- \frac{ D_i-4a}{2}, \frac{D_i-4a}{2} \big ) \times \{0\}$, as soon as 
		\begin{equation*}
			\ell_{\eta, a}:=\frac{a}{\cos (\arcsin   (\frac {2\eta}{\sqrt{\pi ^ 2- c\eta^2}}))}+ \eta   \frac {2\eta}{\sqrt{\pi ^ 2- c\eta^2}}\le 2 a
		\end{equation*}

 We remark that the above condition is satisfied provided $\eta _0$ is small enough.

We proceed to analyse Cases 2.2.1 and 2.2.2. 	 

 \bigskip
 {\it Case 2.2.1}: For 
 a sequence of integers $n$,  there exists 
$  \mathcal  S _n\subset J'_n   \text{ with } {\rm card} ( \mathcal  S_n ) \geq \frac{n}{8}$, such that  
 $$\displaystyle \frac{ \int _{\Om _i}  |\nabla \overline u| ^ 2    u _ 1 ^ 2 }{\int _ {\Om _i}  | \overline u| ^ 2    u _ 1 ^ 2   }      \geq 4 \quad \forall i \in \mathcal  S _n.$$ 
  In this case we can easily conclude by arguing in a similar way as done in Case 1. Indeed we have
 $$\begin{array}{ll}
\mu _ 1 (\Om, u _ 1 ^ 2 )  & \displaystyle  = \frac{1}{n} \sum _{i = 1} ^ n  \frac{ \int _{\Om _i}  |\nabla \overline u| ^ 2    u _ 1 ^ 2 }{\int _ {\Om _i}  | \overline u| ^ 2    u _ 1 ^ 2   } 
\geq  \displaystyle \frac{1}{n}  \Big [ \sum _{i \not \in \mathcal  S _n}  \mu _ 1 (\Om _i, u _ 1 ^ 2 )  +  \sum _{i \in  \mathcal S _n}   \frac{ \int _{\Om _i}  |\nabla \overline u| ^ 2    u _ 1 ^ 2 }{\int _ {\Om _i}  | \overline u| ^ 2    u _ 1 ^ 2   }   \Big ] 
\\  \noalign{\medskip}
&   \geq  \displaystyle \frac{1}{n}  \Big [ 3  \, {\rm card} (  \mathcal S _n^c )   +  4 {\rm card} ( \mathcal  S_n )   \Big ] 
\geq 3   +   \frac{1}{8} \,.
\end{array}
$$
The inequality \eqref{f:qgap} follows,  for a  dimensional constant $\overline c$, provided $w_0$ is small enough. 
 \bigskip

 \item[--]{\it Case 2.2.2}: For 
 a sequence of integers $n$, there exists 
$  \mathcal  S _n' \subset J'_n$  with ${\rm card} (\mathcal   S'_n ) \geq \frac{n}{8 }$, such that 
$$\displaystyle \frac{ \int _{\Om _i}  |\nabla \overline u| ^ 2    u _ 1 ^ 2 }{\int _ {\Om _i}  | \overline u| ^ 2    u _ 1 ^ 2   }     \leq 4 \quad \forall i \in \mathcal  S _n' .$$ 
We observe that  there exists $\mathcal   S _n'' \subset \mathcal  S'_n$,  with  ${\rm card} ( \mathcal  S'' _n ) \geq \frac{n}{16}$, such that  
$$|\Om_i| \leq \frac{17}{n} |\Om| \quad \forall i \in\mathcal   S _n''\,.$$  
 Indeed, otherwise it would be $|\Om_i| \geq \frac{17}{n}$ for at least  $\frac{n}{16}$ cells in  $\mathcal  S'_n $, and the union of such cells
 would have measure strictly larger than the measure of $\Omega$.  Let us now prove that, for every $i \in \mathcal S'' _n$, 
denoting by $K$ the absolute constant appearing in Theorem \ref{t:1dextra} (ii), and by $h _i ^ {max}$ and $h _ i ^ {min}$  the maximum and minimum of $h _ i$ on  $I_{D_i - 4a}$,
 it holds
 \begin{equation}\label{f:suppl}
 \mu _ 1 ( I _ { D_i} , h _ i u _ 1 ^ 2) \geq3+ \frac{8K}{ \pi ^2}  \Big (1 -  \frac{h _ i ^ {max} }{ h _i ^ {min}}  \Big ) ^ 2  \,.
 \end{equation}
 To that aim we apply Lemma \ref{l:stimaPW} with $\om_0 = \Om _i$, $p = u _ 1 ^ 2$, and $v= \overline u$.  Notice that this is possible because we are assuming that $\Om$ is smooth, so that $p$ is uniformly continuous  in $\Om _i$, and 
 $\overline u \in W ^ { 2, \infty} (\Om _i)$.  We obtain
\begin{equation}\label{f:formula}
 \begin{array}{ll}   \displaystyle \frac{ \int _{\Om _i}  |\nabla \overline u| ^ 2    u _ 1 ^ 2 }{\int _ {\Om _i}  | \overline u| ^ 2    u _ 1 ^ 2   }      &  
 \displaystyle \geq   \Big \{ \mu _ 1 ( I _ {D_i}, h_i u _ 1 ^ 2) - \frac{ \alpha (\e)  |\Om_i|   }{\int _{\Om _i}  \overline u ^2   u _ 1 ^ 2  } \Big  [ 1 + \mu _ 1 ( I _ {D_i  } , h_i  u _ 1 ^ 2)  \Big ( 1 +  \beta (\e) |\Om_i|   \Big )  \Big ] \Big \} 
 \\  \noalign{\bigskip} 
&  \displaystyle \geq   \Big \{ \mu _ 1 ( I _ {D_i}, h_i u _ 1 ^ 2) - \frac{  \alpha (\e)  |\Om_i|  }{\frac{1}{n}   } \Big  [ 1 +  \frac{3}{2} \mu _ 1 ( I _ {D_i  } , h_i  u _ 1 ^ 2)  \Big ] \Big \} \,, 
 \end{array} 
\end{equation}  
 where the last inequality holds   provided  $\e$ is so small  that  $\big ( 1 + \beta (\e)  |\Om|  \big ) \leq \frac{3}{2}$
 (which is true as soon as $n$ is sufficiently large). 

 Taking into account that  cells $\Om_i$ for $i \in \mathcal S''_n$  satisfy the condition $|\Om_i| \leq \frac{17}{n} |\Om|$,   the above inequality implies
$$ \displaystyle \frac{ \int _{\Om _i}  |\nabla \overline u| ^ 2    u _ 1 ^ 2 }{\int _ {\Om _i}  | \overline u| ^ 2    u _ 1 ^ 2   }        \geq   
  \Big \{ \mu _ 1 ( I _ {D_i}, h_i u _ 1 ^ 2)  \Big [   1 -   \frac{3}{2}  \cdot 17     |\Om|  \cdot \alpha (\e)   \Big ]   - 
    17   |\Om| \cdot   \alpha (\e) \Big \} \,.
$$
 We infer that, again for $\e$ sufficiently small, 
  it holds 
 $$ \displaystyle \frac{ \int _{\Om _i}  |\nabla \overline u| ^ 2    u _ 1 ^ 2 }{\int _ {\Om _i}  | \overline u| ^ 2    u _ 1 ^ 2   }   \geq   
  \frac{6}{7} \mu _ 1 ( I _ {D_i}, h_i u _ 1 ^ 2) -2  
 \,.  $$     
Recalling that,  for  $i \in \mathcal S''_n$,   the left hand side of the above inequality does not exceed $4$, we infer that  
\begin{equation}\label{f:dieci} 
\mu _ 1 ( I _ {D_i}, h_i u _ 1 ^ 2) \leq 7   \,.
\end{equation}  This enables us to apply Theorem \ref{t:1dextra} (ii), and conclude that \eqref{f:suppl} holds. 

We proceed to analyse Cases 2.2.2.1 and 2.2.2.2.	 
 
\bigskip
  {\it Case 2.2.2.1}: For a sequence of integers $n$, 
there exists  $\mathcal  Z _n \subset \mathcal  S'' _n$,  with  ${\rm card} ( \mathcal  Z_n ) \geq \frac{n}{32}$, such that  
$$ \frac{8K}{ \pi ^2} \Big ( 1- \frac{ h _ i ^ {min}  }{ h _ i  ^ {max} }  \Big ) ^2 \geq c \eta ^ 2  \qquad \forall i \in \mathcal  Z _n .$$ By \eqref{f:suppl} and \eqref{f:dieci},  for $i \in \mathcal Z_n$ it holds 
 \begin{equation*}
 \mu _ 1 ( I _ { D_i} , h _ i u _ 1 ^ 2) \geq3+  c \eta ^ 2  \quad \text{ and }  \qquad \mu _ 1 ( I _ { D_i} , h _ i u _ 1 ^ 2) \leq  7  \,. 
 \end{equation*}
Then, by \eqref{f:formula} we obtain 
$$ \displaystyle \frac{ \int _{\Om _i}  |\nabla \overline u| ^ 2    u _ 1 ^ 2 }{\int _ {\Om _i}  | \overline u| ^ 2    u _ 1 ^ 2   }      \geq    3+  c \eta ^ 2    -  \frac{ 12 \alpha (\e)  |\Om_i|   }{\frac{1}{n}   }   \,. $$ 
Hence, 
 $$\begin{array}{ll}
\mu _ 1 (\Om, u _ 1 ^ 2 )  & \displaystyle 
= \frac{1}{n} \sum _{i = 1} ^ n  \frac{ \int _{\Om _i}  |\nabla \overline u| ^ 2    u _ 1 ^ 2 }{\int _ {\Om _i}  | \overline u| ^ 2    u _ 1 ^ 2   } 
\geq   \displaystyle \frac{1}{n}  \Big [ \sum _{i \not \in \mathcal  Z _n}  \mu _ 1 (\Om _i, u _ 1 ^ 2 )  +  \sum _{i \in  \mathcal Z _n}   \frac{ \int _{\Om _i}  |\nabla \overline u| ^ 2    u _ 1 ^ 2 }{\int _ {\Om _i}  | \overline u| ^ 2    u _ 1 ^ 2   }    \Big ] 
\\  \noalign{\medskip}
&   \geq  \displaystyle \frac{1}{n}  \Big [ 3  \, {\rm card} (  \mathcal Z _n^c )   +   (3+  c \eta ^ 2 )   {\rm card} ( \mathcal  Z_n )   \Big ]  -  12 \alpha (\e)  |\Om | 
\\  \noalign{\medskip}
& 
\geq 3   +   \frac{c}{32} \eta ^ 2 - 12 \alpha (\e) |\Om |   \,,
\end{array}
$$
 and the conclusion follows  (as usual, for a suitable choice of $\overline c$).

\bigskip
 {\it Case 2.2.2.2}: For a sequence of integers $n$, 
there exists $\mathcal Z' _n \subset\mathcal   S''_n$,  with  ${\rm card} ( \mathcal   Z'_n ) \geq \frac{n}{32}$, such that 
$$ \frac{8K}{ \pi ^2}     \Big (1- \frac{ h _ i ^ {min}  }{ h _ i  ^ {max} }  \Big ) ^2 \leq c \eta ^ 2  \qquad \forall i \in \mathcal  Z' _n\,.$$
We shall now prove that, 
provided $\eta_0$ and  $a$ are small enough, and $c$ is well-chosen (as well,  small enough), 
this case cannot occur.

The contradiction argument relies on the following claim:  
denoting by 
$L _ i ^ {min} \leq L _ i  ^ {max}$ the lengths  of the intersections of $\Om _i$ with the lines
$ \{x_1= - \frac{\pi-4a} {2} \}$ and $\{x_1=\frac{\pi-4a} {2} \}$,   we have
\begin{equation}\label{f:sezlunga} 
\text{ for } \eta_0 , a \ll 1 \,, \quad L _ i ^ {min}  \geq \frac{\Lambda}{8} \frac{w}{n} \qquad \forall i \in \mathcal Z'_n\,,
\end{equation}  
where $\Lambda$ is the dimensional constant appearing in Proposition  \ref{p:propsection}.  

We first prove \eqref{f:sezlunga} and then we show how it leads to a contradiction. 

For $i  \in \mathcal Z' _n$,  choosing $c \eta _0 ^ 2 \leq \frac{2K}{ \pi ^ 2}$, it holds
$$\Big (1- \frac{ h _ i ^ {min}  }{ h _ i  ^ {max} }  \Big ) ^2 \leq  \frac{ \pi ^2 c}{8K}   \eta ^ 2    \leq
 \frac{1}{4}\,, $$ 
so that
$$\frac{ h _ i ^ {min}  }{ h _ i  ^ {max} }  \geq \frac{1}{2} \,.$$

The above inequality, combined with the area inequality $\frac{1}{2} h _ {min} D _ i  \leq \eta \pi$, implies  that  $h _{max} \leq \frac{4 \eta \pi}{D_i}$. 
In turn this implies that, 
if $\alpha _i$ is the angle already considered in Case 2.2 formed between the diameters of $\Om$ and $\Om _i$, it holds 
$$h _{max}\cdot \sin \alpha _i  \leq \frac{4 \eta \pi }{\sqrt {\pi ^ 2 - c \eta ^ 2 }}  \frac{2 \eta} {\sqrt {\pi ^ 2 - c \eta ^ 2}}\,. $$
Hence, for $\eta_0$ sufficiently small, we have $h _{max}\cdot \sin \alpha _i \leq a$, which ensures that the two  segments of lengths $L _ i ^ {min}$ and $L _i ^ {max}$ are interior to the trapeze with bases given by  the two segments of lengths 
 $h _ i ^ {min}$ and $h _i ^ {max}$ and oblique sides given by portions of $\partial \Om _i$. 
 This implies via Thales Theorem that 
$$\frac{ L _ i ^ {min}  }{ L _ i  ^ {max} }   
\geq \frac{ h _ i ^ {min}  }{ h _ i  ^ {max} }   
 \geq \frac{1}{2} \, .$$

Denoting by $({\Om}_i ) _ {x_1}$ the intersection of $\Om _i$ with the straight line $\{ x _1 \} \times \R$, we have 
\begin{equation}\label{f:lminlmax} L _ i ^ {min} \geq  \frac{1}{  2 } L _ i ^ {max} \geq \frac{1}{  4   }  \| \mathcal H ^ 1 (({\Om}_i ) _ {x_1}) \| _{L^\infty(I _ \pi)} \geq  \frac{\Lambda}{  4  } \frac{w}{n}\,, \end{equation}   

where the second inequality is satisfied provided $a$ is small  enough, 
 and the third one holds by
Proposition \ref{p:propsection}. This completes the proof of claim \eqref{f:sezlunga}.

Eventually, let us show how \eqref{f:sezlunga} yields a contradiction.  
Let $\Om _i ^ {b}$  and $\Om _i ^ { t}$ be the bottom and the top cell 
in the family of cells $\Om _i$, for $i \in \mathcal  Z'_n$. 
Let $A_1B_1$ and $A_2 B_2$ be two diametral segments  respectively for $\Om _i ^ {b}$  and $\Om _i ^ { t}$, and 
consider the quadrilateral $A_1 B _ 1 A_2 B _ 2$.

We can estimate from below the lengths of the segments $A_1A_2$ and $B _ 1 B _2$  by the lengths of their orthogonal projections on the vertical lines 
$ \{x_1= - \frac{\pi-4a} {2} \}$ and $\{x_1=\frac{\pi-4a} {2} \}$. Applying claim \eqref{f:sezlunga} to all the cells $\Om _i$ for $i \in \mathcal Z'_n$ 
(exception made for  $\Om _i ^ {b}$  and $\Om _i ^ { t}$), and recalling that ${\rm card} ( \mathcal   Z'_n ) \geq \frac{n}{32}$ we get 
$$\min \Big \{ \overline{A_ 1 A_2},  \overline{B_1 B_2} \Big \}  \geq \frac{\Lambda}{  4  } \frac{w}{n}  \Big ( \frac{n}{32} - 2 \Big ) \geq \frac{\Lambda}{  256  } w\,,  $$ 
where the last inequality holds for $n$ large enough. 
On the other hand, since 
$\mathcal Z'_n \subseteq \mathcal S''_n  \subseteq \mathcal S' _n  \subseteq \mathcal J' _n \subseteq \mathcal I ' _n$, it holds
$$\min \Big \{ \overline{A_ 1 B_1},  \overline{A_2 B_2} \Big \}  \geq  \sqrt { \pi ^ 2 - c \eta ^ 2 }\,.$$ 

Then, at least one of the two diagonals of the quadrilateral $A_1 B _ 1 A_2 B _ 2$ turns out to be  larger than the diameter of $\Omega$, yielding a contradiction. Namely, 
since at least one of the inner angles of the quadrilateral is larger than or equal to $\frac{\pi}{2}$, 
we have
$$ \max \Big \{ \overline{A_ 1 B_2},  \overline{A_2 B_1} \Big \}   \geq \Big [ (\pi ^ 2 - c \eta ^ 2  ) + \Big ( \frac{\Lambda}{ 256 } \Big ) ^ 2 w   ^2  \Big ] ^ { \frac{1}{2}}  > \pi \,$$ 
where the last inequality follows by choosing $c$ small enough.  
Having reached a contradiction, our proof is achieved. 
\qed

\bigskip
\subsection{The case $N \geq 3$}

The main tool to prove Theorem \ref{t:quantitative} for a domain $\Omega$ in dimension $N \geq 3$ is the following  quantitative estimate of a weighted two-dimensional Rayleigh quotient of $\overline u$ on  suitable planar sections of $\Omega$.

\begin{theorem}\label{t:quantitative2}  Let $N \geq 3$. 
There exists a dimensional constant  $\overline c>0$ such that: 

\smallskip 
-- for every open bounded convex domain $\Omega$ in $\R ^N$ of diameter $D$ and width $w$, 
with Dirichlet eigenfunctions $u _1, u _ 2$, 

--  for every planar subset $U$ of $\Om$ which, in a suitable orthogonal coordinate system $\{ e _1, \dots, e _N \}$ with  $e_N$ in the direction of the width of $\Om$, can be written as 
$$U = \Big \{ x \in \Om \ :\ x_i  = 0 \ \forall i = 1, \dots, N-2 \,,\ a \leq x _{N-1} \leq b \Big \}\,,$$

-- for every function $h: U \to (0, 1]$ independent of  $x _{N}$ and $\big (\frac{1}{N-2} \big)$-concave, such that
\begin{equation}\label{f:33}
\int_U  | \nabla _U u _ 2  |  ^ 2 h \leq 2 \lambda _ 2 (\Om) \int _ U| u_2  | ^ 2  h\,, 
\end{equation}

\smallskip
if $\overline u := \frac{u _2}{u _ 1} $  satisfies $\int _U \overline u   h u _ 1 ^ 2 = 0$, and $u _ 2$ is not identically zero on $U$, it holds 
$$
\frac{\int_U |\nabla_U \overline u| ^2 h u _ 1 ^ 2 } {\int_U | \overline u| ^2 h u _ 1 ^ 2   } \geq \frac{3 \pi ^ 2}{D^ 2}  + \overline c \frac{w^6}{D^ 8} \,.
$$

\end{theorem} 

 The proof of Theorem \ref{t:quantitative2} demands two preliminary ingredients. The former is the following inequality, 
holding by Proposition \ref{p:inside}:  
there exists an absolute constant $C$ such that, 
if all the assumptions of Theorem \ref{t:quantitative2} are satisfied, and
$\om$ is  a cell of a 
 $(h u _ 1 ^2)$-weighted $L^2$  equipartition of $\overline u$ in $U$, of diameter $d$, it holds
\begin{equation}\label{f:diffd2} 
\frac{\int_\om |\nabla_U \overline u| ^2 h u _ 1 ^ 2 } {\int_\om  | \overline u| ^2 h u _ 1 ^ 2   } \geq   \frac{3 \pi ^ 2}{D^ 2}  +  
 C \frac{(D-d) ^ 3}{D ^ 5} \,.
\end{equation}

 The latter is the following variant of Proposition \ref{p:propsection}.   
\begin{proposition}\label{p:propsection2} 
There exists a  dimensional constant $\Lambda$  such that: if all the assumptions of Theorem \ref{t:quantitative2} are satisfied, and $\omega$ is a cell of a 
 $(h u _ 1 ^2)$-weighted $L^2$  equipartition of $\overline u$ in $U$ composed by $n$ cells, setting 
 $\omega _{x_{N-1}}:= \om \cap ( \{ x _ {N-1} \}    \times \R e_N) $,  
 it holds 
   $$\| \mathcal H ^ 1 (\om _ {x_{N-1}}) \| _{L^\infty(a, b)} \geq  \Lambda \frac{w}{n} \,.$$  
\end{proposition}

\proof 
We consider the restriction of $u _2$ to $U$, and we argue in a similar way as in the proof of Proposition \ref{p:propsection}.  
Denoting by $\eta _U$
the depth of $U$ in direction $e_N$,   
 we have
$u_2 ( x_{N-1}, x_N) = \int  _{- \infty} ^ { x _ N} \frac{\partial u}{\partial x _ N} (x_{N-1}, s ) \, ds$, and 
 $\| \mathcal H ^ 1 ( \om  _{x_{N-1}}) \| _{L^\infty(a, b)} \leq \eta_U$, and hence
$$|u_2 ( x_{N-1}, x_N)| ^ 2 \leq \eta_U \int _{U _{x_{N-1}}} \big  |\frac{\partial u_2}{\partial x _ N} (x_{N-1}, s ) \big | ^ 2 \, ds \qquad \forall (x_{N-1}, x _N) \in U \,.$$
Since $h$ is independent of $x_N$, namely $h = h ( x_{N-1})$,  multiplying the above inequality by $h$ we obtain
the following inequality holding for every  $(x_{N-1}, x _N) \in U$: 
\begin{equation}\label{f:daintegrare} 
|u_2 ( x_{N-1}, x_N)| ^ 2 h ( x_{N-1}) \leq \eta_U \int _{U _{x_{N-1}}} \big  |\frac{\partial u_2}{\partial x _ N} (x_{N-1}, s ) \big | ^ 2  h ( x_{N-1}) \, ds 
\end{equation} 
Integrating over $\om$ and using \eqref{f:33}, we obtain 
$$\begin{array}{ll}
\displaystyle \int _\om |u_2 | ^ 2 h  & \displaystyle \leq \eta _U \int _\om   \Big [ \int _{U _{x_{N-1}}}
\big  |\frac{\partial u _2}{\partial x _ N} (x_{N-1} , s ) \big | ^ 2 h ( x_{N-1}) \, ds   \Big ] \, dx_{N-1} \, dx _ N 
\\ \noalign{\medskip} 
& \displaystyle = \eta_U  \int _U  \chi _{\om} ( x_{N-1}, x_N)    \Big [ \int _{U _{x_{N-1}}} \big  |\frac{\partial u_2}{\partial x _ N} (x_{N-1}, s ) \big | ^ 2 h ( x_{N-1})\, ds   \Big ] \, dx_{N-1} \, dx _ N 
\\ \noalign{\medskip} 
& \displaystyle = \eta_U  \int _a ^ b  \mathcal H ^ 1 (\om_{ x_{N-1} })   \Big [ \int _{U_{x_{N-1}}} \big  |\frac{\partial u_2}{\partial x _ N} (x_1, s ) \big | ^ 2 h ( x_{N-1})\, ds   \Big ] \, dx_{N-1}  
\\ \noalign{\medskip} 
& \displaystyle \leq \eta _U \cdot  \| \mathcal H ^ 1 (\om _{x_{N-1}}) \| _{L^\infty(a, b)}   \cdot   2\,  \lambda _ 2 (\Om) \int _ U| u_2  | ^ 2  h 
\\ \noalign{\medskip} 
&
\displaystyle \leq \eta _U \cdot  \| \mathcal H ^ 1 (\om _{x_{N-1}}) \| _{L^\infty(a, b)}   \cdot   2\,  \frac{\Lambda}{w_\Om ^ 2}  \int _ U| u_2  | ^ 2  h 
\,, 
\end{array} 
 $$ 
where, for the sake of clearness,  we have indicated by $w _\Om$ the  width of $\Omega$. 
Since $\omega$ is a cell of a 
 $(h u _ 1 ^2)$-weighted $L^2$ equipartition of $\overline u$ in $U$, we infer that 
$$ \| \mathcal H ^ 1 (\om _{x_{N-1}}) \| _{L^\infty(a, b)}  \geq  \frac{1}{2 \Lambda} \frac{1}{n} \frac{w _\Om ^ 2}{\eta _U }
\geq   \frac{1}{4 \Lambda}   \frac{  w  _\Om } {n} \,,$$ 
where the last inequality holds because $\eta _U\leq 2 w _\Om$. 
\qed 

\bigskip
\begin{remark}\label{p:propsection2.rem}
If in the above proposition the diameter of $U$ has length at least $\frac{9}{10}D$, 
the angle it forms with the direction $e_{N-1}$ is at most $\arcsin(\frac{10}{9 D} w_{\Om})$. Then the conclusion of the proposition 
continues to hold, possibly with a different constant $\Lambda$, if the local system of cartesian coordinates is changed into
$(e'_{N-1}, e'_N)$, being $e'_{N-1}$ aligned with the diameter  of $U$. 
\end{remark}

 \bigskip
{\it Proof of Theorem \ref{t:quantitative2}}.  It is not restrictive to prove the statement under the hypotheses that $\Om$ is a smooth domain with diameter $\pi$ and small width.  As above, for the sake of clearness, we  denote by $w _\Om$ the  width of $\Omega$,  
and by $\eta _U$
the depth of $U$ in direction $e_N$. 
 We observe that, thanks to the assumption \eqref{f:33}, $\eta _U$ and $w _\Om$ are comparable. Indeed, $\eta _U\leq 2 w _\Om$. To show the converse, namely that also $w _\Om$ is controlled by $\eta _U$, 
 we start from the pointwise inequality \eqref{f:daintegrare}, which holds on $U$, and we integrate it on $U$. 
 Proceeding as in the proof of Proposition  \ref{p:propsection2}, we arrive at 
 $$ \int _U |u_2 | ^ 2 h  \leq \eta _U \cdot  \| \mathcal H ^ 1 (U _{x_{N-1}}) \| _{L^\infty(a, b)}   \cdot   2\,  \frac{\Lambda}{w_\Om ^ 2}  \int _ U| u_2  | ^ 2  h 
\,.
$$ 
 We infer that 
 $$ \int _U  |u_2 | ^ 2 h  \leq 
 \eta _U ^2 \cdot   2\,  \frac{\Lambda}{w_\Om ^ 2}  \int _ U| u_2  | ^ 2  h 
\,,
$$ 
which shows that 
\begin{equation}\label{f:controlloeta}
\eta _U \geq \frac{1}{\sqrt {2 \Lambda } } w _\Om\,.
\end{equation}

 Once we know that the quantities $w _\Om$ and $\eta _U$ are equivalent, 
we 
 indicate by $w_0$ and $\eta _0$ respectively upper bounds for $w _\Om$ and $\eta _U$, and we  
  proceed by adopting the same proof line of Theorem \ref{t:quantitative}.  
The difference is that we have to work with  $(h  u _ 1 ^ 2)$-weighted $L^2$ equipartition of $\overline u$ in $U$, the unique modification with respect to the proof Theorem \ref{t:quantitative}  being the presence of the extra-weight $h$. 
For every $n$, let $\{ \Om _1, \dots, \Om _n \}$ be 
a  $(h  u _ 1 ^ 2)$-weighted  $L^2$ equipartition of $\overline u$ in $U$. 
We set $  D _i   : = {\rm diam} (\Om _i )$ and $h _i$ the profile function of $\Om _i$ in direction orthogonal to a fixed diameter of $\Om _i$. 

We denote by $c$ some number in $(0, 1)$, and 
we follow step by step the same proof line of Theorem \ref{t:quantitative}, distinguishing the same  cases in cascade. 
Below we limit ourselves to indicate which are the required modifications, all the other cases being completely analogue as in Theorem \ref{t:quantitative}:

\begin{itemize}
\item[--] {\it Case 1:}  In place of  Proposition \ref{p:inside}, apply   the consequent inequality 
\eqref{f:diffd2}.   
\item[--] {\it Case 2.2.2:}  Lemma \ref{l:stimaPW} is now applied with $p = h u _ 1 ^ 2$ (notice that such $p$ is still uniformly continuous  on $\Om _i$, since $h$ is continuous on $\overline U$). 
\item[--] {\it Case 2.2.2.2:}  In place of Proposition \ref{p:propsection}, apply Proposition \ref{p:propsection2}  (taking also into account 
Remark \ref{p:propsection2.rem}). \end{itemize} 
\qed

\bigskip

{\it Proof of Theorem \ref{t:quantitative} in dimension $N \geq 3$}.  It is not restrictive to prove the statement under the hypotheses that $\Om$ is smooth and strictly convex, $D = \pi$, and $w$ is small and attained in direction $e _N$. For every $n\in \N$, let $\{ \Om _1, \dots, \Om _n \}$ be a $u _ 1 ^ 2$-weighted $L^2$ equipartition  of $\overline u$ in $\Om$, 
obtained  by the procedure described in Remark \ref{r:PW}, namely  
using  a family of hyperplanes, all parallel to $e _N$, with normals of the type 
$(\cos \alpha_1, \sin \alpha_1, 0, \dots, 0)$, $(0, \cos \alpha_2, \sin \alpha_2, 0, \dots, 0)$, $\dots$, $ (0, \dots, 0, \cos \alpha _{N-2}, \sin \alpha _{N-2}, 0 )$. 
For $n$ large, all the cells become narrow in $(N-2)$-directions, so that they become arbitrarily close to a convex set having at most Hausdorff dimension $2$. 
Since by construction for every cell it holds $\int _{\Om _i} u _ 2 ^ 2 = \frac{1}{n}$, for at most $\frac{n}{2}$ cells it holds
$\int _{\Om _i} |\nabla u _ 2 | ^ 2 \geq 2 \lambda _ 2 (\Om) \int _{\Om_i } u _ 2 ^ 2$. Equivalently,  
\begin{equation}\label{f:16} 
\int _{\Om _i} |\nabla u _ 2 | ^ 2 \leq 2 \lambda _ 2 (\Om) \int _{\Om_i } u _ 2 ^ 2  \qquad \forall i \in  \mathcal I   _n \subset \{1, \dots, n\} \text{ with } {\rm card} ( \mathcal I _n ) \geq \frac{ n}{2}\,.
\end{equation}  

For every $i \in \mathcal I _n$, by
same argument used
to prove the inequality \eqref{f:controlloeta} in the proof of Theorem \ref{t:quantitative2}, we obtain that 
 the depth $\eta _{\Om _i}$ in direction $e _N$ satisfies \begin{equation}\label{f:controlloeta2}
\eta _{\Om _i} \geq \frac{1}{  \sqrt {2 \Lambda } } w _\Om \qquad \forall i \in \mathcal I _n\,.
\end{equation} 

We denote by $\mathcal I _n '$ the subfamily  of $\mathcal I _n$  such that the diameter  $D _i$  of $\Om_i$ satisfies
$ D _ i  \leq \frac{9}{10} \pi $, and we set $\mathcal I _n'' = \mathcal I _n \setminus \mathcal I ' _n$.  
Applying Proposition \ref{p:inside}, we get
$$\mu _ 1 (\Om _ i, u _ 1 ^ 2) \geq 3 \quad \forall i \not \in \mathcal I _n \qquad \text{ and } \qquad  \mu _ 1 (\Om _ i, u _ 1 ^ 2) \geq 3  + C \big ( \frac{\pi}{10} \big ) ^ 3 
\quad \forall i \in \mathcal I' _n \,.$$ 
By Lemma \ref{l:goodparti} it follows that
$$\begin{array}{ll}
\mu _ 1 (\Om, u _ 1 ^ 2 )   & \displaystyle =  \frac{  \int_\Om |\nabla   \overline u  | ^ 2 u _ 1 ^2  }{
 \int_\Om |\overline u  | ^ 2 u _ 1 ^ 2} \geq  \frac{1}{n} \sum _{i = 1} ^ n  \frac{  \int_{\Om_i}  |\nabla   \overline u  | ^ 2 u _ 1 ^2  }{
 \int_{\Om_i}  |\overline u  | ^ 2 u _ 1 ^ 2}  
\\  \noalign{\medskip}
 & \geq  \displaystyle \frac{1}{n}  \Big [ \sum _{i \not \in \mathcal  I _n}  \mu _ 1 (\Om _i, u _ 1 ^ 2 )  +  \sum _{i \in  \mathcal I _n'}  \mu _ 1 (\Om _i, u _ 1 ^ 2 )   +
   \sum _{i \in  \mathcal I _n''}  \frac{  \int_{\Om_i}  |\nabla   \overline u  | ^ 2 u _ 1 ^2  }{
 \int_{\Om_i}  |\overline u  | ^ 2 u _ 1 ^ 2}  
     \Big ] 
\\  \noalign{\medskip}
&   \geq  \displaystyle \frac{1}{n}  \Big [ 3   \big ( n - {\rm card}( \mathcal I _n ) \big  )  +  \big  ( 3 + C \big ( \frac{\pi}{10} \big ) ^ 3 \big )     {\rm card}( \mathcal I' _n ) 
+  \min _{i \in \mathcal I _n ''} \Big [  \frac{  \int_{\Om_i}  |\nabla   \overline u  | ^ 2 u _ 1 ^2  }{
 \int_{\Om_i}  |\overline u  | ^ 2 u _ 1 ^ 2}  
 \Big ]  {\rm card}( \mathcal I _n '')  
  \Big ]\,.
 \end{array}
$$

Therefore, to conclude the proof we are reduced to show that 
$$\liminf _{n \to + \infty} \,  \min _{i \in \mathcal I _n ''}  \Big [  \frac{  \int_{\Om_i}  |\nabla   \overline u  | ^ 2 u _ 1 ^2  }{
 \int_{\Om_i}  |\overline u  | ^ 2 u _ 1 ^ 2}  
 \Big ]    \geq 3 +  c w _\Om ^ 6 \,.$$ 
Indeed in this case we have
$$ \mu _ 1 (\Om, u _ 1 ^ 2 ) \geq 3 +  \frac{1}{2} \Big [ C \big ( \frac{\pi}{10} \big ) ^ 3  \wedge c \Big ] w_\Om ^ 6\,. $$ 
Let  $i _n \in \mathcal I '' _n$ be such that 
 $$\min _{i \in \mathcal I _n ''}  \Big [  \frac{  \int_{\Om_i}  |\nabla   \overline u  | ^ 2 u _ 1 ^2  }{
 \int_{\Om_i}  |\overline u  | ^ 2 u _ 1 ^ 2}  
 \Big ]  =   \frac{  \int_{\Om_{i_n}}  |\nabla   \overline u  | ^ 2 u _ 1 ^2  }{
 \int_{\Om_{i_n}}  |\overline u  | ^ 2 u _ 1 ^ 2}  
  \,.$$
 In the sequel we write for brevity $\Om^n$ in place of $\Om _{ i _n}$. So our target is to show that
 \begin{equation*} \liminf _{n \to + \infty} \frac{  \int_{\Om^{n}}  |\nabla   \overline u  | ^ 2 u _ 1 ^2  }{
 \int_{\Om^{n}}  |\overline u  | ^ 2 u _ 1 ^ 2}   \geq 3 +  c w _\Om ^ 6 \,.
 \end{equation*}

 Let $H^n$ be the intersection of closed halfspaces parallel to $e _N$ such that $  \overline {\Omega ^n} = \overline \Om \cap H ^n$. 
Up to subsequences, and up to changing the coordinate system, we can assume that $H^n$ converge in Hausdorff distance to the hyperplane 
$$\Pi := \big \{ x \ :\ x_i  = 0 \ \forall i = 1, \dots, N-2 \big  \} \,.$$  
In the sequel, a point $(0, \dots, 0, x _{N-1}, x _N) \in \Pi$ will be identified with the pair $(x_{N-1}, x_N)$. 

Accordingly, the sequence $\Om ^{n}$ converge in Hausdorff distance to the set $U := \Om \cap \Pi$, which is of the kind 
$$U = \Big \{ (x_{N-1}, x _N)  \ :\ x_{N-1} \in (a, b) \,, \ x _N \in \Om \cap ( x _{N-1} + \R e _N ) \Big \} \,.$$
We remark that $U$  is nondegenerate, namely it has positive two-dimensional measure.   Indeed, the depth of $U$ in direction $e_N$ is strictly positive because, 
by \eqref{f:controlloeta2},
the depth of $\Om ^{n}$ in direction $e_N$ is uniformly bounded from below. 
Moreover, the length of $(a, b)$ is strictly positive, because
the diameter of $\Om^n$ is not smaller than $\frac{9\pi}{10}$, 
and cannot be attained in direction $e_N$  (which is the direction of the width, assumed to be small).  

Moreover $U$ cannot lie on $\partial \Om$, thanks to our initial assumption of strict convexity on $\Om$. 

For every $(x _{N-1}, x _N) \in  \Pi$, we define the function 
$$h _n ( x_{N-1}) = \mathcal H ^ { N-2} \big (   \big (  \Pi ^\perp_{(x_{N-1}, x_N)} \big  ) \cap H { ^n}  \big ) \,,$$
 where $ \Pi ^\perp_{(x_{N-1}, x_N) }$ denotes the $(N-2)$-affine space passing through  $(x_{N-1}, x_N)$ and orthogonal to $\Pi$. 
 
Up to subsequences, $\frac {h _n }{\| h _n \|_\infty}$   converges a.e.\ on $\R e _{N-1}$ to a function $h$
such that $h = 0 $ on $(- \infty , a) \cup (b, + \infty)$. Moreover, since 
  by the Brunn-Minowski Theorem $h _n$ is $\big ( \frac{1}{N-2}\big)$-concave on $(a, b)$, the convergence is 
locally uniform on $(a, b)$, $h$ is itself $\big ( \frac{1}{N-2}\big)$-concave on $(a, b)$, and consequenty satisfies also $\| h \| _\infty = 1$.  
  
  We claim that
\begin{equation}\label{f:claimU1} \lim _{n \to + \infty}\frac {1}{\| h _n \|_\infty}   {  \int_{\Om^{n}}  f } = \int _U f h 
\qquad \text{ for every } f \in \mathcal C ^ \infty (\overline \Om) \,.
\end{equation} 
and
\begin{equation}\label{f:claimU2} 
u _2 \text{ is not identically zero on } U\,.
\end{equation} 
Assume by a moment these two claims hold true.  
Then we infer that $\Om$, $U$, and $h$ satisfy all the assumptions of Theorem \ref{t:quantitative2}. 
Indeed,  recall
that $\Om ^n$ belongs to the family of cells satisfying \eqref{f:16}, and pass to the limit as $n \to + \infty$: 
by using \eqref{f:claimU1} with $f = |\nabla u _ 2 | ^ 2$ and with $f = u _ 2 ^2$,  
it follows that  assumption \eqref{f:33} is fulfilled. Similarly,
recalling that $\int _{\Om ^n}   \overline u u _ 1 ^ 2 = \int _{\Om ^n}   u _ 1 u _ 2 = 0$, and using  \eqref{f:claimU1} with $f = u _ 1 u _ 2$, it follows that also 
the assumption that $\int _U \overline u h u _ 1 ^ 2= 0$ is satisfied. Finally, 
 the assumption that $u _ 2$ does not vanish identically on $U$ is satisfied by  \eqref{f:claimU2}.   Then, we have
  \begin{equation*}\lim _{n \to + \infty} \frac{  \int_{\Om^{n}}  |\nabla   \overline u  | ^ 2 u _ 1 ^2  }{
 \int_{\Om^{n}}  |\overline u  | ^ 2 u _ 1 ^ 2}  =  \frac{  \int_{U}  |\nabla   \overline u  | ^ 2 h u _ 1 ^2  }{
 \int_{U}  |\overline u  | ^ 2  h u _ 1 ^ 2}   \geq  3  + \overline c {w^6} \,,
 \end{equation*} 
 where the first equality is obtained applying again \eqref{f:claimU1} with $f = |\overline u  | ^ 2 u _ 1 ^ 2 = u _ 2 ^ 2$ and with $f= |\nabla \overline u  | ^ 2  h u _ 1 ^ 2$, 
 and the second inequality follows from Theorem \ref{t:quantitative2} (since $ |\nabla   \overline u  | ^ 2 \geq |\nabla_U   \overline u  |  ^2$).  

To conclude our proof, we now give the proofs of claims \eqref{f:claimU1} and \eqref{f:claimU2}. 

\smallskip
$\bullet$ Proof of claim \eqref{f:claimU1}: Let $\delta _n$  denote the Hausdorff distance between $\Om ^n$ and $U$, and set
$$U ^{-\delta_n} =  \big \{ x \in U \ :\ {\rm dist} (x, \partial U )  > \delta _n \big \}  \,, \qquad 
U ^{\delta_n} =  \big \{ x \in \Pi \ :\ {\rm dist} (x, \partial U )  <\delta _n \big \} \,.$$ 
We have
$$\int_{\Om^{n}}  f  = 
\int _\Pi \int _{\Pi ^\perp_{(x_{N-1}, x_N)}} \!\!\!\!\!\!\!\!\chi _{\Om ^n} f =
\int _{ U ^{-\delta_n}} \int _{\Pi ^\perp_{(x_{N-1}, x_N)}} \!\!\!\!\!\!\!\!\chi _{\Om ^n} f  + \int _{ U ^{\delta_n} \setminus  U ^{-\delta_n}}  \int _{\Pi ^\perp_{(x_{N-1}, x_N)}} \!\!\!\!\!\!\!\!\chi _{\Om ^n} f   \,,
$$ 
where the integrals over $\Pi ^\perp_{(x_{N-1}, x_N)} $ are made with respect to  $x' = (x_1, \dots, x_{N-2})$, while the integrals over 
 $U ^{-\delta_n}$,   $U ^{\delta_n} \setminus  U ^{-\delta_n}$ are made with respect to $(x_{N-1}, x_N)$. 

For $(x_{N-1}, x_N) \in U ^{-\delta_n} $,  setting $h = (x', 0, 0 )$, we have 
 $$f ( x', x_{N-1}, x_N) = f (0, x_{N-1}, x_N )  + \nabla f (0, x_{N-1}, x_N )   \cdot h + o ( |h |)\, ;$$ 
 for $(x_{N-1}, x_N) \in (U ^{\delta_n}  \setminus  U ^{-\delta_n}) $, if $(z_{N-1} , x_N)$ is the projection   of
  $(x_{N-1}, x_N)$ onto $\overline U$ parallel to  $e_{N-1}$, setting   
 $\tilde h = (x', x_{N-1}- z _{N-1}, 0 )$ we have 
 $$f ( x', x_{N-1}, x_N) = f (0, z_{N-1}, x_N )  + \nabla f (0, z_{N-1}, x_N )   \cdot  \tilde h + o ( | \tilde h |)\,.$$ 
 Accordingly,  we have  
 $$\begin{array}{ll} &\displaystyle \int _{ U ^{-\delta_n}} \int _{\Pi ^\perp_{(x_{N-1}, x_N)}} \!\!\!\!\!\!\!\!\chi _{\Om ^n} f  = 
 \int _{ U ^{-\delta_n}}   h _n (x_{N-1}) ( f + O (\delta _n))\,,  \\ 
 \noalign{\smallskip}
&   \displaystyle \int _{ U ^{\delta_n} \setminus  U ^{-\delta_n}} \!\! \int _{\Pi ^\perp_{(x_{N-1}, x_N)}}  \!\!\!\!\!\!\!\!\chi _{\Om ^n} f  = 
 \int _{ U ^{\delta_n} \setminus  U ^{-\delta_n} }   \psi _n (x_{N-1}, x _N) ( f + O (\delta _n))\,,
\end{array}
 $$ 
where $\psi_n$ is  defined by 
 $$\psi _n ( x_{N-1}, x_N) := \mathcal H ^ { N-2} \big (   \big (  \Pi ^\perp_{(x_{N-1}, x_N)} \big  ) \cap \Om { ^n}  \big ) \,.$$
We now divide by $\|h _n\|_\infty$ and pass to the limit as $n \to + \infty$.  
 Taking into account that 
$\psi _n \leq h _n$, and  that $\mathcal H ^ {N-2} \big (U ^{\delta_n} \setminus  U ^{-\delta_n} \big )$ is infinitesimal, we conclude that
$$ \lim _{n\to + \infty}   \frac {1}{\| h _n \|_\infty}   {  \int_{\Om^{n}}  f } =  \lim _{n\to + \infty}  \frac {1}{\| h _n \|_\infty}   
\int _{ U ^{-\delta_n}}     f  {h _n}   = \int _U f h \,. $$ 

\smallskip
$\bullet$ Proof of claim \eqref{f:claimU2}:  Assume by contradiction that $u _2$ is identically zero on $U$. 
We are going to show that this implies
\begin{equation*}\frac{  \int_{\Om ^n} | \nabla u _ 2  | ^ 2  }{ \int _{\Om ^n} |u _ 2  | ^ 2}  \to + \infty \qquad \text{ as } n \to + \infty\,,
\end{equation*}
against \eqref{f:16}. To that aim we are going to use an argument which amounts roughly to control the value of the function by the value of its gradient, in the same spirit of  the {\L}ojasiewicz inequality \cite{Lo59}.
   
  For every $i \in \N$ let us denote by $D_{x_U}  ^ {(i)} u _ 2$  the $i$-th order differential of $u _2$ computed at a point
 $x_U=(0, \dots, 0,x_{N-1}, x_N)$ of $U$. 
   Let $k \in\N \setminus \{0\} $ be the smallest natural number such that, for some point $x_U$ of $U$, 
  $D_{x_U}  ^ {(k)} u _ 2 \neq 0$. 
  Clearly such $k$ exists:
 otherwise, by the analyticity of $u _2$ inside $\Om$, $u _2$ would be identically zero.

  In order to estimate the Rayleigh quotient of $u _ 2$ over $\Om ^n$, we distinguish  between points $x = (x' , x_{N-1}, x_N ) \in \Om ^n$ such that $(x_{N-1}, x_N) \in U$ and such that $(x_{N-1}, x_N) \not \in U$. 

For $x =(x', x_{N-1}, x_N) \in \Om ^ n$ such that $(x_{N-1}, x_N) \in U$, setting  $ x_U= (0,\dots, 0, x_{N-1}, x_N)$ and $ \xi   = (x', 0, 0)$,  we have
$$\begin{array}{ll}
& \displaystyle u _ 2 ( x) = \sum _{i= 0}^k \frac{1}{i !}  D^{( i)} _ {x_U} u _ 2   [ \xi  ^ {(i)}]  + o (| \xi | ^ {k}) 
\\ \noalign{\medskip}
& \displaystyle \nabla  u _ 2   ( x) = \sum _{i= 0}^{ k-1} \frac{1}{i !}  D^{( i)} _ {x_U} \nabla u _ 2   [ \xi  ^ {(i)}]  + o (| \xi | ^ {k-1})  \,.
\end{array}
 $$ 
Since
$$D  ^{(k)} _{{x_U}}   u _ 2  [  \xi   ^{(k)}] =D  ^{(k-1)} _{{x_U}}  (\nabla   u _ 2  \cdot  \xi ) [  \xi   ^{(k-1)}  ]
=  \xi  \cdot \big ( D  ^{(k-1)} _{x_U}  (\nabla   u _ 2  ) [ h ^ {(k-1)}  ]  \big )\, , $$ 
by applying Cauchy-Schwarz inequality we obtain
$$| D  ^{(k)} _{x_U}   u _ 2  [  \xi   ^{(k)} ]|  \leq |  \xi  |    |   D  ^ {(k-1)}  _{x_U}  (\nabla    u _ 2  ) [  \xi ^ {(k-1)}  ]   | \,.$$ 

It follows that
\begin{equation} \label{abf01}
|u _ 2 ( x) | ^ 2 \leq | \xi | ^ 2   |  \nabla  u _ 2   ( x)  | ^ 2 +  o (| \xi | ^ {2k} )   \quad \forall x \in \Om ^ n \ :\  (x_{N-1}, x_N) \in U\,.\end{equation}
 For $x \in \Om ^n$ such that $(x_{N-1}, x_N) \not \in U$, we let $(z_{N-1}, x_N)$ be the point introduced in the above proof of claim \eqref{f:claimU1}, 
we set $  \tilde  \xi  = (x', x_{N-1}-z_{N-1}, 0)$, and we argue in a similar way as above. We obtain  
\begin{equation} \label{abf02}
|u _ 2 ( x) | ^ 2 \leq | \tilde  \xi | ^ 2   |  \nabla  u _ 2   ( x)  | ^ 2 + o (|   \tilde  \xi | ^ {2k} ) \quad \forall x \in \Om ^ n \ :\  (x_{N-1}, x_N) \not \in U \,. 
\end{equation}
Integrating $|u_2|^2$ over $\Om^n$ and taking into account that, by the $(\frac{1}{N-2})$-concavity of $h_n$,  for small $\delta_n$ we have   $ \int_U h_n \ge \frac{|\Om^n|}{2}$, we get, for some positive constants $C_1$ and $C_2$, 
 $$C_1 |h|^{2k} \int_U h_n \ge \int _{\Om ^n} |u _ 2  | ^ 2 \ge  \frac{1}{(k !)^2} \int_U h_n\big( |D^ {(k)}  _ {x_U} u _ 2 [h^{(k)}]|^2 + o (|h| ^ {2k}) \big) \ge C_2 |h|^{2k} \int_U h_n\,.$$

Now, summing \eqref{abf01}-\eqref{abf02} over $\Om_n$, we obtain 
$$
\frac{\delta_n^2 \int _{\Om ^n} |\nabla u _ 2  | ^ 2 }{\int _{\Om ^n} |u _ 2  | ^ 2 }\ge 1+  {o(1)}\,.$$
This is not possible since $\delta_n \to 0$ and, from \eqref{f:16}, the ratio $\frac{ \int _{\Om ^n} |\nabla u _ 2  | ^ 2 }{\int _{\Om ^n} |u _ 2  | ^ 2 }$ is bounded from above. As a conclusion, our assumption that  $u _2$ is identically zero on $U$ fails to be true,
yielding \eqref{f:claimU2}.  \qed 

\section{The Neumann gap}\label{sec:Neumann}

  Through  the addition of  few specific new results, the approach developed in the previous sections
leads to a quantitative lower bound for the  first nontrivial Neumann eigenvalue 
$\mu _ 1 (\Om, \phi)$ defined according to   \eqref{def:mu}, 
where $\phi$ is a generic positive weight in $L ^ 1 (\Om)$    which is no longer related to the first eigenfunction, but is power-concave, a particular case being $\phi \equiv 1$. 

Actually,  the statement of Theorem \ref{t:quantitativeN} holds more generally with $\mu _1 (\Om)$ replaced by $\mu _ 1 (\Om, \phi)$, being $\phi$
any weight which is $(\frac{1}{m})$-concave for some $m \in \N \setminus \{ 0 \}$. Such assumption  is needed not only for the existence of an eigenfunction $\overline u$ (see Proposition \ref{lemma371} in Appendix), but also for the control of the constant $\overline c$ in Theorem \ref{t:quantitativeN}.

The proof proceeds along the same line as Theorem \ref{t:quantitative}, 
being considerably simpler and yet demanding some nontrivial new ingredients.  
The main difficulty arises from the fact that an eigenfunction $\overline u$ can no longer be identified with 
$\frac{u _ 2}{u _ 1}$.  We point out that  this identification was crucial  to obtain  Proposition \ref{p:propsection}, 
which in turn allowed to reach a contradiction in Case 2.2.2.2 of our proof of Theorem \ref{t:quantitative} for $N = 2$. 

To overcome such difficulty, we  manage to acquire a control on the Lebesgue measure of the cells of the partition 
in terms of the width of $\Om$.  This will be possible thanks to
 the following new  geometrically explicit  $L ^ \infty$ estimate for Neumann eigenfunctions   (see   \cite{CM13} for global Lipschitz regularity results):

\begin{proposition}\label{ambu25}
There exists a  positive constant  $C$ depending only on $N+m$ such that, 
for every open bounded convex domain $\Omega\subset \R ^N$ with diameter $D_\Om$  
and any positive  $\frac1m$-concave weight $\phi$,     a first eigenfunction 
$ \overline u$  associated with the Neumann eigenvalue $\mu_1 (\Omega, \phi)$ satisfies
\begin{equation}\label{f:inftybound}
	\norma{ \overline u}_{L ^ \infty(\Om)} \leq C \mu _1 (\Omega, \phi)^{\frac{N+m}{2}} \frac{D_\Om^{N+m}}{{\left (\int_\Omega 
\phi \right )}^{1/2}} \norma{\overline u}_{L^2(\Omega,\phi)}.
	\end{equation} 
	\end{proposition}

\begin{remark}

	(i) In the particular case $\phi \equiv 1$,   the estimate \eqref{f:inftybound}  reads
	\begin{equation}
	\label{ambu40}
		\norma{\overline u}_{L^\infty(\Om)} \leq C \mu _1 (\Om)^{\frac{N}{2}}  \frac{D _\Om ^N}{\abs{\Omega}^{\frac{1}{2}}}
		\norma{\overline u}_{L^2(\Om)}  \,. 
	\end{equation}
		\smallskip 
		
  (ii) As it can be seen for the proof below in case $\phi =1$, the  result continues to hold for any Neumann eigenpair. 	  

	\smallskip 
	
	(iii)      Combined with 
	 the result proved by Maz'ya in \cite[Section 2]{Mazia}, 
	 and with the upper bound for  the relative isoperimetric constant  stated in  \cite[Theorem 1.4]{relative_ok},
	 the inequality  \eqref{ambu40},   written for an arbitrary Neumann eigenpair $(\mu _ k (\Om), u _k)$   gives the following gradient estimate  
\begin{equation*}
	\norma{\nabla u_k}_{L^\infty (\Om)} \leq  C \mu_k (\Om)^{\frac{N}{2}+1}  
	\frac{D_\Om^{N+1}}{\abs{\Omega}^{\frac{1}{2}}}   \norma{u_k}_{L^2 (\Om)} .
\end{equation*}	
\end{remark}

\begin{proof}  Let us start by proving the result  for $\phi \equiv 1$.      Let $a_1 \geq a_2 \geq \dots \geq a_N $ denote the semi-axes of the John ellipsoid of $\Om$. 
		Since \eqref{ambu40} is invariant under scaling, we are going to assume without loss of generality that $a_1= 1$. 
In order to reduce ourselves to work with domains having the unit ball as John ellipsoid, we perform the change of variables $X=T(x)$, with 
	$$
	X_1=x_1, \quad X_2= \frac{x_2}{a_2}  , \quad \dots\quad,  X_N= \frac{x_N}{a_N} \,. 
	$$
Taking into account that $a_1$ is comparable to $D_\Om$, it is readily checked that, in terms
	of the function $v$ defined on $T (\Om)$ by 
	$$
	v(X_1,\dots , X_N) := \overline u(X_1, a_2  X_2, \dots , a_N X_N ),
	$$
	the inequality we want to prove becomes:
	\begin{equation}
		\label{intermedio.1}
		\norma{v}_{L^\infty(T(\Om))} \leq C \mu_1 ^{\frac N 2}(\Omega) \norma{v}_{L^2(T(\Om))} \,.
	\end{equation}

Setting for brevity $A: = T (\Om)$ and $c = ( c_1, \dots, c _N)  :=  ( \frac{1}{a_1^2}, \dots, \frac{1}{a_N^2} ) $, it holds 
	$$ \mu _ 1 (\Om):= 
	\frac{{\int_{\Omega} \abs{\nabla \overline u}^2 }}{{\int_{\Omega}  \overline u^2} }= 
	\frac{\sum_{i=1}^N\int_{A} \big| \frac{\partial v}{\partial X_i} \big | ^2 c_i  \, dX}{{\int_{A}v^2  \, dX}  }= : \mu _ c ( A)  \,,
	$$
	 and the optimality of $\overline u$ can be reformulated  as the following variational property of $v$ 
	\begin{equation}\label{f:weakv} \sum_{i= 1} ^N  \int_{A}  c_i \frac{\partial v}{\partial x_i}\frac{\partial \vphi}{\partial x_i}= \mu_c (A) \int_A v \vphi \qquad  \forall \vphi \in H^1(A) \,.
	\end{equation} 
	We now use the Moser iteration scheme: we claim that, 
	if a solution $v$ to \eqref{f:weakv} belongs to $L^p(A)$, then 	it belongs also to $L^{\frac{p N}{N-1}}(A)$;  more precisely, 
for a dimensional constant $C$,  it holds 
			\begin{equation}\label{f:moser} 
	\Big( \int_{A} |v|^{\frac{
  pN  }{N-1}}\Big )^{\frac{N-1}{N}} \leq 
	C  \mu_c (A) 
	\frac{ p^2}{2 (p-1)}	\int_{A} |v|^p\,. 
	\end{equation} 
  
To prove the claim we observe that, if $v \in L ^ p (A)$, we can  take 	
	$|v|^{p-2}v$ as test function in \eqref{f:weakv}. This gives 
	$$
	\sum_{i=1}^N \int_{A} c_i \frac{\partial v}{\partial x_i}\frac{\partial |v|^{p-2}v}{\partial x_i} = \mu_c (A) \int_A |v|^p.
	$$
	Since $c_i \geq 1$, we infer that
	$$
	\frac{p-1}{\left( \frac{p}{2}\right)^2}\sum_{i=1}^N  \int_{A}  \Big (\frac{\partial |v|^{\frac{p}{2}}}{\partial x_i}\Big )^2 \leq \mu_c  (A) \int_A |v|^p,
	$$
	and hence
	$$
	\int_A \Big | \nabla |v|^{\frac{p}{2}}  \Big | ^2=	\sum_{i=1}^N \int_{A } \Big (\frac{\partial |v|^{\frac{p}{2}}}{\partial x_i}\Big )^2 \leq 	\mu_c (A) \frac{p^2}{4(p-1)}   \int_A |v|^p.
	$$
	
Since
$$
		\int_{A} |v|^p +  \int_A \abs{\nabla |v|^{p}}= 	\int_{A} |v|^p +\int_U \abs{\nabla (|v|^{\frac{p}{2}})^2}\\
		= 	\int_{A} |v|^p  +2 \int _A |v|^{\frac{p}{2}} \Big | \nabla |v|^{\frac{p}{2}}\Big | \,,
		$$
		for every positive constant $\alpha>0$ we have 
		\begin{equation}\label{f:alfa} 
		\begin{array}{ll} \displaystyle \int_{A} |v|^p +  \int_A \abs{\nabla |v|^{p}} & \displaystyle \leq 
		 	\int_{A} |v|^p+ 	\alpha  \int_{A} |v|^p +  \frac 1\alpha \int_A \abs{\nabla |v|^{\frac{p}{2}}}^2
			\\  \noalign{\medskip}
		 &  \displaystyle \leq   \Big (\alpha  +1+  \frac 1\alpha \mu_c (A) \frac{ p^2}{4(p-1)}  \Big )	\int_{A} |v|^p\,.  \end{array} 
\end{equation}

	Next we observe that  there exist  positive dimensional constants $C'$ and $C''$ such that, 
	for every function $w \in W ^ { 1, 1} ( A)$ (extended to zero outside $\overline A$),  it holds 
	
	\begin{equation}
		\label{sobol} \Big ( \int_A  \abs{w}^{\frac{N}{N-1}} \Big ) ^{\frac{N-1}{N}} \leq C' \Big (  \int_A \abs{\nabla w} + \int_{\partial A} \abs{w} \Big )  \leq C'' \Big  ( \int_A \abs{\nabla w} +\int_A \abs{ w} \Big )\,.\end{equation}
		Here the first inequality is due to the continuity of the embedding of $BV ( \R ^N)$ into $L ^ {\frac{N}{N-1}} (\R^N)$,   and the second one to the continuity of the trace operator from $W ^{1, 1} (A)$ to $L ^ 1 (\partial A)$. 
	Notice in particular that the fact that $C''$ is purely dimensional is due to the condition that the John ellipsoid of $A$ is the unit ball  (so that the outradius and the inradius of $A$ are controlled respectively from above and from below). 

		By applying \eqref{sobol}  with $w = |v|^p$,  we infer  from \eqref{f:alfa}  that 
		$$
	\Big ( \int_{A} |v|^{\frac{p N}{N-1}}\Big )^{\frac{N-1}{N}} \leq  C '' \Big (\alpha  +1+ \frac 1\alpha \mu_c(A) \frac{  p^2}{4(p-1)}  \Big )	\int_{A} |v|^p.
	$$
The validity of our claim \eqref{f:moser} (with  $C = \frac{C''}{\alpha}$) 
follows from the above inequality after noticing that, by Payne-Weinberger inequality, and since we have fixed $a_ 1= 1$, $\mu_c (A)=  \mu _ 1 (\Om)$ is bounded from below by 
a dimensional constant so that,  for  $\alpha$ below a dimensional threshold, 
$$\alpha +1 \leq \frac 1\alpha \mu_c(A)  \frac{  p^2}{4(p-1)}  \,. $$

	We now apply \eqref{f:moser} recursively: setting $p_0=2$ and $p_{k+1}= \frac{N}{N-1} p_{k}$ (i.e., $p _k=  2 \big ( \frac{N}{N-1}  \big ) ^ k$), this gives
	$$
	\norma{v}_{L^{p_{k+1}} (A)} \leq  (C \mu _ c (A)) ^{\frac{1}{p_k} }    
	\left(\frac{ p_k^2}{2 (p_k-1)} \right)^{\frac{1}{p_k} }	\norma{v}_{L^{p_{k}}(A)}\,.
	$$
	In the limit as $k \to + \infty$ we obtain 
		$$
	\norma{v}_{L^\infty(A)} \leq \prod_{k=0}^\infty (C  \mu_c(A))^{\frac{1}{p_k} }  \Big (\frac{ p_k^2}{2(p_k-1)} \Big )^{\frac{1}{p_k} }\norma{v}_{L^2(A)}.
	$$
	
The validity of the required inequality \eqref{intermedio.1} (for some dimensional constant $C$) follows by noticing that  
		$$
	\sum_{k \geq 0} \frac{1}{p_k} = \frac{1}{2} \sum_{k \geq 0} \left(\frac{N-1}{N}\right)^k= \frac{N}{2},
	$$
	and
	$$
	\prod_{k=0}^\infty  \Big (\frac{p_k^2}{2(p_k-1)}\Big)^{\frac{1}{p_k} }  \leq\prod_{k=0}^\infty  p_k^{ \frac{1}{p_k} }<+\infty.$$
	
	 The case of a general $(\frac{1}{m})$-concave weight $\phi$ is obtained by collapsing, relying on the results proved in the Appendix. Precisely,
we write  inequality \eqref{ambu40}  on the domain   $\widetilde{\Omega }_\vps \subset \R^{N+m}$ defined in  \eqref{lemma371.4}, 
and we obtain the validity of \eqref{f:inftybound} in the limit as $\vps\ra 0$, by using the lower semicontinuity inequality \eqref{lemma371.3} from Proposition \ref{lemma371}.

\end{proof}

 \bigskip
 
We are now in a position to give the 
proof of Theorem \ref{t:quantitativeN}.   The idea is the same as for the Dirichlet case: reduce the $N$-dimensional gap estimate to a 
two-dimensional one. However, the two-dimensional one will involve a geometric weight, and so we are going to prove  it  more in general replacing $\mu _ 1 (\Om)$ by  $\mu _ 1 (\Om, \phi)$, where $\phi$ is a fixed positive weight which is $(\frac{1}{m})$-concave for some $m \in \N \setminus \{ 0 \}$.  

To avoid cumbersome overlaps, we are going to outline the parts which closely follow the proof of Theorem \ref{t:quantitative}, focusing our attention on the 
differences. 

As in the case of the Dirichlet fundamental gap, by the  Payne-Weinberger reduction argument, estimating from below $\mu _ 1 (\Om, \phi)$ amounts to estimating from below  a weigthed one-dimensional Neumann eigenvalue of the type
$\mu _ 1 (I, p)$. 
The  difference  must be searched in the weight: now  $p= h \phi$, where
$\phi$ is the preassigned power-concave weight (which replaces $u _ 1 ^ 2$), while  $h$ is still, as in the Dirichlet case,  the  power-concave function giving 
the $(N-1)$-dimensional measure of the cell's section orthogonal to $I$.  As a consequence, the one dimensional part of the proof is much simpler. 
Indeed, by the log-concavity of $p$, 
Lemma \ref{l:ml} still holds, so that
$$\mu _ 1 (I, p ) \geq \lambda _ 1 (I, m_p)\,, \qquad \text{ with }  m _p
: =  \Big [ \frac{3}{4} \Big (\frac{p'}{p} \Big ) ^ 2 - \frac{1}{2} \frac{p''}{p} \Big ] \,. $$
Now, the counterpart of the sharp one dimensional lower bound  for $\lambda _ 1 ( I _, m _p)$ given in Theorem \ref{t:1d}
reads simply as follows: since $p$ is log-concave we have that 
 $m _p$ is a positive measure and hence, 
working for definiteness of the interval $I _\pi$,   we have
$$\lambda _ 1 (I_\pi, m_p) \geq \lambda _ 1 ( I_\pi)= 1\,.$$

By analogy with Theorem \ref{t:1dextra}, 
keeping the notation 
$I _ d= ( - \frac d 2, \frac d 2)$ for any $0<d < \pi$, the above inequality can be refined in two distinct directions: 

\smallskip
(i)  There exists an  absolute constant $C>0$  such that 
\begin{equation} \label{f:stima3lengthbis}  
\lambda_ 1  (I _{d}, m _p) \geq  \lambda_1(I_d,0) = \frac{\pi^2}{d^2} \geq  1  + C(\pi - d)   \,. 
\end{equation} 

(ii) There exists an  absolute constant $K$  such that, 
if $p=h \phi$, with  $h$ log-concave and 
affine on an interval $[a,b]$  with 
$I_{ \frac{\pi}{4} } \subseteq [a, b] \subseteq I _d$,
 the following implication holds:  
\begin{equation}\label{f:stima3deltabis} 
\lambda_1 ( I _ d , m _p) \leq 2 \ \Rightarrow \ \lambda _1 ( I _ d, m_p ) \geq 1 +
  \Big [ 1 - \frac{ \min\{ h ( a), h (b)  \}   }{ \max \{ h ( a), h (b)  \} } \Big ] ^2    
 \,.
 \end{equation}
 
We remark that the lower bound \eqref{f:stima3lengthbis}, 
in which the exponent $1$ replaces the exponent $3$ appearing in 
\eqref{f:stima3length}, is a straightforward consequence of  the Taylor expansion of $\frac{\pi ^ 2}{(\pi - \vps  )  ^ { 2}} $ as $ \vps   \to 0^+$. 
On the other hand, 
the lower bound  \eqref{f:stima3deltabis} can be proved in similar way as \eqref{f:stima3delta}. 

Passing to higher dimensions,  
the results in Section \ref{sec:partitions}  about weighted measure or $L ^2$ equipartitions
remain unaltered working with our new weight.

Then, by using measure equipartitions, 
as a counterpart to Proposition \ref{p:inside},  
we obtain the following  localized version of  the (weighted) 
Payne-Weinberger inequality: 
   \begin{equation}\label{f:diffdbis} \mu _ 1 (\om, \phi ) \geq \frac{ \pi ^ 2}{D^2}  + C \frac{(D-d) }{D ^ 3} \,,
\end{equation}
where $\om \subseteq \Om\subset \R ^ N$ ($N \geq 2$) are  open bounded convex sets of diameters $d<D$.  

With this inequality at our disposal, we proceed to  prove the quantitative inequality. 

  Notice in particular that it is not necessary to prove preliminarily  rigidity as  in the Dirichlet case: if the quantitative form of the Payne-Weinberger inequality holds true 
for cells with diameter larger than $\frac\pi2$ and second John semi-axis smaller than a fixed threshold $\overline a_2$,
rigidity follows as direct  consequence. 
Indeed, we consider a partition of  $\Om$ into cells with second John semi-axis smaller than $\overline a_2$. 
If a cell has diameter smaller than $\pi$, we get  rigidity.  Otherwise, all cells have diameter equal to $\pi$, 
and we get as well rigidity by writing   the quantitative form of the inequality for a single cell, which has a strictly positive second John semi-axis. (Notice that, in the Dirichlet case, due to the presence of the weight $u _ 1 ^ 2$, 
a quantitative inequality for convex sets with small width does not imply a quantitative inequality on  a single cell).

The proof of Theorem \ref{t:quantitativeN} is carried over first in case $N = 2$ and then for $N \geq 3$. 
In dimension $N = 2$, we are going to use, in place of  Proposition \ref{p:propsection}, the following result obtained via the $L ^ \infty $ estimate given in Proposition \ref{ambu25}.

\begin{proposition}\label{p:proparea} 
There exists a positive constant $\Lambda$, depending only on $N+m$, such that: if $\Om \subset \R^N$ is an open bounded convex set, 
$\phi$ is a positive $\big ( \frac 1m\big )$-concave weight in $L ^ 1 (\Om)$ and  $\overline u$ is a first eigenfunction for $\mu _ 1 (\Om, \phi)$ normalized in $L ^ 2 (\Om, \phi)$,  
and $\om$ is a cell of a $\phi$-weighted $L^2$ equipartition of $\overline u$ in $\Om$ composed by $n$ cells, it holds 
	\begin{equation*}
		|\omega| \geq \Lambda \frac{|\Omega|}{n}.\;\;\;
	\end{equation*}
\end{proposition}

\proof  By Proposition \ref{ambu25}, we have
	\begin{equation*}
		 \frac{1}{n\int _{\omega}\phi }= \frac{\int_{\omega} \phi \overline u^2 }{\int _{\omega} \phi } \leq \norma{\overline u}^2_{L^\infty (\Om)} \leq  C \mu_1(\Omega,\phi )^{N+m} \frac{D_\Om ^{2(N+m)}}{  \int_\Om \phi}.
	\end{equation*}

Since $\mu_1(\Om,\phi ) D_\Om ^{2(N+m)}$ is bounded  from above by a constant $K$ depending only on $N +m$ (see \cite{kro,Henrot-Michetti} for $\phi \equiv 1$ and Proposition \ref{lemma371} for the general case),  we conclude that
$$
\int _{\omega}\phi   \geq \frac{K}{n}  \int _\Om \phi  \,.$$ 

On the other hand, we have   $\|\phi\|_{L^\infty(\Om)}<+\infty$ and, assuming with no loss of generality that  $0 \in \overline \Om$ is  a maximum point for $\phi|_\Om$, the  $\big ( \frac 1m)$-concavity of $\phi$ implies that
$\phi (x) \ge  \frac {1}{2^m} \|\phi\|_{L^\infty(\Om)}$ for every  $x \in \frac 12 \Om $. Hence 
$$\int _\Om \phi  \geq \frac{1}{2 ^ {N+m} } \|\phi \|_{L^\infty(\Om)} \abs{\Om} \,.$$

We deduce that 
	$$\abs{\omega} \|\phi \|_{L^\infty(\Om)} \ge \int _{\omega} \phi  \ge  
	\frac{K}{n}  \int _{\Omega} \phi  \ge \frac{ K}{n 2^{N+m}} \|\phi \|_{L^\infty(\Om)} \abs{\Om}\,,$$
	and the result follows with $\Lambda =  \frac{ K}{2^{N+m}}$. 
	\qed

\bigskip

\bigskip 

{\it Proof of Theorem \ref{t:quantitativeN} in dimension $N = 2$}. 
We follow the same geometric construction as in the proof of Theorem \ref{t:quantitative}. 
Relying on the inequalities \eqref{f:stima3deltabis} and  \eqref{f:diffdbis}, all the cases work as previously, 
exception made for the last one, Case 2.2.2.2.
To prove that this case cannot occur, the contradiction argument is still based on the validity of claim \eqref{f:sezlunga}. The proof of such claim proceeds unaltered up to the inequality \eqref{f:lminlmax}. 
At this point we have to argue differently. 
Indeed, in the Dirichlet setting, 
the estimate 
$L _ i ^ {min} \geq  \frac{\Lambda}{n}$  was obtained through the control on the length of the cell's section due to
Proposition \ref{p:propsection} (whose proof 
 is no longer valid in the Neumann setting  
because it is based on the identification of a first eigenfunction with 
$\frac{u _ 2}{u _ 1}$). 
However, the validity of the estimate 
$L _ i ^ {min} \geq  \frac{\Lambda}{n}$ can be recovered thanks to the control on the cell's area due to Proposition \ref{p:proparea}.   
Once we have this estimate,  claim \eqref{f:sezlunga} follows, and the proof can be concluded as in the Dirichlet case. 
\qed

\bigskip \bigskip 
{\it Proof of Theorem \ref{t:quantitativeN} in dimension $N \geq 3$}. 
	 Let us  prove the inequality for $\phi \equiv 1$ (the case of a general $(\frac{1}{m})$-concave weight follows by collapsing, using Proposition \ref{lemma371}). 
		
		 It 
 is not restrictive to prove the statement under the  following assumptions: $D _\Om= \pi$, 	 
		 $\mu_1 (\Omega) \leq \pi^2+1$,
		  and  $\Om$   smooth,  so that a first eigenfunction $\overline u$ for $\mu _ 1 (\Om)$ belongs to $\mathcal C ^ 2 (\overline \Om)$.  We fix a coordinate system $(e_1, \dots, e _N)$ such that 
		  $a_1$ is aligned with $e_1$ and $a _2$ is aligned with $e _N$.  	 
	For every $n\in \N$, we consider a $L^2$ equipartition $\{ \Om _1, \dots, \Om _n \}$ of $\overline u$ in $\Om$, 
obtained  by the procedure described in Remark \ref{r:PW}, namely by using  a family of cutting hyperplanes parallel to $e _N$, 
in such way that,  for $n$ large, any  cell $\Om _i$  becomes narrow in $(N-2)$-directions, i.e., arbitrarily close to a
$2D$-convex section $U _i$. 
Via the usual argument \`a la Payne-Weinberger (namely arguing as in Lemmas \ref{l:stimaPW} and  \ref{l:preinside}, except that now our cells are narrow in $(N-2)$ in place of $(N-1)$-dimensions), we obtain 
$$\mu _ 1 (\Om _i ) \geq \mu _ 1 (U _i, \phi)  + o ( 1)\qquad \forall i = 1, \dots, n \,, $$ 
where 
$\phi$ is a $(\frac{1}{N-2})$-concave weight,  and $o ( 1)$ is an infinitesimal quantity as $n \to + \infty$. 

By the quantitative inequality  already proved for the weighted Neumann eigenvalue in $2D$,  we infer that 
$$\mu _ 1 (\Om _i ) \geq  1 + c w _{U_i} ^ 2  + o ( 1)  \qquad \forall i = 1, \dots, n \,.$$

We now search for a good proportion of cells such that $w_ {U_i}$ controls from above $a_2$. 
It is not restrictive to confine our search among cells whose diameter is sufficiently large (otherwise, if for a fixed proportion of cells the diameter is small, we easily obtain the quantitative inequality from \eqref{f:diffdbis},  by arguing as in the Dirichlet case). 
For cells with large diameter,  $w _{U _i}$ is comparable to   the width of $U _i$ in direction $e_N$, which by construction is equal to the width of $\Om _i$ in direction $e_N$, hereafter denoted by $w ^ N _{\Om _i}$.  
We claim that  there exists a dimensional constant $K$ such that 
\begin{equation}\label{f:16bis} 
w _{\Om _i}  ^N \geq K a _2    \qquad \forall i \in  \mathcal I   _n \subset \{1, \dots, n\} \text{ with } {\rm card} ( \mathcal I _n ) \geq \frac{ n}{2}\,.
\end{equation}

Indeed, for $t>0$, let us denote by $\mathcal I _ n ^ t$ the family of indices  $i \in  \{1, \dots, n \}$ such that 
$w_{\Om _i } ^N < t a _2$.   We have 
$$ |\om _i| \sim w_{\Om _i } ^N \mathcal H ^ { N-1}  (\Pi _ {e_N ^ \perp}  (\om _i )) 
< t a _ 2  \mathcal H ^ { N-1}  (\Pi _ {e_N ^ \perp}  (\om _i )) 
\qquad \forall i \in \mathcal I _n ^t \,.$$ On the other hand,  for every $i = 1, \dots, n$,  
by Proposition   \ref{p:proparea} 
  we have, for dimensional constants $\Lambda$, $\Lambda '$
$$
|\om_i | \geq \frac{\Lambda}{n} |\Om| \geq  \frac{ \Lambda '  }{n} a_1 a_2 \dots a _N \,.
$$
We infer that
$$ \frac{  \Lambda ' }{n}  a_1 a_3 \dots a _N  \leq t   \mathcal H ^ { N-1}  (\Pi _ {e_N ^ \perp}  (\om _i ))  \qquad \forall i \in \mathcal I ^ t _n\,.$$  
Summing over $i \in \mathcal I ^ t _n$, we obtain 
$$ {\rm card} (\mathcal I ^ t _n) \frac{  \Lambda ' }{n}  a_1 a_3 \dots a _N  \leq t  \mathcal H ^ { N-1}  (\Pi _ {e_N ^ \perp}  (\Om ))   \sim t a _ 1 a _ 3 \dots a _N \,,$$ 
which implies 
$$t -  {\rm card} (\mathcal I ^ t _n)  \frac{ \Lambda '}{n} \geq 0 \,.$$ 
Then claim \eqref{f:16bis}  is proved taking $\mathcal I _n := \{ 1, \dots, n\} \setminus \mathcal I _n ^ t$ , with 
$t:= \frac{\Lambda ' }{2}$.

In view of \eqref{f:16bis}, our proof is easily achieved  by applying in the usual way 
Lemma \ref{l:goodparti}: 
$$\begin{array}{ll}
\mu _ 1 (\Om )   & \displaystyle =  \frac{  \int_\Om |\nabla   \overline u  | ^ 2   }{
 \int_\Om |\overline u  | ^ 2 } \geq  \frac{1}{n} \sum _{i = 1} ^ n  \frac{  \int_{\Om_i}  |\nabla   \overline u  | ^ 2   }{
 \int_{\Om_i}  |\overline u  | ^ 2 }   \geq  \displaystyle \frac{1}{n}  \Big [ \sum _{i  \in \mathcal  I _n}  \mu _ 1 (\Om _i )  +  \sum _{i \not \in  \mathcal I _n}  \mu _ 1 (\Om _i )  
     \Big ] 
\\  \noalign{\medskip}
&   \geq  \displaystyle \frac{1}{n}  \Big [  (1 + K a _ 2 ^ 2)   {\rm card}( \mathcal I _n )   +  
   (n - {\rm card}( \mathcal I _n ) )
   \Big ] \geq 1 + \frac{K}{2} a _ 2 ^ 2 \,.
 \end{array}
$$

\qed

\section{Appendix}

  In this section we fix some results about weighted Neumann eigenvalues of the type $\mu_ 1 (\Om, p)$, defined according to \eqref{def:mu}. First we consider the case when $p$ equals $u _ 1 ^ 2$, 
and then the case when $p$ equals a $(\frac{1}{m})$-concave function $\phi$. In both cases we have that an eigenfunction exists: in the former case it can be explicitly identified with the quotient between the second and first Dirichlet eigenfunction, in the latter case it can be obtained by a collapsing procedure, namely as the limit of rescaled Neumann eigenfunctions  in higher dimensional convex sets with suitable profile.

\begin{proposition}\label{p:yau}
Given an open bounded convex domain $\Om$ in $\R ^N$, let $\lambda _ 1 (\Omega), \lambda _ 2 (\Omega)$ and  $u _1, u _ 2$ be  the first two Dirichlet eigenvalues and  eigenfunctions of $\Omega$ (normalized in $L ^2$), and let  $\mu _ 1 (\Om, u_ 1 ^ 2 )$ be defined according to  \eqref{def:mu}.  Then it holds
$$
\mu _ 1 (\Om, u_ 1 ^ 2) = \lambda _ 2 (\Omega) - \lambda _ 1 (\Omega)\,, 
$$ and  an eigenfunction for  $\mu _ 1 (\Om, u_ 1 ^ 2)$ is given by $\overline u := \frac{u _ 2}{ u _ 1}$.  
\end{proposition}

\proof  When $\Omega $ is smooth, the result is well-known  \cite[Section 1.2.2]{carron}; actually, in this case we have that $\overline u = 
\frac{u _ 2}{ u _ 1}$ is smooth up to the boundary of $\Omega$, see \cite[Appendix A]{SWYY}. 
 When $\Omega$ is not smooth, the statement can be obtained by approximation. Consider an increasing sequence of open smooth convex domains $\Om _\e \subset \Om$ 
converging to $\Om$ in Hausdorff distance as $\e \to 0$. 
For $i = 1,2$,  let $\lambda_ i ( \Om _ \e) $  and  $u _ i ^ \e$ denote the first two Dirichlet  eigenvalues and eigenfunctions of $\Om _\e$, normalized in $L ^ 2$, and  extended to $0$ in $\Om \setminus \Om _\e$.
We claim that, in the limit as $\e \to 0$, we have
\begin{equation*}
\lambda _i (\Om _\e) \to \lambda _ i (\Om)\, , \qquad u _i ^ \e \to u _ i \text{ in } H ^ 1_0 (\Om) \qquad i = 1, 2 \,, \qquad u _ 1 ^ \e \to u _ 1 \text{ in } L ^ \infty (\Om)\,.
\end{equation*} 
To prove this claim we observe first that, for $\e$ small enough,  the domains $\Om_\e$ contain a fixed ball. Hence, by the decreasing monotonicity of Dirichlet eigenvalues under domain inclusion,  the sequence $\lambda _i (\Om _\e)$ is bounded. It follows that   $\Delta u_i ^ \e$ is bounded in $L ^ 2 (\Om_\e )$, and hence  $u _ i ^ \e$ is bounded in $H ^ 2 (\Om _\e)$ (see e.g. \cite[Theorem 3.1.2.1]{gris}), so that up to subsequences it converges weakly in $H ^2 (\Om )$ and  strongly in $H ^ 1_0 (\Om)$; 
then, the limits of  
  $\lambda_ i ( \Om _ \e) $ and   $u _ i ^ \e$  can be identified respectively with $\lambda _ i (\Om)$ and $u _i$ (see e.g.\ \cite[Section 4.6]{bubu}). 
  It remains to prove that $u ^ \e _1 \to u _1$ in $L ^ \infty (\Om)$. By \cite[Lemma 3.1]{davies}, there exists a dimensional constant $C$ such that
  $ \|u _1 ^ \e \| _{\infty} \leq C \lambda _ 1 (\Om _\e) ^ { N/ 4}$,   
so that  the sequence $u _1 ^ \e$ remains bounded in $L ^ \infty (\Om_\e)$.

In turn, this implies that $\sup _ \e \|\nabla u _ 1 ^\e \| _{  L ^ \infty ( \Omega_\e) } < + \infty$.  Indeed, 
denoting by $w^\e$ the torsion function on $\Om_\e$ (i.e.  the unique solution in $H ^ 1 _0 (\Om_\e)$ 
to the equation $- \Delta w  = 1 $ in $\Om_\e$), by direct computation the function 
$$P_\e(x)= |\nabla u _ 1 ^\e|^2+ \lambda _ 1 (\Om _\e)(u _ 1 ^\e)^2-2  \lambda _ 1^2 (\Om _\e)\|u _1 ^ \e \| _{\infty} ^2 w^\e$$
is subharmonic in $\Om _\e$ (see also \cite[Section 3]{BMPV}). Hence $P_\e$ attains its maximum at the boundary. 
The uniform boundedness of $\nabla u _ 1 ^\e$ in $L ^ \infty(\Om _\e)$ follows by taking into account that 
$\| w^\e \| _{L ^ \infty ( \Omega_\e)}  $   is  bounded from above  by a constant depending on $|\Om^\e|$ and that, 
by a classical barrier argument,  $\|\nabla u _ 1 ^\e\| _{L ^ \infty (\partial \Omega_\e)} $ is bounded  from above by $\lambda_1(\Om^\e)\| u _1 ^ \e \| _{\infty} $.  

Hence the functions  $u _ 1 ^\e$ are equibounded and equicontinuous in $\Om$, and the uniform convergence of $u ^ \e _1$ to $u _1$  follows from the Ascoli-Arzel\`a theorem.

Now we observe that the function $\overline u$ is admissible in the definition \eqref{def:mu} of $\mu _ 1 (\Om)$, so that 
\begin{equation}\label{f:uno} \mu _ 1 (\Om, u _ 1 ^ 2 )\leq 
  \int _ \Om |\nabla \overline u | ^ 2  u _ 1 ^ 2 \, dx  \,.\end{equation}
By the strong convergences $u _ i ^ \e \to u _ i$  in $H ^ 1 _0(\Omega)$,  for every compact set $K \subset \Omega$, 
setting  $v ^ \e:= \frac{u _ 2 ^ \e}{u _ 1 ^ \e}$,  
  up to subsequences  it holds 
$ |\nabla  v ^ \e | ^ 2 ( u ^ \e _ 1) ^ 2  \to  |\nabla \overline  u | ^ 2  u _ 1 ^ 2$  a.e.\ on $K$.
Hence, by Fatou's lemma, 
\begin{equation}\label{f:due}  \int _ K |\nabla \overline u | ^ 2  u _ 1 ^ 2 \, dx   \leq  \liminf _ {\e \to 0 } \int _ {K} |\nabla  v ^ \e | ^ 2 ( u ^ \e _ 1) ^ 2\leq  \liminf _ {\e \to 0 } \int _ {\Om _\e} |\nabla  v ^ \e | ^ 2 ( u ^ \e _ 1) ^ 2 .
\end{equation}   
By \eqref{f:uno} and \eqref{f:due}, exploiting the arbitrariness of $K$, and  using the statement for the smooth domains $\Om _\e$, we obtain
\begin{equation}\label{f:tre} \mu _ 1 (\Om, u _ 1 ^ 2 )\leq \int _ \Om |\nabla \overline u | ^ 2  u _ 1 ^ 2 \, dx \leq  \liminf _ {\e \to 0 } (\lambda_ 2 (\Om _\e) - \lambda_ 1 (\Om _\e)  ) = \lambda_ 2 (\Om ) - \lambda_ 1 (\Om )   \,. \end{equation}

To conclude the proof, it remains to show that the two inequalities in \eqref{f:tre} are in fact equalities. 
Let $\delta >0$ be fixed, and let  $v _\delta$ be a function in $H ^ 1 _{\rm loc } (\Om)$, with  $\int_\Om v _\delta ^ 2 u _ 1 ^ 2 = 1$ and $\int _\Om v_\delta  u _ 1 ^ 2 = 0$, such that
$\mu _ 1 (\Om, u _ 1 ^ 2 )  \geq \int _\Om |\nabla v _\delta| ^ 2 u _ 1 ^ 2 - \delta$. 
Since $\Om _\e \subset \Om$, we have that $v _\delta \in H ^ 1 _{\rm loc } (\Om_\e)$. 
If the approximating domains $\Om _\e$ are suitably chosen, we have 
\begin{equation}\label{f:3convergenze}
 \begin{array}{ll}
& \displaystyle  \lim _{ \e \to 0 } \int _{\Om _\e}  v _\delta ^ 2 (u _1 ^ \e ) ^ 2 =  \int _\Om  v _\delta ^ 2 u _1   ^ 2 \,  (=1) \, , \quad
 \\ \noalign{\medskip}
 &\displaystyle\lim _{ \e \to 0 }   \int _{\Om _\e} v _\delta  (u _1 ^ \e ) ^ 2  = \int _\Om  v _\delta  u _1   ^ 2 \,  ( = 0) \, , \quad 
\\ \noalign{\medskip}
&\displaystyle \lim _{\e \to 0 } \int _{\Om _\e} |\nabla v _\delta| ^ 2 (u _1 ^ \e ) ^ 2 =  \! \int _\Om |\nabla  v _\delta| ^ 2 u _1   ^ 2\,.$$
\end{array}
\end{equation} 
More precisely, by the uniform convergence of $u _1 ^ \e$ to $u _ 1$ and  Lebesgue dominated convergence theorem, the equalities 
in \eqref{f:3convergenze}
are satisfied  provided $\Om_\e$ is chosen so that, for some $\eta_\e \ra 0$,   $u_1^\e \le (1+\eta_\e) u_1$. 
Such choice is possible thanks to the following argument. Fix the origin at the maximum point of $u_1$, and consider a small contraction $(1-\e)\Om$ of $\Om$ with respect to the origin. The  log-concavity of  $u_1$ allows to order by inclusion 
the level sets of the functions  $u_1(\frac{x}{1-\e})$ and $u_1(x)$, yielding that 
$u_1(\frac{x}{1-\e})\le u_1(x)$ for every $x \in (1-\e)\Om$, with strict inequality except at the origin. 
Then, taking $\Om_\e$ as a smooth convex approximation of  $(1-\e) \Om$, we get the required inequality $u_1^\e \le (1+\eta_\e) u_1$, for some $\eta_\e \ra 0$. 

We finally notice that the convergences in  \eqref{f:3convergenze}  remain true if we replace therein the functions $v _\d$ by their translations and normalizations 
$v _{\d, \e}$, defined by 
$$v _{\d, \e} := \frac{v _\d - \frac{1}{|\Om _\e|}  \int _{\Om _\e} v _\delta  (u _1 ^ \e ) ^ 2       } { \|v _\d - \frac{1}{|\Om _\e|}  \int _{\Om _\e} v _\delta  (u _1 ^ \e ) ^ 2  \|_{L ^ 2  (\Om _\e, (u _1 ^ \e ) ^ 2 )}} \,.  $$
Since $v _{\d, \e}$ is an admissible test function for $\mu _ 1 (\Om _\e,( u _ 1 ^ \e)  ^ 2)$,  by using the statement for the smooth domains $\Om _\e$ we obtain
$$  \mu _ 1 (\Om, u _ 1 ^ 2 )  \geq \limsup _{\e \to 0 }  \int _{\Om _\e} \!\!\! |\nabla v _{\delta, \e} | ^ 2 (u _ 1 ^ \e) ^ 2 - \delta \geq 
 \limsup _{\e \to 0 } [ \lambda_ 2 (\Om _\e) - \lambda_ 1 (\Om _\e)   - \delta]  = \lambda_ 2 (\Om ) - \lambda_ 1 (\Om )    - \delta\,. 
 $$ 
Eventually, by letting $\delta \to 0$, we conclude that  the two inequalities in \eqref{f:tre} hold with equality sign. \qed

\bigskip

 We now turn attention to the case of power-concave weights. For the convenience of the reader we start with the following 
\begin{lemma}\label{lemma371.1}
Let $\Om$ be an open bounded convex domain in $\R ^N$, and let $\phi \in L ^ 1 (\Om)$  
be a positive weight 
which is $(\frac 1 m)$-concave  for some $m \in \N \setminus \{ 0 \}$. Then the
embedding 
$H^1(\Om, \phi) \hookrightarrow L^2(\Om,\phi)$
is compact.
\end{lemma}  
\proof
Let $\widetilde{\Omega} \sq \R^{N+m}$ be defined by
	 $$
	 \widetilde{\Omega } = \left\{ (x,y) \in \R^{N}\times \R^{m} : \, x \in \Omega ,\,  \norma{y}_{\R^m} <\omega_m^{-1/m} \phi^{1/m}(x)\right\}\,. 
	 $$
	 Since $\phi^{1/m}$ is  concave,  $\widetilde{\Omega}$ is  open  and convex, so that 
	 	$H^1(\widetilde{\Omega})$ is compactly embedded into $L^2(\widetilde{\Omega})$. 
Now, if  $\{u_n\}$ is a bounded sequence in $H^1(\Om,\phi)$, setting  
	$
	\tilde{u}_n(x,y):= u_n(x)$
	 we have
	 \begin{align*}
	\int_{\widetilde{\Omega}} \tilde{u}_n^2(x,y) \, dxdy &= \int_{\Omega} u_n^2(x)  \phi (x) \, dx ,\\
		\int_{\widetilde{\Omega}} \abs{\nabla_{x,y}\tilde{u}_n}^2(x,y) \, dxdy &= \int_{\Omega} \abs{\nabla_x u_n}^2(x)  \phi(x) \, dx.
	 \end{align*}
	Then, up to subsequences, $\tilde u _n$ converges weakly in 
	$H^{1}(\widetilde{\Omega}) $ and strongly in $L ^ 2 ( \widetilde \Om)$ to a function $\tilde u \in H^{1}(\widetilde{\Omega})$.  
	 Since $\tilde u$ is constant in the  $y$ variable, the function $u(x): = \tilde u (x, y)$ belongs to $H^1(\Omega,\phi)$, and $u _n$ converges strongly to $u$ in $L ^ 2 (\Om, \phi)$. 
\qed

\bigskip
  Now, under the assumptions of Lemma \ref{lemma371.1}, 
the compactness of the embedding $H^1(\Om, \phi) \hookrightarrow L^2(\Om,\phi)$ ensures that    
the operator mapping a function $f \in L^2(\Om,\phi)$ with  $\int_\Om f \phi =0$ into the unique solution to 
$$
 \inf_{v \in H^1(\Omega,\phi)} \left(\frac{1}{2}\int_{\Omega}  \abs{\nabla v}^2 \phi  - \int_{\Omega}  f v \phi \right)
$$
is positive, self-adjoint, and compact. Then the eigenvalues of the weighted problem 
\begin{equation*}
	\begin{cases}
		- {\rm div} \, (\phi \nabla  u)= \mu (\Om, \phi) \phi u & \text{in } \, \Omega \\
		\displaystyle{ \phi \frac{\partial u}{\partial \nu}=0} & \text{on} \, \partial\Omega
	\end{cases}
\end{equation*}
can be computed by the classical min-max formula. In particular,  we have  that
$$\mu_1(\Om,\phi)= \inf_{u \in H^1(\Om,\phi)\sm\{0\}, \int_\Om u \phi =0} \frac{\int_\Om |\nabla u|^2 \phi }{\int_\Om u^2 \phi },$$
and the infimum is attained.
We now show that the above eigenvalue  and a first eigenfunction associated with it can be obtained by collapsing, i.e.\ by a limiting procedure as as $\e \to 0^+$ starting   from the 
convex hypographs 
	\begin{equation}\label{lemma371.4}
	 \widetilde{\Omega }_\vps := \left\{ (x,y) \in \R^{N}\times \R^{m} : \, x \in \Omega ,\,  \norma{y}_{\R^m} <\vps \omega_m^{-1/m} \phi^{1/m}(x)\right\}\sq \R^{N+m}.
	\end{equation}

	\begin{proposition}\label{lemma371}
  Let $\Om$ be an open bounded convex domain in $\R ^N$, and let $\phi \in L ^ 1 (\Om)$  
be a positive weight 
which is $(\frac 1 m)$-concave  for some $m \in \N \setminus \{ 0 \}$. 
Let $\tilde u_\vps$ be $L^2$-normalized first eigenfunctions associated with $\mu_1(\widetilde{\Omega }_\vps)$. Setting 
$u_\vps(X,Y):= \vps ^{\frac m2} \tilde u_\vps (X, {\vps} Y)$,  up to subsequences we have 
$$\mu_1( \widetilde{\Omega }_\vps) \to \mu_1(\Om, \phi)\,\; \mbox{ and } \;\;    u_\vps  \to \tilde u \mbox{ in } L^2(\widetilde{\Omega }_1),$$
 where $\tilde u(X,Y):=u(X)$,
 $u$ being a $L^2(\Om, \phi)$-normalized eigenfunction associated with $\mu _ 1 (\Om , \phi)$. 
 In particular
\begin{equation}\label{lemma371.3} \|u\|_{L^\infty (\Om)} \le \liminf_{\vps \to 0} \vps^{\frac m2}\|\tilde u_\vps\|_{L^\infty (\widetilde{\Omega }_\vps)} .
\end{equation}
\end{proposition}

\proof 
By the change of variables $X=x, Y= \frac{1}{\vps} y$, we get
$$
\begin{array}{ll} 
& \displaystyle 
 \int_{\widetilde{\Omega }_1} |\nabla_X u_\vps |^2+ \frac{1}{\vps ^2} |\nabla_Y u_\vps|^2   dX dY  =\int_{\widetilde{\Omega }_\vps} |\nabla_x \tilde u_\vps |^2 + |\nabla_y \tilde u_\vps |^2 dx dy= \mu_1( \widetilde{\Omega }_\vps ) \\  \noalign{\medskip}
& \displaystyle \int_{\widetilde{\Omega }_1} u_\vps^2(X,Y) dXdY =1, \; \; \int_{\widetilde{\Omega }_1} u_\vps(X,Y) dXdY=0\,.
\end{array}
$$ 
From Kr\"oger's inequality \cite{kro}, we know that $\limsup \mu_1( \widetilde{\Omega }_\vps) <+\infty$, so that $\{u_\vps\}$ is bounded in $H^1(\widetilde{\Omega }_1)$, and possibly passing to a subsequence it converges weakly in $H^1(\widetilde{\Omega }_1)$ to some function $\tilde u$ with $\nabla _Y \tilde u =0$ in $\widetilde{\Omega }_1$. Setting $u(X) := \tilde u(X,Y)$, we get that 
$\int_\Om u^2 \phi =1$, $\int_\Om u\phi=0$ and
$$ \mu_1(\Om, \phi) \le \int_\Om |\nabla u|^2\phi \le \liminf _{\vps \to 0} \mu_1(\widetilde{\Omega }_\vps).$$
To show the converse inequality, let $v$ be a normalized first eigenfunction for  $\mu _ 1(\Om, \phi)$ and define
$$\tilde v_\vps  \in H^1( \widetilde{\Omega }_\vps), \qquad \tilde v_\vps(x,y):= \vps^{-\frac m2} v(x).$$
We have$$\int _{\widetilde{\Omega }_\vps}   \tilde v_\vps =0, \qquad \int _{\widetilde{\Omega }_\vps}   \tilde v_\vps^2 =1, \qquad \int_{\Om}|\nabla v |^2 \phi = \int _{\widetilde{\Omega }_\vps}|\nabla \tilde v_\vps|^2,   $$
so  that $\tilde v_\vps $ is a test function for $\mu_1(\widetilde{\Omega }_\vps)$ and 
$\mu_1(\Om, \phi) \ge \limsup _{\vps \ra 0} \mu_1(\widetilde{\Omega }_\vps)$. 
We conclude that  the equality 
$\mu_1(\Om, \phi) = \lim_{\vps \ra 0} \mu_1(\widetilde{\Omega }_\vps)$ holds, 
and that the function $u$ above has to be an eigenfunction for $\mu_1(\Om, \phi) $. Together with the convergence $u_\vps \to \tilde u$ in $L^2( \widetilde{\Omega }_1)$, this implies \eqref{lemma371.3}.
 \qed 
 
\bigskip 
 
\noindent

{\bf Acknowledgments.} The authors wish to thank M.\ van den Berg  for fruitful discussions and suggestions and  B. Bogo\c{s}el  for numerical evidence related to Theorem \ref{t:1d}. 

\noindent

{\bf Statements and Declarations.}  The authors have no competing interests.

\bibliographystyle{mybst}

\bibliography{References}
\end{document}